\documentclass[a4paper,11pt,review]{article}
\pdfoutput=1

\title{Error boundedness of Correction Procedure via Reconstruction / Flux Reconstruction}

\author{Philipp \"Offner }
\date{\today}

\usepackage[utf8]{luainputenc}
\usepackage[english]{babel}
\usepackage{csquotes}

\usepackage[a4paper, textheight=622pt, textwidth=468pt, centering]{geometry}

\usepackage[plainpages=false,pdfpagelabels,hidelinks]{hyperref}
\makeatletter
\hypersetup{pdfauthor={\@author}}
\hypersetup{pdftitle={\@title}}
\makeatother

\usepackage{cite}

\usepackage{amsmath}
\allowdisplaybreaks
\usepackage{amssymb}
\usepackage{commath}
\usepackage{mathtools}

\usepackage{siunitx}

\usepackage{amsthm}
\theoremstyle{plain}
  \newtheorem{thm}{Theorem}[section]
  
  \theoremstyle{definition}
  \newtheorem{ex}[thm]{Example}
\newtheorem{re}[thm]{Remark}
\newtheorem{as}[thm]{Asumption}
\newtheorem{Result}[thm]{Result}
  
\usepackage{color}
\usepackage{graphicx}
\usepackage[small]{caption}
\usepackage{subcaption}

\usepackage{pgfplots}
\pgfplotsset{compat=1.11}
\usepgfplotslibrary{external} 
\begingroup\expandafter\expandafter\expandafter\endgroup
\expandafter\ifx\csname pdfsuppresswarningpagegroup\endcsname\relax
\else
\fi

\usepackage{booktabs}
\usepackage{rotating}
\usepackage{multirow}

\usepackage{calc}
\usepackage{xparse}

\renewcommand{\vec}[1]{\underline{#1}}
\NewDocumentCommand{\mat}{mo}{%
  \IfValueTF{#2}{%
    \underline{\underline{#1}}{#2}
  }{%
    \underline{\underline{#1}}\,
  }%
}

\newcommand{\diag}[1]{\operatorname{diag}\left(#1\right)}

\renewcommand{\d}{\operatorname{d}}

\renewcommand{\L}{\mathbf{L}}
\renewcommand{\P}{\mathbb{P}}

\newcommand{\Ep}{\mathbb{E}}
\newcommand{\Eps}{\mathbf{E}}

\newcommand{\fnum}{f^{\mathrm{num}}}
\newcommand{\vecfnum}{\vec{f}^{\mathrm{num}}}
\newcommand{\vecfnumk}{\vec{f}^{\mathrm{num},k}}

\newcommand{\K}{\tilde{\operatorname{K}}}

\newcommand{\vecepsilon}{\vec{\epsilon}_1^{\mathrm{num},k}}

\newcommand{\Ip}{\mathbb{I}}

\renewcommand{\epsilon}{\varepsilon}
\renewcommand{\phi}{\varphi}
\newcommand{\N}{\mathbb{N}}

\newcommand{\Ol}{\mathcal{O}}

\renewcommand{\r}{\right}
\renewcommand{\l}{\left}

\newcommand{\est}[1]{\left\langle#1\right\rangle}

\newsavebox{\DelimiterBox}
\newlength{\DelimiterHeight}
\newlength{\DelimiterDepth}
\newsavebox{\ArgumentBox}
\newlength{\ArgumentHeight}
\newlength{\ArgumentDepth}
\newlength{\ResizedDelimiterHeight}

\newcommand{\mean}[1]{\overline{#1}}

\newcommand{\jump}[1]{%
  \savebox{\ArgumentBox}{$ \displaystyle #1 $}%
  \settoheight{\ArgumentHeight}{\usebox{\ArgumentBox}}%
  \settodepth{\ArgumentDepth}{\usebox{\ArgumentBox}}%
  \savebox{\DelimiterBox}{$ [\!\![ $}%
  \settoheight{\DelimiterHeight}{\usebox{\DelimiterBox}}%
  \settodepth{\DelimiterDepth}{\usebox{\DelimiterBox}}%
  \setlength{\ResizedDelimiterHeight}{\maxof{1.2\ArgumentHeight}{\DelimiterHeight}}
  \!
  \resizebox{\width}{\ResizedDelimiterHeight}{ [\![ }
  \mkern-6.5mu
  #1
  \mkern-6.5mu
  \resizebox{\width}{\ResizedDelimiterHeight}{ ]\!] }
  \!\!
}

\usepackage{lineno}
\newcommand*\patchAmsMathEnvironmentForLineno[1]{%
  \expandafter\let\csname old#1\expandafter\endcsname\csname #1\endcsname
  \expandafter\let\csname oldend#1\expandafter\endcsname\csname end#1\endcsname
  \renewenvironment{#1}%
     {\linenomath\csname old#1\endcsname}%
     {\csname oldend#1\endcsname\endlinenomath}}%
\newcommand*\patchBothAmsMathEnvironmentsForLineno[1]{%
  \patchAmsMathEnvironmentForLineno{#1}%
  \patchAmsMathEnvironmentForLineno{#1*}}%
\AtBeginDocument{%
\patchBothAmsMathEnvironmentsForLineno{equation}%
\patchBothAmsMathEnvironmentsForLineno{align}%
\patchBothAmsMathEnvironmentsForLineno{flalign}%
\patchBothAmsMathEnvironmentsForLineno{alignat}%
\patchBothAmsMathEnvironmentsForLineno{gather}%
\patchBothAmsMathEnvironmentsForLineno{multline}%
}

\definecolor{darkspringgreen}{rgb}{0., 0.55, 0.3}
\definecolor{dartmouthgreen}{rgb}{0.05, 0.5, 0.06}
\definecolor{etonblue}{rgb}{0.59, 0.78, 0.64}
\definecolor{airforceblue}{rgb}{0., 0.4, 0.66}
\definecolor{arylideyellow}{rgb}{0.91, 0.84, 0.42}
\definecolor{emerald}{rgb}{0.31, 0.78, 0.47}
\definecolor{uclagold}{rgb}{1.0, 0.7, 0.0}
\definecolor{cadmiumorange}{rgb}{0.93, 0.53, 0.18}

\begin{document}

\maketitle

\begin{abstract}

We study the long-time error behavior of correction procedure via 
reconstruction / flux reconstruction (CPR/FR) methods for  linear hyperbolic conservation 
laws. We
 show that not only the choice of the numerical flux (upwind or central)
affects the growth rate and asymptotic value of the error, but  that the 
selection of bases (Gauß-Lobatto or Gauß-Legendre) is even more important.
Using a Gauß-Legendre basis, the error reaches the asymptotic value faster
and to a lower value than when using a Gauß-Lobatto basis. 
Also the differences in the error caused by the 
numerical flux are  not essential for low resolution computations
in the Gauß-Legendre case. This behavior is better seen on  a particular FR scheme 
which has a strong connection with the discontinuous Galerkin framework, but holds also 
for other flux reconstruction schemes with low order resolution computations.

\end{abstract}

\section{Introduction}\label{sec1:Intro}

There exists plenty of examples in the literature where stable approximations
of hyperbolic conservation laws exhibit
a linear error growth (or nearly linear growth) in time, 
even though stability of the numerical schemes should guarantee that the solution 
remains bounded, see
\cite{hesthaven2002nodal, cohen2006spatial, nordstrom2003high}.  
The reason behind this is the following:
the error equation for the time variation contains a forcing term generated by the
approximation or truncation errors and
this forcing term can trigger the unbounded growth of the error.\\
Simultaneously, there are other examples where the temporal error growth is
bounded \cite{abarbanel2000error, koley2009higher}.
Finally, in \cite{nordstrom2007error}, the author gives an explanation 
under what conditions the error is or is not bounded in time. 
The author works with SBP-SAT (Summation-by-Parts/Simultaneous-Approximation-Term) 
finite difference approximation and deduces
that the error behavior depends only on the choice of boundary condition of the problem. \\
If one considers waves in cavities or with periodic boundary conditions,
 linear growth is observed as it is investigated in \cite{hesthaven2002nodal},
where for inflow-outflow problems one obtains  boundedness.
In other words, if an appropriate boundary condition (sufficiently dissipative) 
is applied, the error is bounded. In this framework, the error behavior
does not depend on the internal 
discretization. 
In \cite{kopriva2017error}, the authors analyze the long-time behavior of the error
for discontinuous Galerkin spectral element methods (DGSEM).
They confirm the conclusion from \cite{nordstrom2007error} that the bounded error 
property is due to the dissipative boundary conditions, but in contrast
to \cite{nordstrom2007error}, in the DGSEM framework the internal approximation 
has indeed an influence on the behavior of the error. 
The choice of the numerical flux (upwind or  central) is essential for 
the magnitude of the error and the  speed at which the asymptotic error is 
reached. With the upwind flux one obtains better results.  \\
In this paper, we examine the long-time error behavior for a recent 
correction procedure via reconstruction (CPR) / flux reconstruction (FR) method.
The CPR/FR is a unifying  framework for several high-order methods such as discontinuous
Galerkin (DG), spectral difference (SD), spectral volume (SV), and
the residual distribution methods \cite{abgrall2018connection, huynh2014high},
and  includes through special choice of the nodal basis and of the correction matrix 
the DGSEM of \cite{kopriva2017error}. 
Here, we investigate not only the numerical flux, but also the selection of nodal 
basis (Gauß-Legendre and Gauß-Lobatto). We recognize that  the selection 
of the flux function is less  important than the choice of the nodal basis
for the error behavior. 
Using Gauß-Legendre basis in the approximation, the error is lower than in the Gauß-Lobatto
case  and the speed of attaining the error asymptotics  is even faster 
for most of the problems under consideration. 
The selection of the numerical flux 
has less influence  on the error behavior when using Gauß-Legendre nodes
than when applying a Gauß-Lobatto 
basis. Our investigation leads us to the conclusion that for many problems 
Gauß-Legendre points are the right choice, especially for low-order resolutions.
Our analysis extends and completes the investigation / predictions from 
\cite{kopriva2017error}
to a more general framework. \\ 
The paper is organized as follows:
in the second section, we repeat the main ideas of the SBP-CPR/FR methods 
and demonstrate the connection 
between CPR/FR and the DG framework. Then, in the section \ref{sec:model}
we present the model problem under consideration. In the next section
 \ref{sec:stability} we provide some approximation results and repeat the stability 
 analysis of the SBP-CPR/FR methods. 
Then, we extend the results from \cite{kopriva2017error}
to the linearly stable one-parameter family of Vincent \cite{vincent2011newclass}
and also consider in our investigation Gauß-Legendre nodes in the section \ref{sec:Error_equation}.
These nodes do not contain the boundary values 
in one element and this yields  a further error term in our error equation.
We focus on this additional error term and give an interpretation for it. 
We confirm our theoretical investigation by numerical tests in the section \ref{6_Numerical} which includes also one
example from \cite{kopriva2017error}
for comparison.
We mention some limitations of our results and finally, we summarize and
discuss these limitations. 
In the appendix, we show the relation between FR and DG and focus on stability conditions for FR methods
as described in \cite{vincent2011newclass}.

\section{Correction Procedure via Reconstruction/Flux Reconstruction using Summation-by-Parts Operators}
\label{sec2:CPR}

In the first part, we shortly repeat the main idea of 
CPR/FR methods using Summation-by-parts Operators (SBP).
For the rest of this work, we call them FR methods. 
We follow the introduction and notation given in the articles 
\cite{ranocha2016summation,ranocha2017extended}.\\
We consider a one-dimensional  scalar conservation law
{\small 
\begin{equation}
\label{eq:scalar_CL} 
  \partial_t u (x,t) + \partial_x f(u(x,t)) = 0, \qquad t>0,\; x\in (0,L)
\end{equation}}
equipped with  adequate initial and boundary conditions.
The domain $(0,L)$ is split into $K$ non-overlapping elements $[0,L]=[x^0,x^1]
\bigcup \cdots\bigcup
[x^{K-1},x^K]$. 
The FR method is a semidiscretization applying a polynomial approximation on elements. 
Each interval $[x^{k-1},x^{k}]$ is transferred  onto a standard element.
In our case we consider $[-1,1]$ and  all calculations are conducted within this reference element.
The term $ \frac{\Delta x_k}{2}=\frac{x^k-x^{k-1}}{2} $ denotes 
the transformation factor.
Let $\P^N$ be the space of polynomials of degree
$ \leq N$, $-1 \leq \xi_i \leq 1$ $(0\leq i\leq N$) be the
 interpolation points in $[-1,1]$, $\Ip^N:C([-1,1])\to \P^N[-1,1]$  
the interpolation operator and $P_{N-1}^m u$  the orthogonal projection of 
$u$ onto $\P^{N-1}$ with respect to the inner product of the Sobolev space $H^m((-1,1))$.
The solution $u$ is approximated by a polynomial $U \in \P^N$. 
A nodal Lagrange 
basis is usually employed\footnote{Modal bases are also possible \cite{ranocha2017extended},
 but we won't consider these in
 this paper.}.
 Instead of working with $U$ one can also express the numerical solution 
 as the vector $\vec{u}$ with coefficients 
$\vec{u}_i = U(\xi_i), i \in \set{0, \dots, N}$.
All the relevant information are stored in these coefficients and one may write
{\small
\begin{equation}\label{eq:Approx}
 u(\xi)\approx U(\xi) =\sum\limits_{i=0}^N \vec{u}_i l_i(\xi),
\end{equation}
}
where $l_i(\xi)$ is the i-th Lagrange interpolation polynomial that satisfies 
$l_i(\xi_j)=\delta_{ij}$.
In finite difference (FD) schemes, it is common to work with the coefficients
only and since we are working with SBP operators with origins lying in the FD
community \cite{kreiss1974finite}, we  utilize the coefficients as well.
The flux $f(u)$ is also approximated by 
a polynomial, where the coefficients are given by
$\vec{f}_i = f \left( \vec{u}_i \right) = f \left( U(\xi_i) \right)$.

With respect to the selected basis (interpolation points), an approximation 
of the derivative is represented by the matrix $\mat{D}$.
Moreover, a discrete scalar product is represented by the symmetric and positive
mass/norm matrix $\mat{M}$. This matrix approximates the usual $\L^2$ scalar product.
It is 
{\small 
\begin{equation}\label{eq:approx_matrix}
 \mat{D}\vec{u}\approx \vec{\partial_x u} \text{ and }
(\vec{u},\vec{v})_N:= \vec{u}^T\mat{M}\vec{v} \approx \int_{x^{k-1}}^{x^{k}} u v\d x.
\end{equation}}
Applying Lagrange polynomials, we obtain $D_{ij}=l_j'(\xi_i)$.
The matrix $\mat{M}=\diag{\omega_0,\cdots,\omega_N}$ is associated as usual with the
quadrature rule 
given by the polynomial basis (Gauß-Lobatto or Gauß-Legendre) where
$\omega_j$ are the quadrature weights associated with the nodes $\xi_j$.
For Gauß-Legendre nodes, $\omega_j=\int_{-1}^1 l_j(x)\d x$. 
Note that in case of Gauß-Lobatto nodes, the mass matrix is in general not exact. 
As described in the review articles \cite{svard2014review, fernandez2014review}
SBP operators are constructed in such way that they mimic integration-by-parts on a discrete
level. Up to now, we have expressions for the derivative as well as for the integration.
Hence, the evaluation on the boundary is missing. Here, we have to introduce 
two different operators. 
First, the restriction operator denoted by the matrix $\mat{R}$ which 
approximates the interpolation of a function to the boundary points 
$\{x^{k-1},x^{k}\}$
Second, the diagonal boundary matrix $\mat{B}=\diag{-1,1}$ that  gives the difference of boundary
values. This means 
{\small 
\begin{equation*}
 \mat{R} \vec{u}\approx \begin{pmatrix}
                         u(x^{k-1})\\
                         u(x^{k})
                        \end{pmatrix} \text{ and } (u_L,u_R) \mat{B}\begin{pmatrix}
                         v_L\\
                         v_R
                        \end{pmatrix}=u_Rv_R-u_Lv_L.
\end{equation*}
}
where $v_{i}$ $(i=L,R)$ describes the position in the element, i.e. the left and right boundary points.
Finally, all operators are introduced and they
have to fulfill the SBP property
{\small 
\begin{equation}
\label{eq:SBP}
  \mat{M} \mat{D} + \mat{D}[^T] \mat{M}
  = \mat{R}[^T] \mat{B} \mat{R},
\end{equation}
}
in order to mimic integration-by-parts on a discrete level
{\small 
\begin{equation}
  \vec{u}^T \mat{M} \mat{D} \vec{v} + \vec{u}^T \mat{D}[^T] \mat{M} \vec{v}
  \approx
  \int_{x^{k-1}}^{x^k} u \, (\partial_x v) \d x+ \int_{x^{k-1}}^{x^k} (\partial_x u) \, v\d x
  = u \, v \big|_{x^{k-1}}^{x^k}
  \approx
  \vec{u}^T \mat{R}[^T] \mat{B} \mat{R} \vec{v}.
\end{equation}
}
As an example, we consider Gauss-Lobatto nodes in $[-1,1,]$ which include the boundary points.
Then, the restriction operators are simply 
{\small 
\begin{equation}
 \mat{R}=\begin{pmatrix}
          1 &0 &\cdots &0&0\\
          0&0&\cdots &0&1 
         \end{pmatrix}, \qquad 
    \mat{R}[^T] \mat{B} \mat{R}=\diag{-1,0,\cdots,0,1}.
\end{equation}
}
The general aspects of SBP operators are introduced 
and we focus on our FR approach now.
Contrary to DG methods, we do not use a variational formulation 
(i.e. weak form) of \eqref{eq:scalar_CL}.
Instead, the differential form is applied, corresponding to a strong form DG method.
To describe the semidiscretisation all operators
are introduced. We apply the
 discrete derivative matrix $\mat{D}$ to $\vec{f}$.
The divergence is $\mat{D} \vec{f}$.
 Since the numerical solutions will probably have discontinuities
across elements, we will have this in the discrete flux, too.
In order to avoid this problem, a numerical flux $\vecfnum$ is introduced
which computes a common flux at the boundary using values from both neighboring elements.
The main idea of the FR schemes is that the numerical flux at the boundaries
will be corrected by
functions in such manner that information of two neighboring elements
interact and basic properties, like conservation, hold also in the  semidiscretisation.
Therefore, we add a correction term using a correction matrix $\mat{C}$ at the boundary nodes.
This gives \emph{Flux Reconstruction} its name.
Hence,  a simple FR method for
{\small 
\eqref{eq:scalar_CL} reads
\begin{equation}
\label{eq:SBP CPR}
  \partial_t \vec{u}
  = - \mat{D} \vec{f}
    - \mat{C}\left( \vecfnum - \mat{R} \vec{f} \right).
\end{equation}
}
A general choice of the correction matrix $\mat{C}$ recovers the linearly stable
flux reconstruction methods of \cite{vincent2011newclass, vincent2015extended},
as presented in \cite{ranocha2016summation}.\\
In our investigation,
we only consider the one parameter\footnote{The results for the multi-parameter family 
are similar to those about the one parameter family, 
since the one parameter family is contained in the extended range of schemes. 
the one parameter family for simplicity.} family of Vincent et al. \cite{vincent2011newclass}.
To describe the setting and to specify the correction matrix,
we introduce a symmetric matrix $\mat{\K}$ 
satisfying $\mat{M}+\mat{\K}>0$, i.e. positive definite.
Then, the correction matrix $\mat{C}=(\mat{M}+\mat{\K})^{-1} \mat{R}[^T] \mat{B}$
is applied in \eqref{eq:SBP CPR}
where $\K$ is defined through:
{\small 
\begin{equation}\label{eq:one_parameter}
  \mat{\K}=\kappa(\mat{D}^N)^T \mat{M}\mat{D}^N, \quad \text{ with }  \mat{\K} \mat{D}=0.
\end{equation}
}
The term $\kappa$ represents the free parameter and  
the selection of  $\kappa$ yields different numerical methods.
In particular, $\kappa\equiv0$ is the canonical choice of the correction matrix
and the resulting scheme corresponds to a strong form of a DG method \cite{gassner2013skew}.
Furthermore, since \eqref{eq:one_parameter} holds ($\mat{D}^{N+1}=0$ 
(polynomials of degree $\leq N$)) and $\K$ is symmetric, we may write
 {\small
\begin{equation}
\label{eq:SBP_2}
  \l(\mat{M}+\mat{\K} \r)\mat{D} + \mat{D}[^T] (\mat{M}+\mat{\K})
  \stackrel{\eqref{eq:one_parameter}}{=} \mat{M} \mat{D} + \mat{D}[^T] \mat{M} 
  \stackrel{\eqref{eq:SBP}}{=} \mat{R}[^T] \mat{B} \mat{R}.
\end{equation}}
Therefore, the SBP property is also valid for $ \mat{M}+\mat{\K}$.

\begin{re}\label{kappa}
 The only condition on $\kappa$ is given by the requirement that $\mat{M}+\mat{K}$ is positive definite.
 It is essential since the term represents a norm in the discrete setting and linear stability will be 
 analyzed in respect to this discrete norm in subsection \ref{subsec:stability_inv}.
 To guarantee that the term $\mat{M}+\mat{K}$  with \eqref{eq:one_parameter} is positive definite,
 $\kappa$ has to be bounded from below. In \cite[Section 3.6]{ranocha2016summation}, the bounds on 
 $\kappa$ are determined and we repeat them here for completeness. The superscript denotes the used
 nodes (Gauß-Legendre with  $G$, Gauß-Lobatto with $L$).
 With $a_N=\frac{(2N)!}{2^N(N!)^2}$, we get the following bounds:
 {\small 
 \begin{equation*}
  \kappa>\kappa_-^G:=-\frac{1}{(2N+1)a_N^2(N!)^2}, \qquad \kappa>\kappa_-^L:=-\frac{1}{Na_N^2(N!)^2}\\
 \end{equation*}
 }
 Furthermore, we like to mention that the investigation is based on an idea of Jameson \cite{jameson2010proof}. Instead 
 of working with the classical $\L^2$ norm, he applies a broken Sobolev norm involving derivatives
 and the argument that in finite-dimensional vector spaces all norms are equivalent. 
In the appendix \ref{sec:appendix}, we give the definition of the used norm together with more explanations
and an example about the connection between DG methods and the FR framework.  
\end{re}
As we already mentioned before, the different selection of $\kappa$ yields various numerical methods.
\begin{table}[ht]
{\footnotesize{
\begin{center}
               \begin{tabular}{l|l||c||r} \hline
$N$ & $\kappa _{SD}$ & $\kappa _{Hu}$ & $\kappa _{DG}$\\  \hline
$2$ & $ 4/135$ & $ 1/15$ &$0$\\
$3$ & $ 1/1050$ & $ 8/4725 $ &$0$ \\ 
$4$ & $ 8/496125$& $1/39690$ & $0$\\
$5$ & $ 1/5893965$& $12/49116375$ & $0$\\\hline
\end{tabular}\caption{Values of $\kappa $ to get different numerical schemes \cite{vincent2011newclass}}\label{ta:correction_terms}
\end{center}
  }}            \end{table}
In table \ref{ta:correction_terms},
we provide the terms $\kappa$ for, in our opinion, the most popular 
FR schemes for different order of accuracy. The exact formulas can be found in the appendix.\\
Theoretically, the parameter $\kappa$ can 
tend to infinity as it is described and analyzed  in the same paper.
However, the numerical results in \cite{vincent2011newclass, ranocha2016summation}
show that the most accurate results are obtained when $\kappa=0$ is used
and significant accuracy is lost for 
$\kappa\to \infty$. 
Hence, we restrict ourself in the investigation and consider schemes
between the range of DG, spectral difference and Huynh scheme \cite{huynh2007flux}.
We assume the following:
\begin{as}\label{As:c_zero}
With increased $N$, 
$\kappa$ tends to zero. It is $\kappa \in \Ol\left(N^{-1}(a_NN!)^{-2}\right)$. 
\end{as}

%
%

\section{Model Problem}
\label{sec:model}
To analyze the long time error behavior of SBP-FR method, we study (similarly to 
\cite{nordstrom2007error, kopriva2017error})
the scalar linear advection equation with 
non-periodic boundary conditions 
{\small
\begin{equation}\label{eq:Model}
\begin{aligned}
 \partial_t u+\partial_x u&=0, \quad x\in [0,L], \quad t\geq 0 \\ 
 u(t,0)&=g(t),\\
 u(0,x)&=u_0(x). 
\end{aligned}
\end{equation}}
We assume also that the initial and boundary values are chosen in such way that 
 $u(t,\cdot)\in H^m_{\kappa,N}((0,L))$ for
$m>1$ and that $||u(t,\cdot)||_{H^m_{\kappa,N}}$ is uniformly bounded in time. 
As it is described in \cite{kopriva2017error},
such conditions  are physically meaningful, because they describe problems
where the boundary input is, for instance, 
sinusoidal. In our numerical tests in section \ref{6_Numerical},
we will present  an example where these conditions are not fulfilled, see 
subsection \ref{subsec:counter}. 
Here, $H_{\kappa,N}^m$ denotes the following function space
{\small
 \begin{equation*}
  H^m_{\kappa,N}((0,L))=\{u\in \L^2((0,L))\;|\; \forall |j| \leq m\;:  u^{(j+N)} \in \L^2((0,L))  \}
 \end{equation*}}
 equipped with the norm 
 {\small
  \begin{equation}\label{Sobolev}
   ||u||_{H_{\kappa,N}^m}:= \left(\sum\limits_{j=0}^m ||u^{(j)}||_{\L^2((0,L))}^2
  +\kappa||u^{(j+N)}||_{\L^2((0,L))}^2 \right)^\frac{1}{2}.
  \end{equation}} 
  In
total, the highest derivative is $m+N$. If $N$ tends to infinity,
$\kappa$ tends to zero rapidly because of our assumption \ref{As:c_zero}.
This means nothing else that in \eqref{Sobolev}, the first term is always 
dominant in the norm calculation\footnote{Instead of working 
with $H^m_{\kappa,N}$ we also may work directly
$H^{m+N}$ in our investigation since we always assume sufficient smoothness of the solution.}.
We will use this fact later in the estimation of the error.
As it is  well-known, the boundary conditions of \eqref{eq:Model} 
have an essential impact on the solution and in \cite{nordstrom2017roadmap},
the author also shows that a correct implementation of the boundary conditions is essential for well-posedness.
We shortly give the following example from \cite{kopriva2017error}  
in the DG context whereas for general FR schemes the analysis  
can be found in \cite{vincent2011newclass, vincent2015extended}.
\begin{ex}\label{example_DG_boundary}
Here, the energy of the solution $u$ of the initial boundary value problem \eqref{eq:Model} is measured by the 
standard $\L^2$-norm 
$||u||^2= \int_{0}^L u^2\d x$. Focusing on the weak formulation of the advection equation \eqref{eq:Model}, 
we multiply with a test function $\phi \in C^1([0,L])$ and  integrate over the domain. We get
 {\small
\begin{equation*}
 \int_0^L u_t\phi \d x +\int_{0}^L u_x \phi \d x =0.
\end{equation*}}
Taking $\phi=u$ and integration by parts yield 
 {\small
\begin{equation*}
 \frac{1}{2} \frac{\d}{\d t} ||u||^2=-\frac{1}{2}\int_0^L u u_x \d x= \frac{1}{2}\l(g^2(t)-u^2(L,t)\r).
\end{equation*}}
Integration in time over an interval $[0,T]$ leads to
 {\small
\begin{equation}\label{eq:stability}
 ||u(T)||^2+\int_0^T u^2(L,t) \d t =||u_0||^2+ \int_0^T g^2(t) \d t.
\end{equation}}
We see that the energy at time $T$ can be expressed by the initial energy plus the energy added at the left side through the boundary
condition minus the energy, which we lose through the right side. 
Therefore, the selection of the boundary conditions is essential and the numerical approximation has to imitate this.
\end{ex}

\section{Approximation Results and Stability of the SBP-FR Methods}\label{sec:stability}
\subsection{Numerical Errors and Approximation Results}
Before we start with our stability analysis of the SBP-FR methods 
and derive the error equations in the next section \ref{sec:Error_equation},
we give an overview of the notation
and some basic approximation properties which will be used later in this paper. 
In table \ref{tab:notation}, we summarize the applied notation for clarification.\\
We analyze stability in the semidiscrete sense.
Therefore, we  divide  the entire interval $[0,L]$ into elements
$e^k=[x^{k-1},x^{k}]$, $k=1,\dots, K$,
where the $x^k$ are the element boundaries, and in particular where $x^0=0$ and $x^K=L$. 
As it was explained in section
\ref{sec2:CPR},
we transform every element to our standard element and use a SBP-FR method. 
 We investigate both Gauß-Lobatto and Gauß-Legendre quadrature.
We can define the discrete inner product by
  {\small
 \begin{equation}\label{eq:Quadrature_Lobatto}
 (U,V)_N:=\sum_{j=0}^N U(\xi_j)V(\xi_j) \omega_j.
 \end{equation}
For  Gauß-Lobatto $(UV \in \P^{2N-1})$ /
Gauß-Legendre $(UV \in \P^{2N+1})$ quadrature, it holds 
\small{
\begin{equation*}
  (U,V)_N=\sum_{j=0}^N U(\xi_j)V(\xi_j) \omega_j=\int_{-1}^1 UV \d \xi
 \quad \forall UV \in \P^{2N-1}/\; \P^{2N+1}.
\end{equation*}
  }}
We choose the numerical flux to have the form 
 {\small
\begin{equation*}
 \fnum(U_L,U_R)=\frac{U_L+U_R}{2}-\frac{\sigma}{2}(U_R-U_L), \quad \sigma \in [0,1],
\end{equation*}}
where $U_L, U_R$ are the states on the left and right.
For $\sigma=0$ we get the central flux and
for $\sigma=1$ the upwind flux is obtained. At the physical boundaries we apply always
the upwind flux together with $g$ at the left boundary and the calculated value
at the right boundary.
We  analyze the temporal behavior of the error which is given by $E^k:=u^k(x(\xi),t)-U^k(\xi,t)$. 
We denote by $u^k$ the restriction of $u$ on the $k$-th interval and term $x(\xi)$ denotes
the transformation from the point $\xi$ in the reference interval to point $x$ in the the $k$-th interval.
We can split the error into two parts:
 {\small
\begin{equation}\label{eq:error_basis}
 E^k=\underbrace{(\Ip^N(u^k)-U^k)}_{=:\epsilon^k_1\in \P^N}+
 \underbrace{(u-\Ip^N(u^k))}_{=:\epsilon_p^k }.
\end{equation}}
We will investigate discrete norms in this context using the discrete inner
products \eqref{eq:Quadrature_Lobatto}. 
The global norm, which depends only on $t$, is defined by 
 {\small
\begin{equation}\label{eq:globalnorm}
 ||U(t)||_N^2 := \sum\limits_{k=1}^K \frac{\Delta x_k}{2}  ||U^k(t)||_N^2 =  
 \sum\limits_{k=1}^K \frac{\Delta x_k}{2}  \vec{u}^{k,T}\mat{M}\vec{u}^k,
\end{equation}}
and $U(0)$ is the interpolant of the initial condition $u_0$. 
In \eqref{eq:globalnorm} we provide for the norm 
both the notation with the coefficients and 
with the polynomial solution. They are identical.
With the triangle inequality, we can bound the error \eqref{eq:error_basis} by
 {\small
\begin{equation}\label{eq:error-esti}
 ||E^k||_N\leq ||\epsilon_1^k||_N+||\epsilon_p^k||_N. 
\end{equation}}
The form $\epsilon_p^k$ is the interpolation error, which is the sum of the series truncation error 
and of the aliasing error. 
As it was already described in 
\cite{canuto2012spectral, funaro2008polynomial, hesthaven2008filtering, glaubitz2018application, offner2013spectral, offner2015zweidimensionale }, 
the continuous norms
converge spectrally fast for the different bases under consideration if $u\in C^\infty((-1,1))$.
We denote by 
 {\small
\begin{equation*}
 |u|_{H^{m;N}((-1,1))}:= \l(\sum\limits_{j=\min(m,N+1)}^m ||u^{(j)}||_{\L^2((-1,1))}^2\r)^\frac{1}{2}
\end{equation*}}
the semi-norms of the Sobolev space 
$H^m((-1,1))$, by $P_N$   the projection 
operator of the truncated Legendre series
\footnote{See section 5.4.2 of \cite{canuto2012spectral} for detail.} and by
$\Ip^N$ the  interpolation operator from section \ref{sec2:CPR}. 
We get:
\begin{itemize}
 \item Gauß-Lobatto/Gauß-Legendre points {\small
 \begin{equation}\label{eq:Gauss-esti}
||u-\Ip^N(u)||_{\L^2((-1,1))} \leq C N^{-m}|u|_{H^{m;N}((-1,1))};
 \end{equation}}
\end{itemize}
where $C$ depends on $m$.
The generalization of these formulas 
\eqref{eq:Gauss-esti} and \eqref{eq:Legendre-esti} for $1\leq l\leq m$ are
\begin{itemize}
 \item Gauß-Lobatto/Gauß-Legendre points  {\small
  \begin{equation}\label{eq:Gauss-esti_all}
||u-\Ip^N(u)||_{H^l((-1,1))} \leq C N^{2l-\frac{1}{2}-m}|u|_{H^{m;N}((-1,1))};
 \end{equation}}
\end{itemize}
In view of our investigation, we need to consider our interpolation error not only in the 
standard interval $[-1,1]$, but in each element $e^k$. 
Therefore, we will transform our estimations \eqref{eq:Gauss-esti},\eqref{eq:Legendre-esti} to every element. 
We get with the interval length $\Delta x^k=x^k-x^{k-1}$: 
\begin{itemize}
 \item Gauß-Lobatto/Gauß-Legendre\footnote{ A more a detailed analysis can be found in \cite{bernardi1989properties,bernardi1992polynomial}.
 } points (Combination of \cite[Theorem 6.6.1]{funaro2008polynomial} and \cite[Section 5.4.4]{canuto2012spectral} )
  {\small
 \begin{equation}\label{eq:Gauss-esti_element}
||\epsilon^k_p||_{H^n(e^k)} \leq  C  \l(\Delta x^k\r)^{n-\min\{m,N\}+\frac{1}{2}} N^{n-m+\frac{1}{2}}|u|_{H^{m;N}(e^k)};
 \end{equation}}
 for $n=0,1$. For Gauß-Lobatto, the exponents in \eqref{eq:Gauss-esti_element} are reduced by   $\frac{1}{2}$.
\end{itemize}
 We have introduced all the needed approximation estimations.
\begin{re}
We want to point out that the following investigation can also be done using a modal Legendre basis. 
 Here, we would assume exact integration and the matrix $\mat{M}$ would also be diagonal,
 see \cite{ranocha2017extended} for details. 
 Hence, the interpolation operator in the equation \eqref{eq:error_basis} can be replaced by the projection operator
 and the \emph{interpolation error} is only the series truncation error. 
The analysis for a modal Legendre basis is similar to the
Gauß-Legendre case and can be transferred with 
equivalent estimations 
to \eqref{eq:Gauss-esti}-\eqref{eq:Gauss-esti_element} for the projection operator.
For example, instead of \eqref{eq:Gauss-esti} we apply
  \begin{equation}\label{eq:Legendre-esti}
||u-P_N(u)||_{\L^2((-1,1))} \leq C N^{-m}|u|_{H^{m;N}((-1,1))}.
 \end{equation}
\end{re}

\subsection{Stability of the SBP-FR Methods}\label{subsec:stability_inv}
We follow the steps from \cite{ kopriva2017error} 
and start by repeating the main aspects of the stability 
analysis of  the SBP-FR methods, see \cite{ranocha2016summation} for details.
Then, we derive an error equation for the SBP-FR methods 
for the model problem \eqref{eq:Model} in the following section \ref{sec:Error_equation}. 
In \cite{kopriva2017error} the authors 
investigate the long-time error behavior for 
the Discontiniuous-Galerkin-Spectral-Element Method (DGSEM) using 
Gauß-Lobatto nodes. Here, we make two extensions to their  investigation.
First, we also consider Gauß-Legendre nodes and secondly, 
we also investigate the long-time error behavior of the 
one-parameter family of Vincent et al. where the DGSEM is
included.\\
Instead of using the discrete norm which is represented by $\mat{M}$ 
and corresponds to the continuous $\L^2$-norm we are 
applying  $\mat{M}+\mat{\K}$ analogously to \cite{ranocha2016summation} 
and introduced in section \ref{sec2:CPR}.

We are studying  the change of the discrete norm
 {\small
\begin{equation}\label{eq_discrete_norm_K}
 ||\vec{u}^k||^2_{M+\K}=(\vec{u}^k, \vec{u}^k)_{M+\K}= 
\vec{u}^{k,T} (\mat{M}+\mat{\K}) \vec{u}^k
\end{equation}}
for the total energy.
We multiply $ \vec{\phi}^{k,T}\l(\mat{M}+\mat{\K} \r)$ to equation \eqref{eq:SBP CPR}.
Here, the term $k$ describes the element and $T$ means only the transposed vector. We get
 {\small
\begin{equation}\label{eq:stabilit_gleichung_Vincent}
  \frac{\Delta x_k}{2}  \vec{\phi}^{k,T}\l(\mat{M}+\mat{\K} \r) \partial_t \vec{u}^k=-
  \vec{\phi}^{k,T}\l(\mat{M}+\mat{\K} \r)\mat{D}\vec{u}^k-\vec{\phi}^{k,T}
  \l(\mat{M}+\mat{\K} \r)\mat{C}\l( \vecfnumk-\mat{R}\vec{u}^k
  \r).
\end{equation}}
With $\vec{\phi}^k=\vec{u}^k$ and $\mat{C}=(\mat{M}+\mat{\K})^{-1} \mat{R}[^T] \mat{B}$ it comes:
 {\small
\begin{equation}\label{eq:stabilit_gleichung5.1}
  \frac{\Delta x_k}{2}  \vec{u}^{k,T}\l(\mat{M}+\mat{\K} \r)\partial_t \vec{u}^k=-
  \vec{u}^{k,T}\l(\mat{M}+\mat{\K} \r)\mat{D}\vec{u}^k-\vec{u}^{k,T}
  \mat{R}[^T] \mat{B} \l( \vecfnumk-\mat{R}\vec{u}^k
  \r),
\end{equation} }
where the numerical flux is given by $\vecfnumk= (f^{\mathrm{num,k}}_L, f^{\mathrm{num,k}}_R )^T$.
With $\mat{\K}\mat{D}=0 $ and the SBP property \eqref{eq:SBP_2},
the above equation \eqref{eq:stabilit_gleichung5.1}
can be written as
 {\small
\begin{equation}\label{eq:stabilit_gleichung5.2}
\begin{aligned}
  \frac{\Delta x_k}{2}  \vec{u}^{k,T}\l(\mat{M}+\mat{\K} \r)\partial_t \vec{u}^k 
  = \vec{u}^{k,T} \mat{D}^{T}\mat{M}\vec{u}^k- \vec{u}^{k,T} \mat{R}^T\mat{B}
  \mat{R}\vec{u}^k -\vec{u}^{k,T}  \mat{R}[^T] \mat{B}\l( \vecfnumk-\mat{R}\vec{u}^k\r). 
\end{aligned}
\end{equation}}
Adding  \eqref{eq:stabilit_gleichung5.1} with \eqref{eq:stabilit_gleichung5.2},
using the symmetry of the scalar product induced by $\mat{M}$ and dividing by two
yields 
 {\small
\begin{equation}\label{eq:stability5.3}
\begin{aligned}
\frac{\Delta x_k}{4}\frac{\d}{\d t} ||\vec{u}^k||_{M+\K}^2 
=- \vec{u}^{k,T} 
 \mat{R}^T\mat{B}\l(\vecfnumk- \mat{R}\vec{u}^k \r)
-\frac{1}{2} \vec{u}^{k,T}\mat{R}^T\mat{B}\mat{R} \vec{u}^k. 
\end{aligned}
\end{equation}
}
The rate of change of the total energy is the sum over all elements.
 {\small
\begin{equation}\label{eq:energy5.1}
 \frac{1}{2} \frac{\d}{\d t} \sum\limits_{k=1}^K \frac{\Delta x_k}{2} || \vec{u}^k||^2_{M+\K}
 =-\sum\limits_{k=1}^K \vec{u}^{k,T} \mat{R}^T\mat{B} \l(\vecfnumk-\frac{1}{2}\mat{R} \vec{u}^k
 \r).
\end{equation}}
If we now split the sum into three parts and use the fact that the numerical flux 
is unique at the interface of two elements, we can rewrite\footnote{Details can be
found in \cite{ranocha2016summation, kopriva2017error} as well as in section \ref{sec:Error_equation}.}
\eqref{eq:energy5.1} as
 {\small
\begin{align*}
\frac{1}{2} \frac{\d}{\d t}\sum\limits_{k=1}^K \frac{\Delta x_k}{2}  ||\vec{u}^k||_{M+\K}^2 
=& \frac{1}{2} \frac{\d}{\d t}\sum\limits_{k=1}^K \frac{\Delta x_k}{2}  ||U^k(t)||_{M+\K}^2 =
-\sum\limits_{k=1}^K \vec{u}^{k,T}  \mat{R}[^T] \mat{B}  \l(\vecfnumk-\frac{1}{2}\mat{R} \vec{u}^k \r) \\
=& \frac{1}{2}g(t)^2-\frac{1}{2} \l(U^1_L(t)-g(t) \r)^2-  \frac{1}{2}\l(U^K_R(t)\r)^2- 
\frac{\sigma}{2}\sum\limits_{k=2}^K \l(\jump{U^k(t)} \r)^2,
\end{align*}}
where $U_i$ ($i=L,R$) 
describe the approximated solution \eqref{eq:Approx} and where the indices give
 the position in the elements. The term $\jump{U^k}:=U_R^{k-1}-U_L^k$ is the jump.
 We define  the global norm corresponding to 
  $\mat{M}+\mat{\K}$ by 
$
 ||U(t)||^2_{\K_M}:= \sum\limits_{k=1}^K \frac{\Delta x_k}{2} || U^k(t)||^2_{M+\K}
$ and set $U(0)$ the interpolant of the initial condition $u_0$. The integration from zero to $T$
yields  {\small
 \begin{align*}
  ||U(T)||^2_{\K_M}+&\int\limits_0^T \l(U^K_R(t) \r)^2 \d t +\int\limits_0^T  \l(U^1_L(t)-g(t) \r)^2 +
  \sigma \int\limits_0^T \sum\limits_{k=2}^K \l(\jump{U^k}\r)^2 \d t=||U(0)||_{\K_M}^2 + \int\limits_0^T g^2(t)\d t,
 \end{align*}}
 which also satisfies 
  {\small
\begin{equation}\label{eq:discretstability5.2}
   ||U(T)||^2_{\K_M}+\int\limits_0^T \l(U^K_R(t) \r)^2 \d t \leq ||U(0)||_{\K_M}^2 + \int\limits_0^T g^2(t)\d t.
 \end{equation}}

This has already been proven more generally  in \cite[Theorem 5]{ranocha2016summation}.
Let us mention that the norms under consideration are fully discrete.
However, the same schemes are analyzed as described in \cite{vincent2011newclass}.
Before starting with the analysis in the next section \ref{sec:Error_equation}, 
we shortly repeat  again the notation which will be used in this
paper in table \ref{tab:notation}  for clarification, and we also repeat the assumptions  which are made. 
We introduce some 
abbreviations which are used in the next section.

\renewcommand{\arraystretch}{1.2}
\begin{table}[!ht]
\centering
 {\scriptsize{
\begin{tabular}{l|l}
 {\textbf{\textit{Notation}}} & \textbf{\textit{Interpretation / Assumptions  }}\\ \hline \hline
  $u$ & is the solution of \eqref{eq:Model}. It is $u\in H^m_{\kappa,N}$.\\
  $U$ & is the spatial approximation of $u$ given by \eqref{eq:Approx}. It is $U\in \P^N$. \\
  $\vec{u}$ & \multirow{2}{*}{\shortstack[l]{are the coefficients 
  of $U$, evaluated at the \\ \quad interpolation / quadrature nodes.}}\\
  & \\
  $\mat{D}$ & is the discrete derivative matrix. \\
  $\mat{R}$ & \multirow{2}{*}{\shortstack[l]{is the restriction operator performing interpolation \\ \quad to the boundary.}}\\
  & \\
  $\mat{M}$ & is the diagonal mass / norm matrix.\\
  $\mat{\K}$ & 
  \multirow{2}{*}{\shortstack[l]{ is a symmetric matrix to build different FR schemes. \\ 
  \quad It is defined in \eqref{eq:one_parameter}.  }} \\   & \\
   $\kappa$ & 
   \multirow{2}{*}{\shortstack[l]{  Free parameter to define the different FR schemes from 
  \eqref{eq:one_parameter}, \\ 
  \quad bounded from below and tends to zero for increasing polynomial order $N$. }} \\   & \\
  
  $\est{\cdot,\cdot}$ & is the usual $\L^2$ scalar product.\\
  $||\cdot|| $ & is the norm induced by the $\L^2$ scalar product.\\
  $(\cdot,\cdot)_N$ & is the discrete scalar product given by
  \eqref{eq:Quadrature_Lobatto}.\\
  $||\cdot||_N$ & is the norm induced by the discrete scalar product from above.\\
   $||\cdot||_{M+\K}$ & is the norm 
   \eqref{eq_discrete_norm_K} induced by the discrete scalar product with respect 
   to $\mat{M}+\mat{\K}$.\\
    $||\cdot||_{\K_M} $ &  is global norm corresponding to $\mat{M}+\mat{\K}$. \\
  $\Ip^N$  &  is the interpolation operator. \\
  $P^m_{N-1}(u)$ & \multirow{2}{*}{\shortstack[l]{is the orthogonal projection of $u$
  onto $\P^{N-1}((-1,1))$ using \\ \quad the inner product of $H^m((-1,1))$. } } \\
  & \\
  $E^k= u^k-U^k$ & is the total error in the $k$-th element.\\
  $\epsilon_1^k:=\Ip^N(u^k)-U^k$ & \multirow{2}{*}
  {\shortstack[l]{is the difference between interpolation and spatial \\ 
  \quad approximation in the $k$-th element.}} \\
  & \\
  $\epsilon_p^k=u^k-\Ip^N(u^k)$ & is the interpolation error. 
\end{tabular}
}}
\caption{Summary of the notations used in this article.}
\label{tab:notation}
\end{table}
\renewcommand{\arraystretch}{1}


\section{Long-time Error Behavior for SBP-FR}\label{sec:Error_equation}

The  error is investigate in respect to the global discrete norm $||\cdot||_{\K_M}$.
It is given by $E^k= u^k(x(\xi),t)-U^k(\xi,t) $ within any element. 
Using \eqref{eq:error_basis} and the triangle inequality, we can bound the error by
\begin{equation}\label{eq:error-esti_Vincent}
 ||E^k||_{M+\K}\leq ||\epsilon_1^k||_{M+\K}+||\epsilon_p^k||_{M+\K}. 
\end{equation}

Instead of focusing on both terms in \eqref{eq:error-esti_Vincent},
we derive the error equation only for 
 $\epsilon_1^k$ since this 
term has the major influence on the error behaviors.
First, $\epsilon_p^k$ is simply an interpolation error of
the exact solution, and so there is actually no process /mechanism 
in the procedure which might lead to any accumulation of the 
$\epsilon_p^k$ errors.\\
Secondly, by following the ideas / steps of \cite{kopriva2017error},
we do not further consider  the terms $\epsilon_p^k$ due to the following fact.
We are working in a finite dimensional normed vector space. 
All norms are equivalent in this vector space and this allows us to bound 
the discrete norm 
{\small 
\begin{equation}\label{eq_discrete_norm_K_2}
 ||\vec{\epsilon_p}^k||^2_{M+\K}=(\vec{\epsilon_p}^k, \vec{\epsilon_p}^k)_{M+\K}= 
\vec{\epsilon_p}^{k,T} (\mat{M}+\mat{\K}) \vec{\epsilon_p}^k
\stackrel{\eqref{eq:one_parameter}}{=}
\vec{\epsilon_p}^{k,T}\mat{M} \vec{\epsilon_p}^k 
+\vec{\epsilon_p}^{k,T}\kappa(\mat{D}^N)^T \mat{M}\mat{D}^N \vec{\epsilon_p}^k
\end{equation}}
in terms of the continuous ones from the Sobolev space $H^m$. 
Since $u$ is sufficiently smooth (i.e. $u\in H_{\kappa,N}^m$ with $m>1$), 
the continuous norms can be estimated by 
\eqref{eq:Gauss-esti} - \eqref{eq:Gauss-esti_element} and 
we obtain a spectral decay for the interpolation error $\epsilon_p^k$.
In other words, we have to investigate the error behavior of $\epsilon_1^k$ in detail.
Here, we follow the approach from \cite{kopriva2017error}. \\
Before we derive the error equation for $\epsilon_1^k$, we give the main 
result along with a short description of the main steps of the analysis
\begin{Result}
The factor $\eta(t)$ given in $\eqref{eq:Errorbound2}$ depends on $\epsilon_1$.
 If the mean of $\eta$ can be bounded 
 from below by a positive constant $\delta_0$, i.e. $\mean{\eta}\geq \delta_0>0$, then there exists a constant $C_1$ such that the errors $\epsilon_1^k(t)$ of 
 \eqref{eq:error-esti_Vincent} satisfies the inequality
 {\small
\begin{equation}\label{eq:Errorbound3}
 ||\epsilon_1(t)||_{\K_M} \leq \frac{1-\exp(-\delta_0 t)}{\delta_0} C_1
\end{equation}}
in the  discrete norm $||\cdot||_{\K_M}$. The total error is bounded in time.
\end{Result}
In the following, the exact conditions for the above inequality to be fulfilled will be derived 
and we specify in detail what factors play a key role in the definition of $\eta$ and $\delta_0$.
We outline the major steps of our analysis:
\begin{enumerate}
 \item Inserting the error $E^k$ into the continuous equation in every element leads us to an error 
 equation for $\epsilon_1^k$. 
 \item Adding zero in a suitable way gives us the possibility to split the equations into 
 a continuous and a discrete part.
 \item We add both parts for every element and obtain the error behavior for the total domain.
 \item By estimating the continuous terms 
 we get an inequality for the error $\epsilon_1$ in the discrete norms
 and with some assumptions we retrieve the long-time error behavior of $\epsilon_1$.
\end{enumerate}

We derive now the error equation for $\epsilon_1^k=\Ip^N(u^k)-U^k$.

We are searching the solution $u$ of the continuous equation
 {\small
\begin{equation}\label{eq:Continuous}
 \frac{\Delta x}{2}\est{\partial_t u,\phi^k}+u\phi\Big|_{-1}^1-\est{u,\partial_\xi \phi^k}=0,
\end{equation}
}
where $\est{u,\phi^k}:=\int_{-1}^1 u\phi^k \d \xi$ defines the inner product. The equation \eqref{eq:Continuous} can be derived from the advection equation
\eqref{eq:Model} by multiplying  with the test function $\phi$, 
integrating over the standard element and using 
integration-by-pars. 
With $\phi^k \in \P^N\subset \L^2 $ and $u^k=\Ip^N(u^k)+\epsilon_p^k$,
we get for the continuous equation
 {\small
\begin{equation}\label{eq:Continuous2}
 \frac{\Delta x_k}{2} \est{\partial_t \Ip^N(u^k),\phi^k }+
 \Ip^N(u^k)\phi^k \Big|_{-1}^1- 
 \est{\Ip^N(u^k), \partial_\xi\phi^k }
 =-\frac{\Delta x_k}{2} \est{\partial_t \epsilon_p^k, \phi^k } -\epsilon_p^k\phi^k \Big|_{-1}^1 +
\est{\epsilon_p^k, \partial_\xi\phi^k }.
 \end{equation}
 }
\begin{re}
For Gauß-Lobatto nodes  it holds $\epsilon_p^k=0$ at the endpoints
 because the interpolant is equal to the solution there. 
 Thus,  $\epsilon_p^k\phi^k\Big|_{-1}^1=0$.  
\end{re}
Using integration-by-parts for  $ \est{ \epsilon_p^k,\partial_\xi \phi^k} $ 
yields
 {\small
\begin{equation}\label{eq:Error1}
 \frac{\Delta x_k}{2} \est{\partial_t \Ip^N(u^k),\phi^k }+
 \Ip^N(u^k)\phi^k \Big|_{-1}^1- \est{\Ip^N(u^k), \partial_\xi\phi^k }
 =-\frac{\Delta x_k}{2} \est{\partial_t \epsilon_p^k, \phi^k } 
 -\est{ \partial_\xi (\epsilon_p^k), \phi^k}.
\end{equation}}
Applying now interpolation, the discrete norm gives for the first term
 {\small{
\begin{equation}\label{eq:interpolation2_Vincent}
 \est{\partial_t \Ip^N(u^k),\phi^k } = 
 \l( \partial_t \vec{\Ip^N(u^k)}, \vec{\phi}^k \r)_{M+\K}
 +\l\{ \est{ \partial_t \Ip^N(u^k), \phi^k}-
 \l( \partial_t \vec{\Ip^N(u^k)}, \vec{\phi}^k \r)_{ M+\K}  \r\}.
\end{equation}}}

Since $\phi^k \in \P^N$ and the exactness of the quadrature formulas, we obtain 
for the volume  term in \eqref{eq:Error1}
 {\small{
\begin{equation}\label{eq:interpolation_Vincent}
\begin{aligned}
 \est{ \Ip^N(u^k), \partial_\xi \phi^k}
 =  \l(  \vec{\Ip^N(u^k)}, \partial_\xi \vec{\phi}^k \r)_{N}.
\end{aligned}
\end{equation}
Finally, the values of the interpolation polynomial at the
boundaries of the element ($-1$ and $1$) can be approximated 
by a limitation process from 
the left side $\Ip^N(u^k)^{-}$ and right side $\Ip^N(u^k)^{+}$. 
To simplify the notation we define
 {\small{
\begin{equation}
 \vecfnumk \l(\Ip^N(u^k)^{-},\Ip^N(u^k)^{+}\r) :=
 \l(\fnum \l(\Ip_R^N(u)^{k-1},\Ip_L^N(u)^{k}\r), 
 \fnum \l(\Ip_R^N(u)^{k},\Ip_L^N(u)^{k+1}\r)      \r)^T.
\end{equation}}}
We obtain for the approximation 
 {\small
\begin{align*}
 \Ip^N(u^k)\phi^k\bigg|_{-1}^1=& \vec{\phi}^{T,k} \mat{R}^T \mat{B} 
 \vecfnumk \l(\Ip^N(u^k)^{-},\Ip^N(u^k)^{+}\r) \\
 &+ \l( \phi^k \Ip^N(u^k)\bigg|_{-1}^1  
 -\vec{\phi}^{T,k} \mat{R}^T \mat{B} \vecfnumk \l(\Ip^N(u^k)^{-},\Ip^N(u^k)^{+}\r)  \right).
\end{align*}}
$u$ is continuous $(m>1)$. Using Gauß-Lobatto points the error term in the 
braces is zero, because the interpolation polynomial is 
evaluated at these boundaries and the numerical flux is unique.
For Gauss-Legendre points, we get an additional error term
which corresponds to an interpolation error (in the pointwise sense) at these end points. 
The numerical flux is again unique and so the error term  reads 
 {\small
\begin{equation}\label{eq:error_epislon_2}
\epsilon_{2,\phi}^k:=  \l( \phi^k \Ip^N(u^k)\bigg|_{-1}^1  
 -\vec{\phi}^{T,k} \mat{R}^T \mat{B} \vecfnumk \l(\Ip^N(u^k)^{-},\Ip^N(u^k)^{+}\r)  \right)  
\end{equation}}
in the $k$-th element. Finally, using \eqref{eq:interpolation_Vincent}-\eqref{eq:error_epislon_2}
in \eqref{eq:interpolation2_Vincent} we obtain 
 {\small
\begin{equation}\label{eq:Error2_Vincent}
\begin{aligned}
 &\frac{\Delta x_k}{2} \left(\partial_t \vec{\Ip^N(u^k)},\vec{\phi}^k \right)_{M+\K}+
 \vec{\phi}^{T,k} \mat{R}^T \mat{B} \vecfnumk \l(\Ip^N(u^k)^{-},\Ip^N(u^k)^{+}\r) 
 -\l( (\vec{\Ip^N(u^k)}, \partial_\xi \vec{\phi}^k \r)_{M+\K}\\
 =&-\frac{\Delta x_k}{2} \est{\partial_t \epsilon_p^k, \phi^k } 
 -\est{ \partial_\xi (\epsilon_p^k), \phi^k }
 -\frac{\Delta x_k} {2} \l\{ \est{ \partial_t \Ip^N(u^k), \phi^k}
 -\l( \partial_t \vec{\Ip^N(u^k)}, \vec{\phi}^k \r)_{M+\K}  \r\} 
 -\epsilon_{2,\phi}^k.
\end{aligned}
\end{equation}}
Adding zero to the terms in the curly braces and using
 {\small
\begin{align*}
  \est{ \partial_t \Ip^N(u^k), \phi^k}- \l( \partial_t \vec{\Ip^N(u^k)}, 
  \vec{\phi}^k \r)_{M+\K} =& \est{\underbrace{ \partial _t
  \l( \Ip^N(u^k)-P^m_{N-1} \l(\Ip^N(u^k)\r)\r)}_{=:Q(u^k)}, \phi^k }\\
  &- \l( \partial_t \l( \vec{\Ip^N(u^k)}-\vec{P^m_{N-1} \l(\Ip^N(u^k)\r) }  \r), \vec{\phi}^k   \r)_{M+\K} ,
\end{align*}}
with $P^m_{N-1}$  the orthogonal projection operator\footnote{
The projection operator is defined by the 
classical truncated Fourier series $P^m_{N-1}u=\sum_{j=0}^{N-1} \hat{u}_j\Phi_j$ up to order $N-1$ where
(broken) Sobolev type orthogonal polynomials $\{\Phi_k\}$ are used as basis function in the 
underlying space. Essential is that it projects $u$ unto $\P^{N-1}$ 
and having the representation 
\eqref{eq:one_parameter} in mind. For more details about the projection operator and about 
approximation results, we strongly recommend \cite{canuto2012spectral}.
 }
of $u$ onto $\P^{N-1}$ 
yields in \eqref{eq:Error2_Vincent}
 {\small
\begin{equation} \label{eq:Error3_Vincent}
\begin{aligned}
 &\frac{\Delta x_k}{2} \left(\partial_t \vec{\Ip^N(u^k)},\vec{\phi}^k \right)_{M+\K}+
 \vec{\phi}^{T,k} \mat{R}^T \mat{B} \vecfnumk \l(\Ip^N(u^k)^{-},\Ip^N(u^k)^{+}\r) 
 -\l( \vec{\Ip^N(u^k)}, \partial_\xi \vec{\phi}^k \r)_{M+\K} \\
 =&\frac{\Delta x_k}{2} \est{(T^k(u),\phi^k }+\frac{\Delta x_k }{2}  
 \left( \vec{Q(u^k)},\vec{\phi}^k \r)_{M+\K}-\epsilon_{2,\phi}^k,
\end{aligned}
\end{equation}}
where $
 T^k(u)= -\left\{ \partial_t \epsilon_p^k+ \partial_x \epsilon_p^k+ Q(u^k) \r\}.
$\\
$Q$ measures the projection error of a polynomial
of degree $N$ to a polynomial of degree $N-1$.
Since $u$ is bounded, this 
value has also to be bounded. 
Since \eqref{eq:Gauss-esti_element} and $\kappa\to 0$, 
the interpolation error $\epsilon_p^k$ converges in $N$ to zero, provided that $m>1$ and that
the Sobolev norm of the solution is uniformly bounded in time.  
Therefore, we also need  the initial and boundary conditions in the
model problem \eqref{eq:Model}.
For the time derivative, we get the boundedness of the norm by the relation $\partial_t u=-\partial_x u$. 
The term $\epsilon^k_{2,\phi}$ is bounded, because $u$ is bounded and also continuous. 
For the numerical fluxes, this value 
describes the error between the interpolation polynomial at $-1$ and $1$, as well as 
the numerical approximation by the numerical flux function at these boundaries. 
From a different perspective, this value can also be interpreted as the additional
dissipation which is added in the Gauß-Legendre case, since for Gauß-Lobatto nodes
this error term is zero. This discussion yields that the right side of 
\eqref{eq:Error3_Vincent} is well-defined. \\
Now, we derive the error equation for $\epsilon_1^k$. 
We apply the SBP property \eqref{eq:SBP},  
$\mat{C}=\l(\mat{M}+\mat{\K} \r)^{-1}\mat{R}^T\mat{B} $ 
and the property of  \eqref{eq:one_parameter}
to equation \eqref{eq:stabilit_gleichung_Vincent}
and obtain
 {\small
\begin{equation*}
 \frac{\Delta x_k}{2} \left(\partial_t \vec{u}^k   ,\vec{\phi}^k \right)_{M+\K}+
 \vec{\phi}^{T,k} \mat{R}^T \mat{B}\vecfnumk \l(\l(U^k\r)^-,\l(U^k\r)^+\r) \\
 -\l( \vec{u}^k , \partial_\xi \vec{\phi}^k \r)_{M+\K}
 =0.
\end{equation*}
}
We subtract this equation from \eqref{eq:Error3_Vincent} and by
the linearity of the numerical flux, we get an equation for the error
$\epsilon_1^k=\Ip^N(u^k)-U^k$. It becomes:
 {\small
\begin{align*}
 &\frac{\Delta x_k}{2} \left(\partial_t \l( \vec{\Ip^N(u^k)}-\vec{u}^k  \r) ,\vec{\phi}^k \right)_{M+\K}+
 \vec{\phi}^{T,k} \mat{R}^T \mat{B}\vecfnumk \l(\l( \Ip^N(u^k)-U^k\r)^-,\l(\Ip^N(u^k)-U^k\r)^+\r) \\
 -&\l( \l( \vec{\Ip^N(u^k)}-\vec{u}^k \r), \partial_\xi \vec{\phi}^k \r)_{M+\K}
 =\frac{\Delta x_k}{2} \est{T^k(u),\phi^k }+\frac{\Delta x_k }{2} 
 \left( \vec{Q(u^k)},\vec{\phi^k} \r)_{M+\K}-\epsilon_{2,\phi}^k,\\ 
 \Longleftrightarrow& \frac{\Delta x_k}{2} 
 \left(\partial_t \vec{\epsilon}_1^k ,\vec{\phi}^k \right)_{M+\K}+
 \vec{\phi}^{T,k} \mat{R}^T \mat{B} \vecfnumk \l(\l( \epsilon_1^k\r)^-,\l(\epsilon_1^k\r)^+\r) 
 -\l(  \vec{\epsilon}_1^k , \partial_\xi \vec{\phi}^k \r)_{M+\K} \\
 =&\frac{\Delta x_k}{2} \est{T^k(u),\phi^k }+\frac{\Delta x_k }{2} 
 \left( \vec{Q(u^k)},\vec{\phi}^k \r)_{M+\K}-\epsilon_{2,\phi}^k ,
\end{align*}}
where we  can write for the term 
$\l(  \vec{\epsilon}_1^k , \partial_\xi \vec{\phi}^k \r)_{M+\K}=
\l(  \vec{\epsilon}_1^
k , \partial_\xi \vec{\phi}^k \r)_{N}$ because  of \eqref{eq:one_parameter}.
Putting $\phi^k=\epsilon_1^k$, we obtain the  energy equation 
 {\small
\begin{align*}
 &\frac{\Delta x_k}{4} \frac{\d}{\d t} || \epsilon_1^k||_{M+\K}^2+ \vec{ \epsilon}_1^{k,T} \mat{R}^T \mat{B} \vecfnumk \l(\l( \epsilon_1^k\r)^-,\l(\epsilon_1^k\r)^+\r)
 -\l(\vec{\epsilon}_1^k, \partial_\xi \vec{\epsilon}_1^k \r)_N\\
 =& \frac{\Delta x_k}{2} \est{T^k(u),\epsilon_1^k }+\frac{\Delta x_k }{2} 
 \left( \vec{Q(u^k)},\vec{\epsilon}_1^k \r)_{M+\K}-\tilde{\epsilon}_2^k,
\end{align*}}
with $\tilde{\epsilon}_2^k =  \l( \epsilon_1^k \Ip^N(u^k)\big|_{-1}^1  
 -\vec{ \epsilon}_1^{k,T}\mat{R}^T \mat{B} \vecfnumk \l(\Ip^N(u^k)^{-},\Ip^N(u^k)^{+}\r)  \right)$.
 Summation-by-parts yields for $\epsilon_1^k$
  {\small
 \begin{align*}
  \l(\vec{\epsilon}_1^k, \partial_\xi \vec{\epsilon}_1^k \r)_N =\vec{\epsilon_1}^{T,k}\mat{M}\mat{D} \vec{\epsilon_1}^k&= \vec{\epsilon_1}^{T,k} \mat{R}^{T} \mat{B}\mat{R} \vec{\epsilon_1}^k-
  \vec{\epsilon_1}^{T,k}\mat{D}^T\mat{M}\vec{\epsilon_1}^k,\\
  \Longleftrightarrow \quad \vec{\epsilon_1}^{T,k}\mat{M}\mat{D} \vec{\epsilon_1}^k&= \frac{1}{2}\vec{\epsilon_1}^{T,k} \mat{R}^{T} \mat{B}\mat{R} \vec{\epsilon_1}^k,
 \end{align*}}
and we get  {\small
 \begin{align*}
 &\frac{\Delta x_k}{4} \frac{\d}{\d t} || \epsilon_1^k||_{M+\K}^2+ \vec{ \epsilon}_1^{k,T} \mat{R}^T \mat{B} \l( \vecfnumk \l(\l( \epsilon_1^k\r)^-,\l(\epsilon_1^k\r)^+\r) -\frac{1}{2} \mat{R}\vec{\epsilon_1}^k
 \r)\\ =& \frac{\Delta x_k}{2} \est{T^k(u),\epsilon_1^k }+\frac{\Delta x_k }{2}  \left( \vec{Q(u^k)},\vec{\epsilon}_1^k \r)_{M+\K}-\tilde{\epsilon}_2^k.
\end{align*}}
We have split our equation into a continuous and a discrete part. Coming to step three
of our investigation, we sum up over all elements and obtain
 {\small
\begin{equation}\label{eq:energy2_Vincent}
\begin{aligned}
  \frac{1}{2} \frac{\d}{\d t} \sum\limits_{k=1}^K \frac{\Delta x_k}{2} ||\epsilon_1^k||^2_{M+\K}
  &+\sum\limits_{k=1}^K \vec{ \epsilon}_1^{k,T}
  \mat{R}^T \mat{B} \l(\vecepsilon -\frac{1}{2} \mat{R}\vec{\epsilon_1}^k \r)\\
 &= \sum\limits_{k=1}^K \l(\frac{\Delta x_k}{2} 
 \l( \est{T^k(u),\epsilon_1^k }+\left( \vec{Q(u^k)},\vec{\epsilon}_1^k \r)_{M+\K} \r)-\tilde{\epsilon}_2^k  \r),
\end{aligned}
\end{equation}
}
where $\vecepsilon:=\vecfnumk \l( \l( \epsilon_1^k\r)^-,\l(\epsilon_1^k\r)^+\r)$. 
This equation has the same form as  \eqref{eq:energy5.1} except the right hand side. 
We estimate the bracket on the right hand side by the Cauchy-Schwarz inequality. 
It is 
 {\small
 \begin{align*}
 \tilde{R}=& \sum\limits_{k=1}^K \frac{\Delta x_k}{2} 
 \l( \est{T^k(u),\epsilon_1^k }+\left( \vec{Q(u^k)},\vec{\epsilon}_1^k \r)_{M+\K} \r) \\
 \leq& \sqrt{\sum\limits_{k=1}^K \frac{\Delta x_k}{2}  ||T^k(u)||^2 }  
 \sqrt{\sum\limits_{k=1}^K \frac{\Delta x_k}{2}  ||\epsilon_1^k||^2 }   
+\sqrt{\sum\limits_{k=1}^K \frac{\Delta x_k}{2}  
||\vec{Q(u^k)}||_{M+\K}^2 } 
\sqrt{\sum\limits_{k=1}^K \frac{\Delta x_k}{2}  ||\vec{\epsilon}_1^k||_{M+\K}^2 }  . 
\end{align*}}
With the global norm over all 
elements and the equivalence between the continuous and discrete norms, we obtain 
 {\small
\begin{equation}\label{eq:R_absch_Vincent}
 \tilde{R}\leq \l\{ c||T|| +||Q||_{\K_M}  \r\} ||\epsilon_1||_{\K_M}
=: \tilde{\Ep}(t) ||\epsilon_1||_{\K_M}.
\end{equation}
}

Using  estimation \eqref{eq:R_absch_Vincent} in \eqref{eq:energy2_Vincent}, we get an inequality for the 
global energy equation for the total error. It is
 {\small
\begin{equation}\label{eq:error_equationglobal_Vincent}
  \frac{1}{2} \frac{\d}{\d t} ||\epsilon_1||^2_{\K_M} +\sum\limits_{k=1}^K \vec{ \epsilon}_1^{k,T}
  \mat{R}^T \mat{B} \l(\vecepsilon -\frac{1}{2} \mat{R}\vec{\epsilon_1}^k \r) \leq \tilde{\Ep}(t)  
  ||\epsilon_1||_{\K_M} -\underbrace{\sum_{k=1}^K\tilde{\epsilon}_2^k}_{:=\Theta_2} 
\end{equation}}
with $\tilde{\epsilon}_2^k =  \l( \epsilon_1^k \Ip^N(u^k)\big|_{-1}^1
 -\vec{ \epsilon}_1^{k,T}\mat{R}^T \mat{B} \vecfnumk \l(\Ip^N(u^k)^{-},\Ip^N(u^k)^{+}\r)  
  \right)$.
Applying the same approach like in \cite{kopriva2017error, offner2019error} and splitting the sum on 
the left side into three parts
(one for the left physical boundary,
one for the right physical boundary) and  summing over the internal element endpoints),
we get
 {\small
\begin{align*}
 &\sum\limits_{k=1}^K \vec{ \epsilon}_1^{k,T}
  \mat{R}^T \mat{B} \l(\vecepsilon -\frac{1}{2} \mat{R}\vec{\epsilon_1}^k \r) =
 \sum\limits_{k=1}^K \vec{ \epsilon}_1^{k,T}
  \mat{R}^T \mat{B} \l(\vecfnumk \l( \l( \epsilon_1^k\r)^-,\l(\epsilon_1^k\r)^+\r) -\frac{1}{2} \mat{R}\vec{\epsilon}_1^k \r)\\
 & =- \Eps_L^1 \l( f^{\mathrm{num},1}_L -\frac{1}{2} \Eps_L^1 \r)+ \sum\limits_{k=2}^K \l( f^{\mathrm{num},k}_L -\frac{1}{2} \left(\Eps_R^{k-1}+\Eps_L^{k}  \right) \r)
  \l( \Eps_R^{k-1}-\Eps_L^{k}  \r) \\
  & + \Eps_R^K \l(f^{\mathrm{num},K}_R -\frac{1}{2} \Eps^K_R \r).
\end{align*}
}
We describe with $\Eps_i$ ($i=L,R$) the approximated error $\epsilon_1$
 and we have 
 $f^{\mathrm{num},k}_L:=f^{\mathrm{num},k} 
 \l(\Eps^{k-1}_R,\Eps^{k}_L \r)$, $ f^{\mathrm{num},1}_L:=f^{\mathrm{num},1} \l(0,\Eps^{1}_L \r)$ 
 and $f^{\mathrm{num},K}_R:= f^{\mathrm{num},1} \l(\Eps^K_R,0 \r) $.
 The external states for the physical boundary contributions are zero
 because $\Ip^N(u^1)=g$ at the left boundary. The external state for $U^1$ 
 is set to $g$. 
 At the right boundary, where the upwind numerical flux is used,
 it does not matter what we set for the external state
 because its coefficients in the numerical solution is zero.
 We get for the inner element with 
 $\jump{\Eps^k}=\Eps_R^{k-1}-\Eps_L^{k}$; 
  {\small
 \begin{equation*}
 \begin{aligned}
  \sum\limits_{k=2}^K \l( f^{\mathrm{num},k}_L -\frac{1}{2} \left(\Eps_R^{k-1}+\Eps_L^{k}  \right) \r)
  \l( \Eps_R^{k-1}-\Eps_L^{k}  \r)& = \sum\limits_{k=2}^K \frac{\sigma}{2}\l(\jump{\Eps^{k} }\r)^2\geq0,
  &\text{with}  \begin{cases}
                                                                                                                        \sigma=0 \quad \text{central flux, }\\
                                                                                                                        \sigma=1 \quad \text{ upwind flux. }                                                                                                       \end{cases}
 \end{aligned}
 \end{equation*}
 }
For the left and right boundaries, we finally get
 {\small
\begin{align*}
 &\text{left:} &-\Eps_L^1 \l( f^{\mathrm{num},1}_L -\frac{1}{2} \Eps_L^1 \r)=-\Eps_L^1  \left( \left( \frac{0+\Eps_L^1 }{2}-\sigma \frac{\Eps_L^1 }{2} \right)-\frac{\Eps_L^1 }{2} \right)= \frac{\sigma}{2}\l(\Eps_L^1\r)^2,  \\
 &\text{right:}&\Eps_R^K \l(f^{\mathrm{num},K}_R -\frac{1}{2} \Eps^K_R \r)= \Eps^K_R  \l( \l( \frac{0+\Eps^K_R }{2} +\frac{1}{2} \sigma \Eps^K_R  \r)-\frac{\Eps^K_R }{2}\r)= \frac{\sigma}{2} \l(\Eps^K_R  \r)^2.
\end{align*}
}
Therefore, the energy growth rate is bounded by 
 {\small
 \begin{equation}\label{eq:Errorbound}
  \frac{1}{2} \frac{\d}{\d t} ||\epsilon_1||^2_{\K_M} +  \underbrace{\frac{\sigma}{2} \l( \l(\Eps^K_R  \r)^2 +\l(\Eps_L^1\r)^2  \r) 
  +\frac{\sigma}{2} \sum\limits_{k=2}^K \l(\jump{\Eps^{k} }\r)^2 }_{BTs } 
  \leq \tilde{\Ep}(t)||\epsilon_1||_{\K_M}-\Theta_2.
\end{equation}
}
The term  $BTs$ is bigger than or equal to zero depending on the used fluxes. 
The energy growth inequality \eqref{eq:Errorbound}
is similar to the inequality in \cite{kopriva2017error}.
The differences are the used norms and the term $\Theta_2$ which will yield 
a smaller upper bound under the condition $\Theta_2\geq0$. 
We follow the steps in 
\cite{kopriva2017error, nordstrom2007error}
and get 
\begin{equation}\label{eq:Errorbound2}
\frac{\partial}{\partial t} ||\epsilon_1||_{\K_M} +  
\underbrace{ \frac{BTs+\Theta_2}{  ||\epsilon_1||_{\K_M}^2} }_{\eta(t) } ||\epsilon_1||_{\K_M} 
\leq \tilde{\Ep}(t).
\end{equation}
Like it was described in \cite{nordstrom2007error},
we have to assume that the mean value of $\eta(t)$ is bounded by a
positive constant $\delta_0$ from below.
This means that 
$\overline{\eta}\geq \delta_0>0$. Under the assumption on $u$, the 
right hand side $\tilde{\Ep}(t)$ is bounded in time and we can put 
$\max\limits_{t\in [0,\infty)} \tilde{\Ep}(t)\leq C_1<\infty$.
Applying these facts in \eqref{eq:Errorbound2},
we integrate over the time and get the following inequality for the error 
\begin{equation*}
 ||\epsilon_1(t)||_{\K_M} \leq \frac{1-\exp(-\delta_0 t)}{\delta_0} C_1, \tag{\ref{eq:Errorbound3}}
\end{equation*}
see \cite[Lemma 2.3]{nordstrom2007error} for details. 
If $\Theta_2>0$, the existence of $\delta_0$ is actually met without 
restrictions.\\
On the physical boundaries, we apply always the upwind flux ($\sigma=1$).
Therefore,
we may modify the $BTs$ term in \eqref{eq:Errorbound} and we have
\begin{equation*}
 BTs= \frac{1}{2} \l( \l(\Eps^K_R  \r)^2 +\l(\Eps_L^1\r)^2  \r) 
  +\frac{\sigma}{2} \sum\limits_{k=2}^K \l(\jump{\Eps^{k} }\r)^2 >0.
\end{equation*}
In case that $BT>-\Theta_2$ the assumption on the existence of a positive constant $\delta_0$ 
is therefore always fulfilled. 
We transferred the results from \cite{kopriva2017error, nordstrom2007error}
to the more general case of the one parameter family of Vincent et al. \eqref{eq:one_parameter}
and extended the basis also to Gauß-Legendre. 
We may conclude:\\
\emph{If the truncation error is bounded, the dissipative boundary conditions 
keep also the error bounded in time for both Gauß-Lobatto as well as Gauß-Legendre nodes.
The selections of basis and numerical fluxes have an essential influence on the error behaviour.
}\\
In \cite[p.325]{kopriva2017error} their model\footnote{Using Gauß-Lobatto nodes and
investigating the DGSEM (i.e. $\kappa=0$ in \eqref{eq:one_parameter}).}
\eqref{eq:Errorbound3} yields the authors 
the following predictions:
\begin{enumerate}
 \item[P1] \emph{Using the upwind flux at the physical boundaries and either the upwind
 flux or the central
 flux at the interior element interfaces, the error growth is bounded asymptotically in time.}
 \item[P2] \emph{Using the upwind flux $\sigma = 1$ in the interior will lead to a smaller asymptotic error than
using the central flux, $\sigma = 0$. This will be especially true for under-resolved approximations.}
\item[P3] \emph{As the resolution increases, the difference between the asymptotic error from the central
and upwind fluxes should decrease.}
\item [P4] \emph{The error growth rate will be larger when the upwind flux is used compared to when the
central flux is used. Equivalently, the upwind flux solution should approach its asymptotic
value faster than the central flux solution.}
\end{enumerate}
Through our investigation in this section and by the model \eqref{eq:Errorbound3}, these predictions are also valid
using Gauß-Legendre nodes and we can  extend these by the following: 
\begin{enumerate}
 \item [P5] \emph{The error should be larger when Gauß-Lobatto nodes are used compared to when 
 a Gauß-Legendre basis is applied.}\\
 Applying Gauß-Legendre nodes is more accurate than using Gauß-Lobatto nodes. 
 Therefore, the norm of the $\epsilon_1$-error is smaller by utilizing a Gauß-Legendre 
 basis. This leads directly to a bigger value of $\mean{\eta}$ and thus 
 to a bigger $\delta_0$ in \eqref{eq:Errorbound3}. 
 Furthermore, the $\epsilon_2$-error may have a positive effect on the error behavior. 
 \item [P6] \emph{By using Gauß-Legendre nodes, the choice of the numerical fluxes should be less important then 
 in the Gauß-Lobatto case. }\\
 This is a direct consequence of P3 and P5. 
 \item [P7] \emph{By applying a FR scheme with $\kappa\neq0$, the errors should show some oscillations.}\\
 For $\kappa\neq 0$, the correction term with $\K\neq 0$ works directly on the highest
 degree monomials in $U$ with different strength which leads to oscillations. 
 As it is shown in \cite{allaneau2011connections},
 the correction functions \eqref{eq:one_parameter} correspond 
 to top-mode filters of different strengths.
 \end{enumerate}

\begin{re}\label{remark_epsilon}
 If $\Theta_2\geq0$, the term can be seen as 
an additional dissipation term which is added in the Gauß-Legendre case
and lead to a smaller upper bound. Furthermore, because of the construction 
of  of $\tilde{\epsilon_2}$ and the evaluation 
at the boundaries, we suppose that the $\epsilon_1$ error has a direct influence
on the behavior of $\tilde{\epsilon_2}$ and we also get some noisy behavior.
Finally, it is 
\begin{equation*}
 \tilde{\epsilon}_2^k=\vec{ \epsilon}_1^{k,T}\mat{R}^T \mat{B}\l(\mat{R} 
  \vec{\Ip^N(u^k)} -   \vecfnumk \l( \Ip^N(u^k)^{-}, \Ip^N(u^k)^{+} \r) \r) .
\end{equation*}
The term in the bracket describes the error between the flux function 
and the numerical flux function at the element boundaries. 
In case of a consistent numerical flux these terms tend to zero
under mesh refinement and/or increasing the polynomial order of the 
approximation. This can be shown  by Taylor series expansion but is not the topic here.
In the next section, the term $\tilde{\epsilon}_2$ will be investigated numerically.\\
In our investigation, we apply the discrete norms \eqref{eq_discrete_norm_K}
of the one parameter family of Vincent et al. and we assume in \ref{As:c_zero} that $\kappa$ in 
\eqref{eq_discrete_norm_K} 
tends rapidly to zero if we increase the polynomial order $N$. 
The nowadays main used FR schemes (presented in table \ref{ta:correction_terms})
fulfill the assumption \ref{As:c_zero}. 
The first term in the norm is the essential one.
 For the norm of the interpolation error $||\vec{\epsilon_p}^k||_{\K+M}$,
we estimate the first part by 
\eqref{eq:Gauss-esti} -\eqref{eq:Gauss-esti_all}
and  due to assumption \ref{As:c_zero}
 we may neglect the terms of the interpolation errors in the investigation. 
Nevertheless, the 
stability analysis of \cite{vincent2011newclass} allows  $\kappa$ to tend to infinity. 
Therefore, in the case $\kappa\to \infty$ a more detailed analysis is necessary.
We may estimate the second part of the norm \eqref{eq_discrete_norm_K}
using Bernstein inequality together with estimation \eqref{eq:Gauss-esti}, or 
directly \eqref{eq:Gauss-esti_all}.\\ 
A further investigation about the approximation 
behavior concerning  limit processes of $K,N,m, \kappa \to \infty$
would be indeed desirable for this case. 
Both is beyond the scope 
of this paper where the influence of the flux functions and nodal
bases is investigated.\\
We want to point out that spectral convergence is always investigated under the requirement
that our solution  $u$ is $C^{\infty}$ and $m$ from \eqref{eq:Gauss-esti_all}
tends to infinity for all orders of accuracy \cite{canuto2012spectral}.
It is clear that our approach can be easily transformed to multidimensional
problems using a tensor product structure on structured grids.
\end{re}

\section{Numerical Tests}\label{6_Numerical}

In this section we consider numerical tests which demonstrate both the error 
bound \eqref{eq:Errorbound3} and confirms our predictions.
The usage of a tensor product structure to consider
multidimensional problems does not lead to more 
information or further observations for our model problem as the one-dimensional setting.
%
This is the reason why we 
limit ourself to the one-dimensional case.
We do not only apply Gauß-Lobatto nodes, but also employ
a Gauß-Legendre basis. Our numerical simulations confirm our observation from remark \ref{remark_epsilon} 
that the error term  $\tilde{\epsilon}^k_2$ 
may have a positive effect on the numerical scheme, and we get more accurate solutions
by using a Gauß-Legendre basis.
Also the influence of the different numerical fluxes is less important than in the Gauß-Lobatto 
case, especially
if we chose the SBP-FR methods with $\mat{C}=\mat{M}[^{-1}] \mat{R}[^T] \mat{B}$ for the space 
discretization. \\
Simultaneously, results about the error behavior of several other correction terms $\mat{C}$ 
will be given.
The numerical schemes of table \ref{ta:correction_terms} and multiples of
them will be considered. 
We present several examples which justify our observations, but also show some limitations 
of our results. \\
%
%
We use an upwind flux (dotted lines) and  central flux (straight lines) at the interior 
element interfaces\footnote{We apply always
an upwind flux at the physical boundaries.}.\\
For time integration we use in all numerical examples a SSPRK(3,3) where the time step is chosen in 
such way  that the time integration error is negligible.
All elements are of uniform size. 
 \begin{figure}[!htp]
\centering 
  \begin{subfigure}[b]{0.45\textwidth}
    \includegraphics[width=\textwidth]{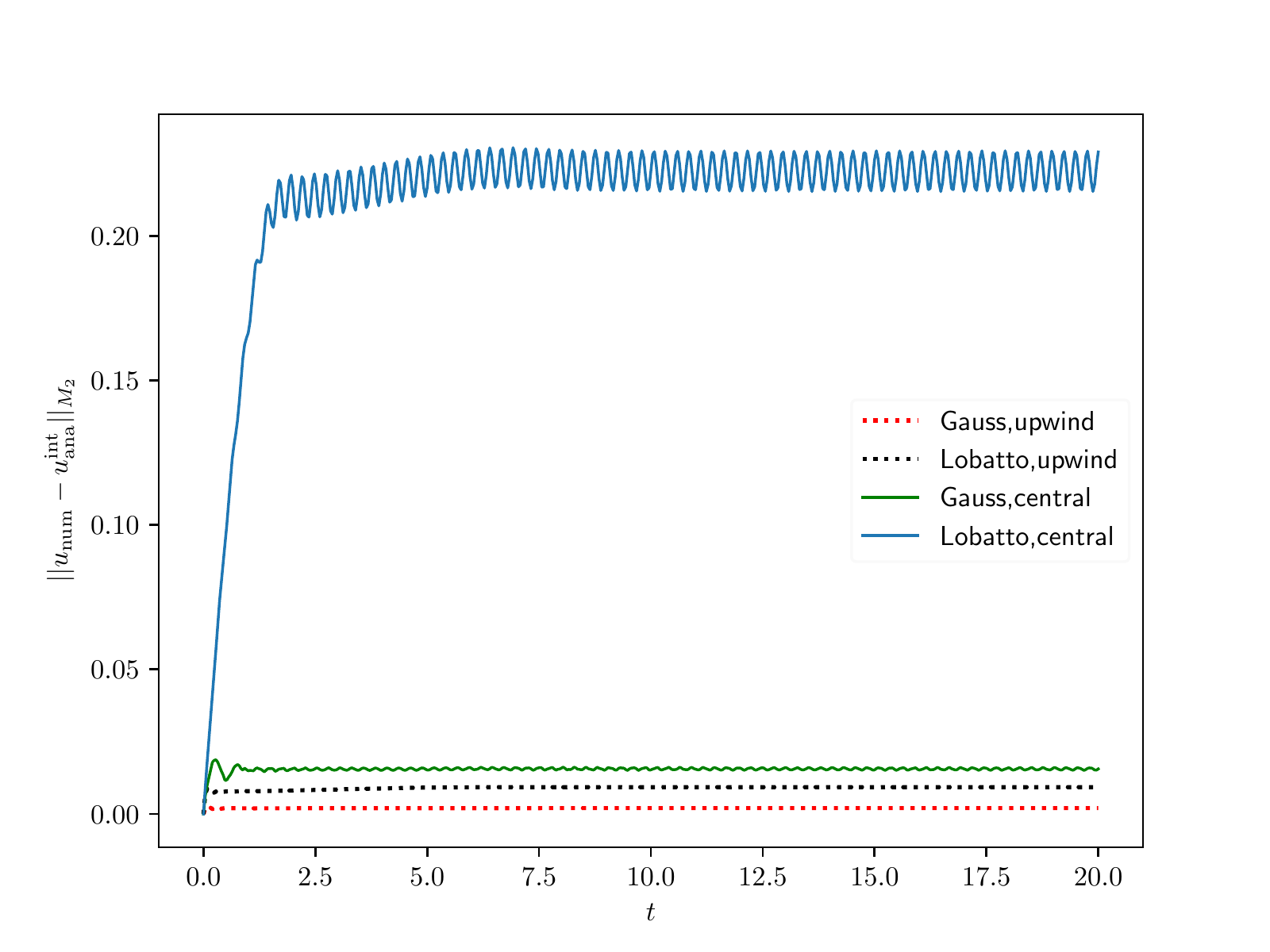}
    \caption{$N=4, K=30, t=20$}
  \end{subfigure}%
  ~
  \begin{subfigure}[b]{0.45\textwidth}
    \includegraphics[width=\textwidth]{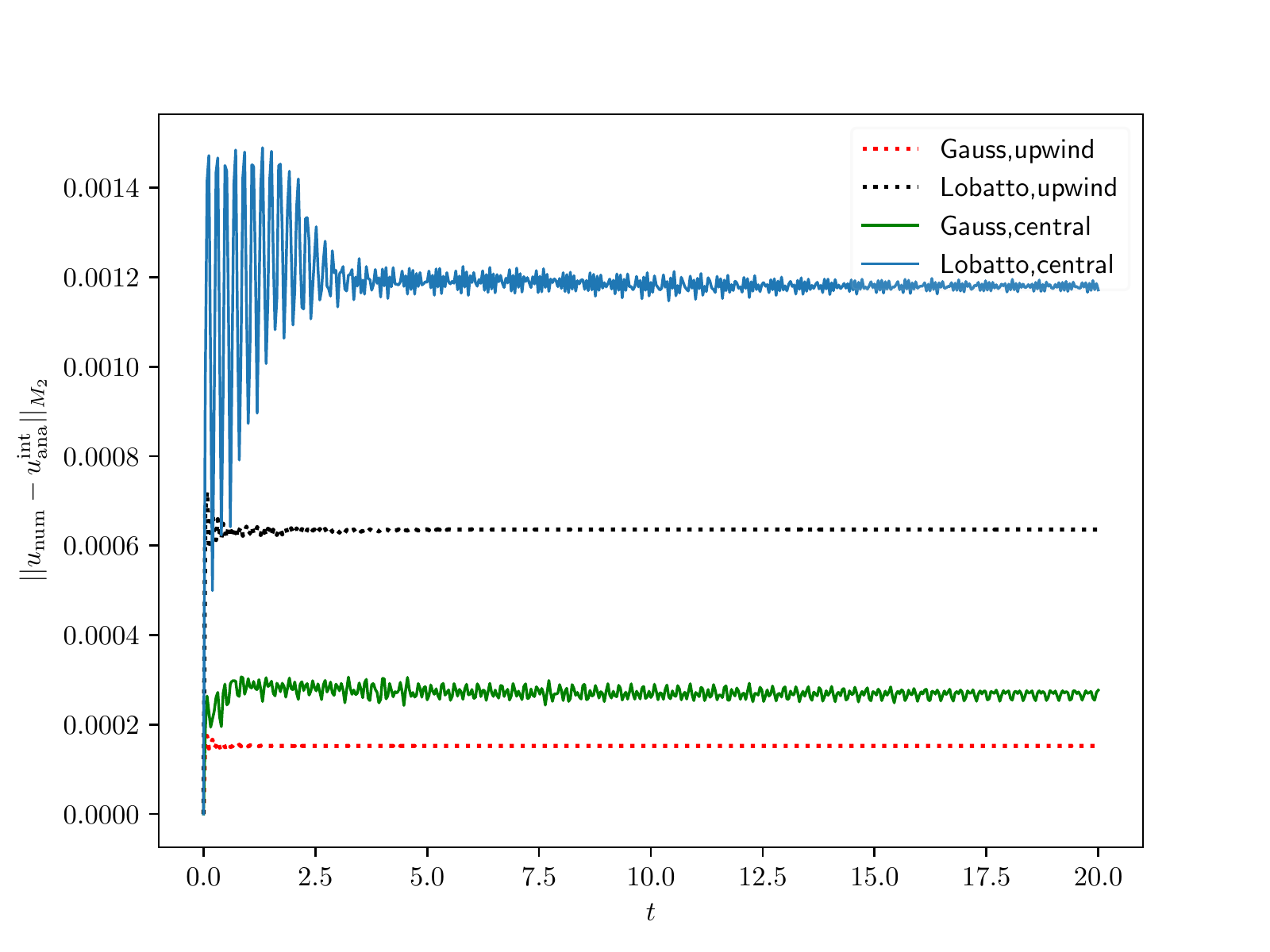}
    \caption{$N=4, K=50, t=20$}
  \end{subfigure}%
   \\~
   \begin{subfigure}[b]{0.45\textwidth}
    \includegraphics[width=\textwidth]{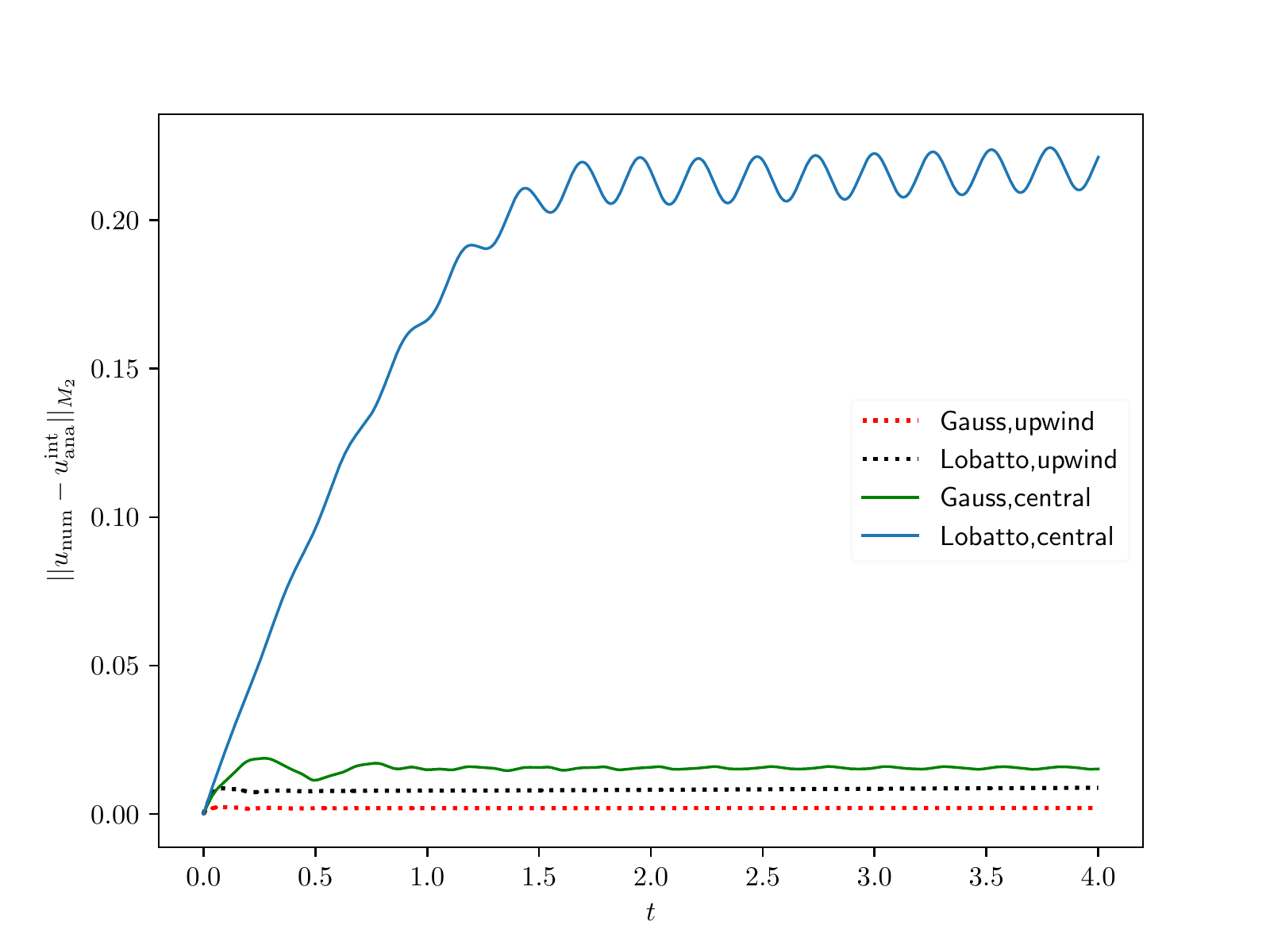}
    \caption{$N=4, K=30, t=4$}
  \end{subfigure}%
  ~
 \begin{subfigure}[b]{0.45\textwidth}
    \includegraphics[width=\textwidth]{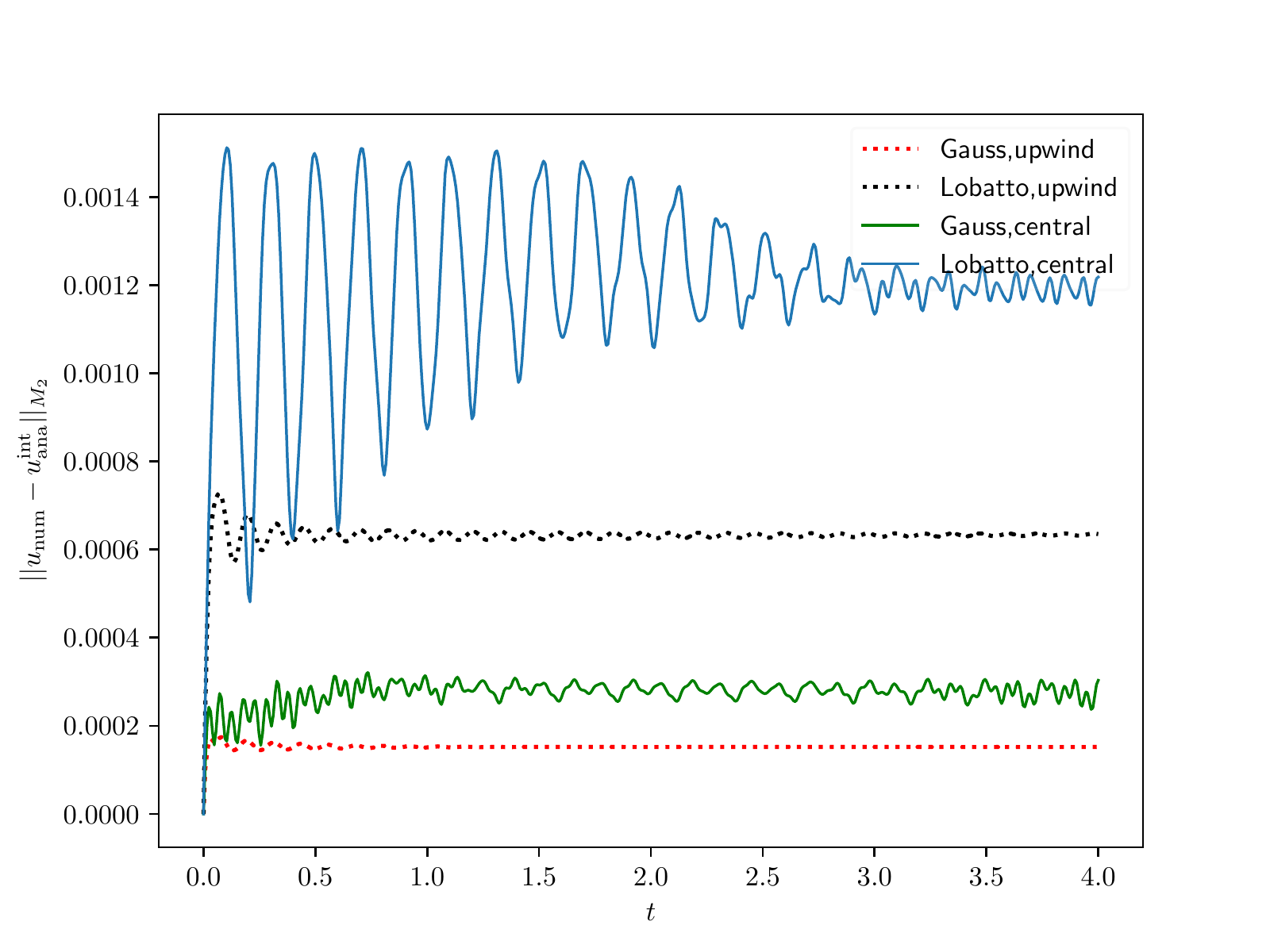}
    \caption{$N=4, K=50,t=4$}
  \end{subfigure}%
  \caption{Error  as a function in time. The dashed lines are always the calculation with 
  the upwind flux. Left side less elements than right. (c) and (d) early time behavior.  }
  \label{fig:sin_test}
\end{figure}

\subsection{Error behavior }

\subsubsection*{Sine-Testcase}
We start our numerical section with the 
same example as in \cite[p. 32]{kopriva2017error}.
We analyze the error behavior for $L=2\pi$ and the
initial condition $u_0=\sin (12(x-0.1))$, with the boundary
condition $g(t)$ chosen so that the exact solution is $u(x,t)=\sin (12(x-t-0.1))$.
In figure \ref{fig:sin_test} 
we illustrate the discrete errors over time for different number of elements
with a fourth order polynomial approximation.
The errors are always bounded in time for all combinations (P1)
(upwind / central flux and Gauß-Lobatto / Gauß-Legendre basis).
We  realize that the upwind flux errors reach 
its asymptotic values faster than the central flux errors in all cases (P4).
Simultaneously, the error bounds for the central flux are larger than for the upwind flux (P2). 
These results have already been seen in  
\cite{kopriva2017error} (in case of Gauß-Lobatto nodes), together with the fact that the central flux errors
are noisier than the upwind flux in all observations for all of the meshes and polynomial orders.
%
Here, we make the following two \emph{new} observations which confirms our (new) predictions.\\
The error bounds using Gauß-Lobatto points are larger 
than in the Gauß-Legendre case (P5), and secondly the influence
of the different numerical fluxes is less important than in 
the Gauß-Lobatto case (P6). 
Besides the accuracy properties of the different bases,
the error term $\tilde{\epsilon_2}$ has also a positive effect on 
the total error, and we get a more accurate solution in this case, see remark \ref{remark_epsilon}.
Gauß-Legendre nodes do not include the points at the element interfaces.
The additional dissipation comes from his fact and so the influence of the dissipation from 
the upwind flux is less important compared to the Gauß-Lobatto case. 
If we increase the order of approximation, the error bounds of the 
different combinations should coincide.
The figure \ref{fig:sin_test_polynome} justifies this prediction  (P3).  
 \begin{figure}[!htp]
\centering 
  \begin{subfigure}[b]{0.45\textwidth}
    \includegraphics[width=\textwidth]{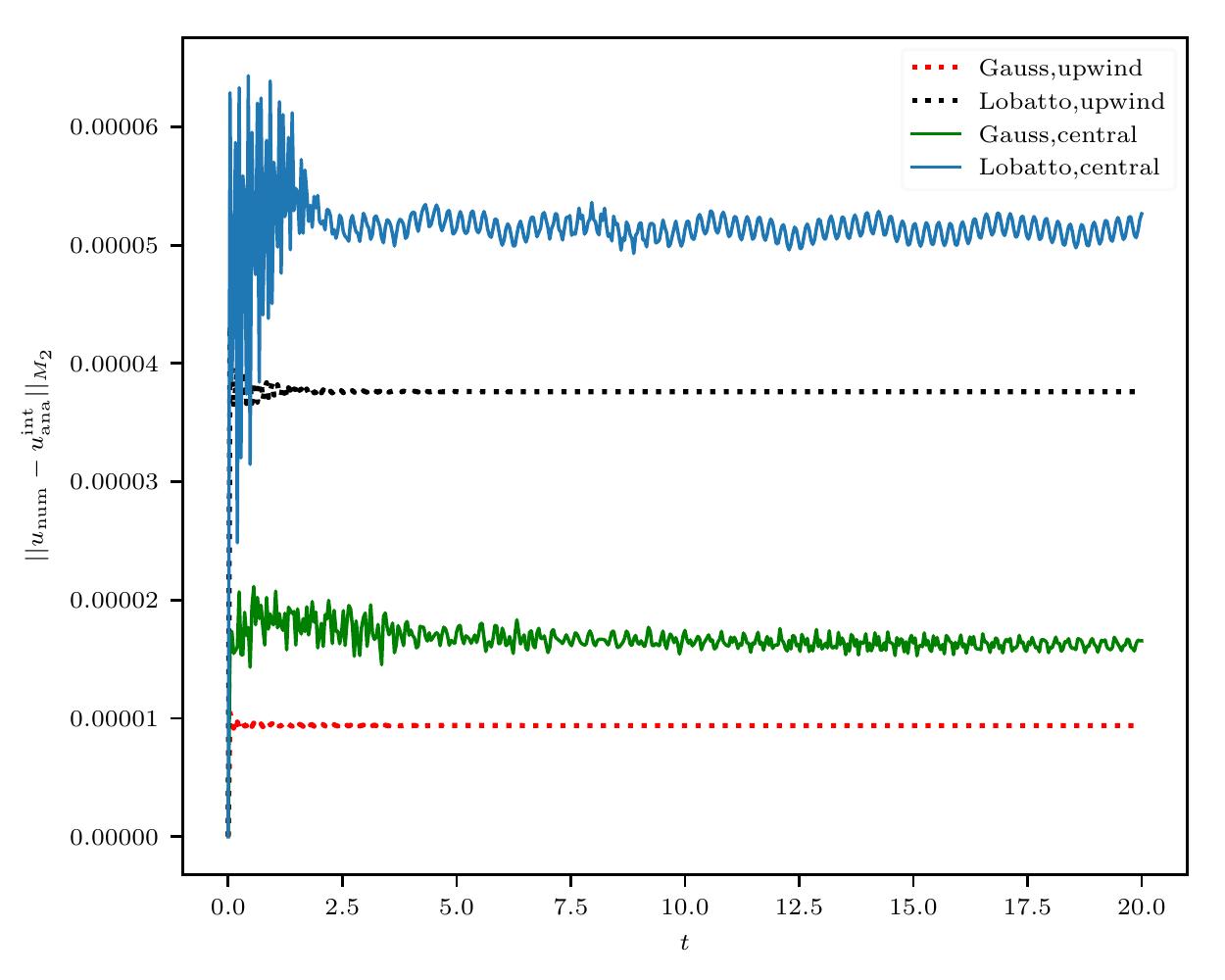}
    \caption{$N=5$}
  \end{subfigure}%
  ~
  \begin{subfigure}[b]{0.45\textwidth}
    \includegraphics[width=\textwidth]{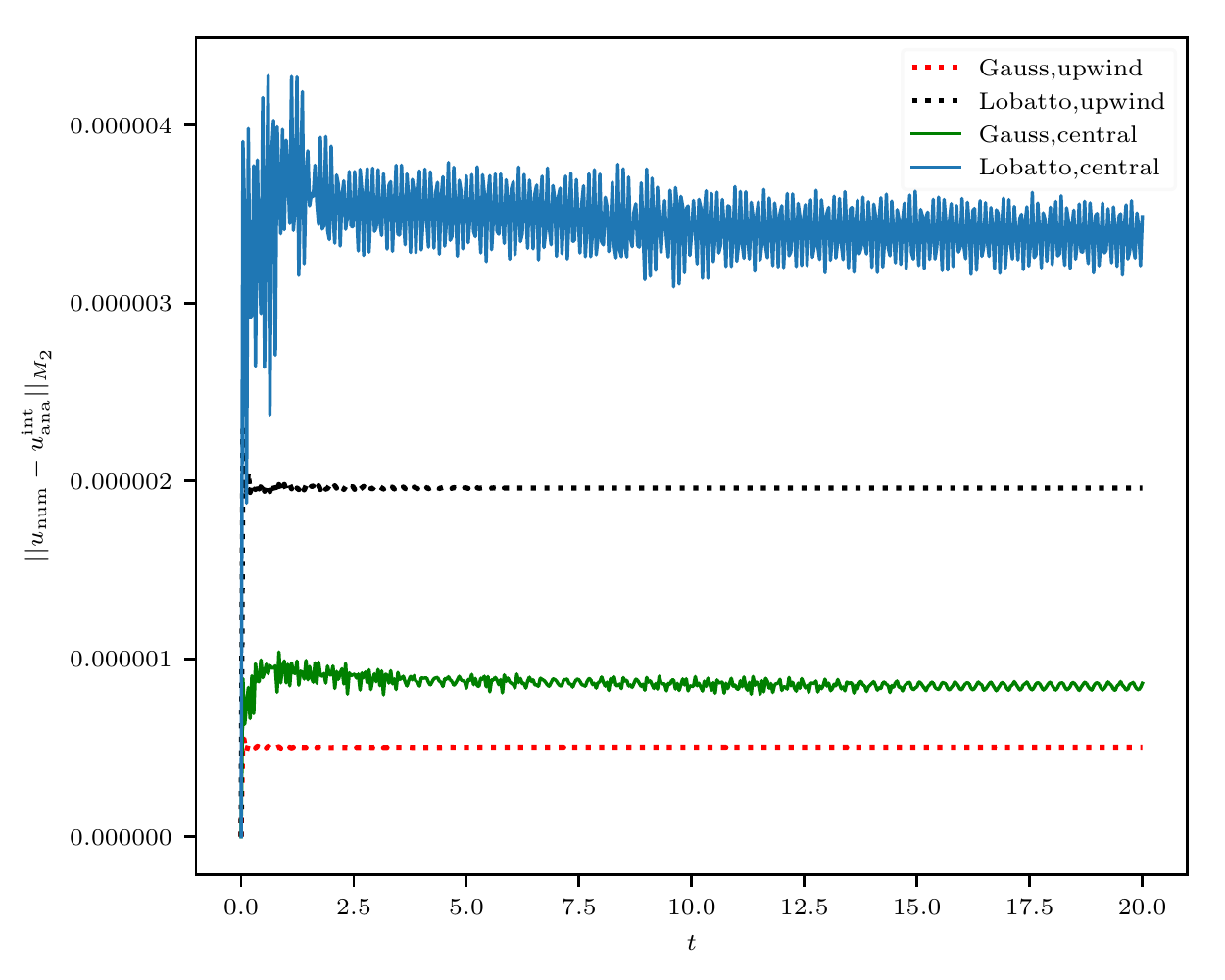}
    \caption{$N=6$}
  \end{subfigure}%
  \\~ 
   \begin{subfigure}[b]{0.45\textwidth}
    \includegraphics[width=\textwidth]{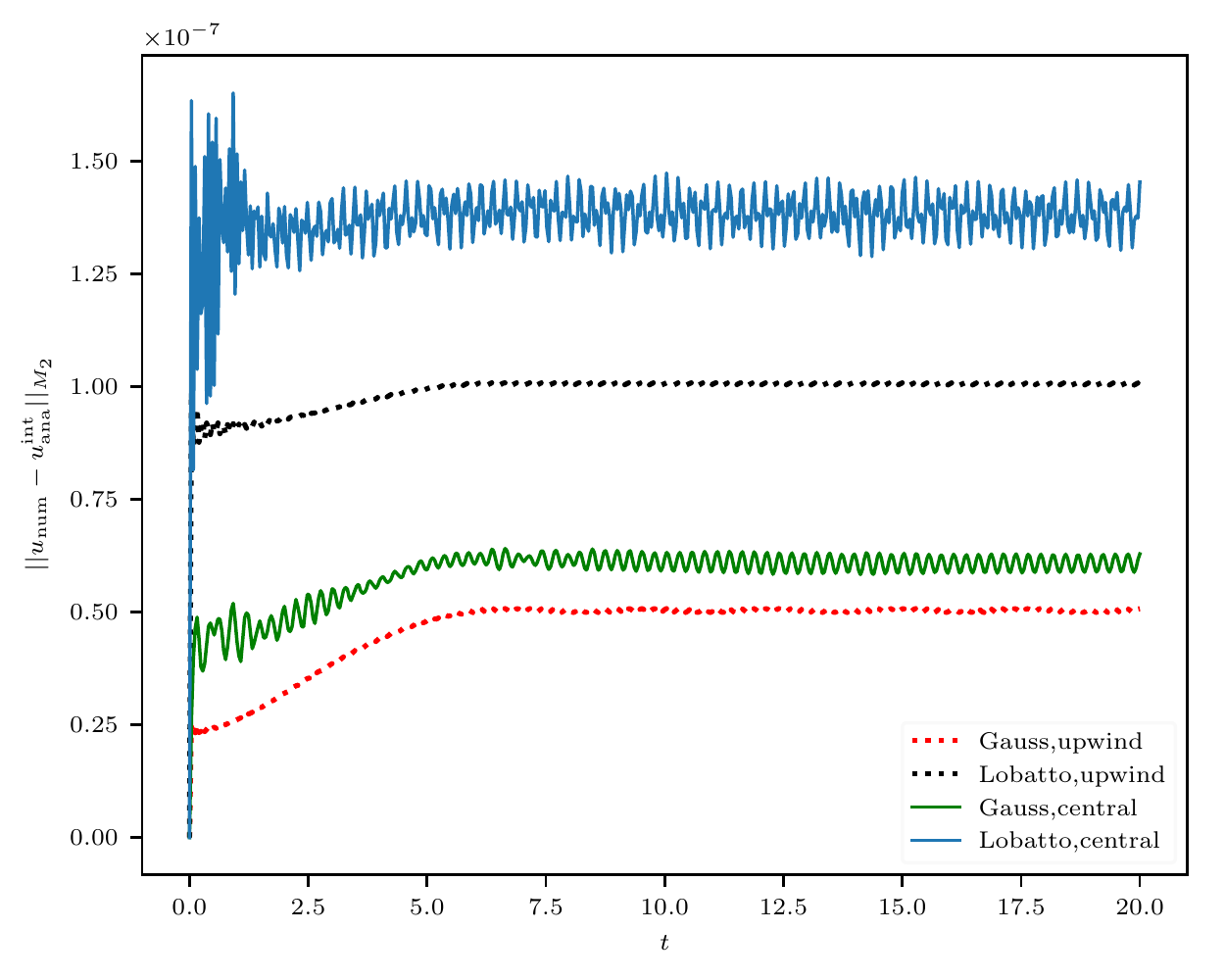}
    \caption{$N=7$}
  \end{subfigure}%
  ~
 \begin{subfigure}[b]{0.45\textwidth}
    \includegraphics[width=\textwidth]{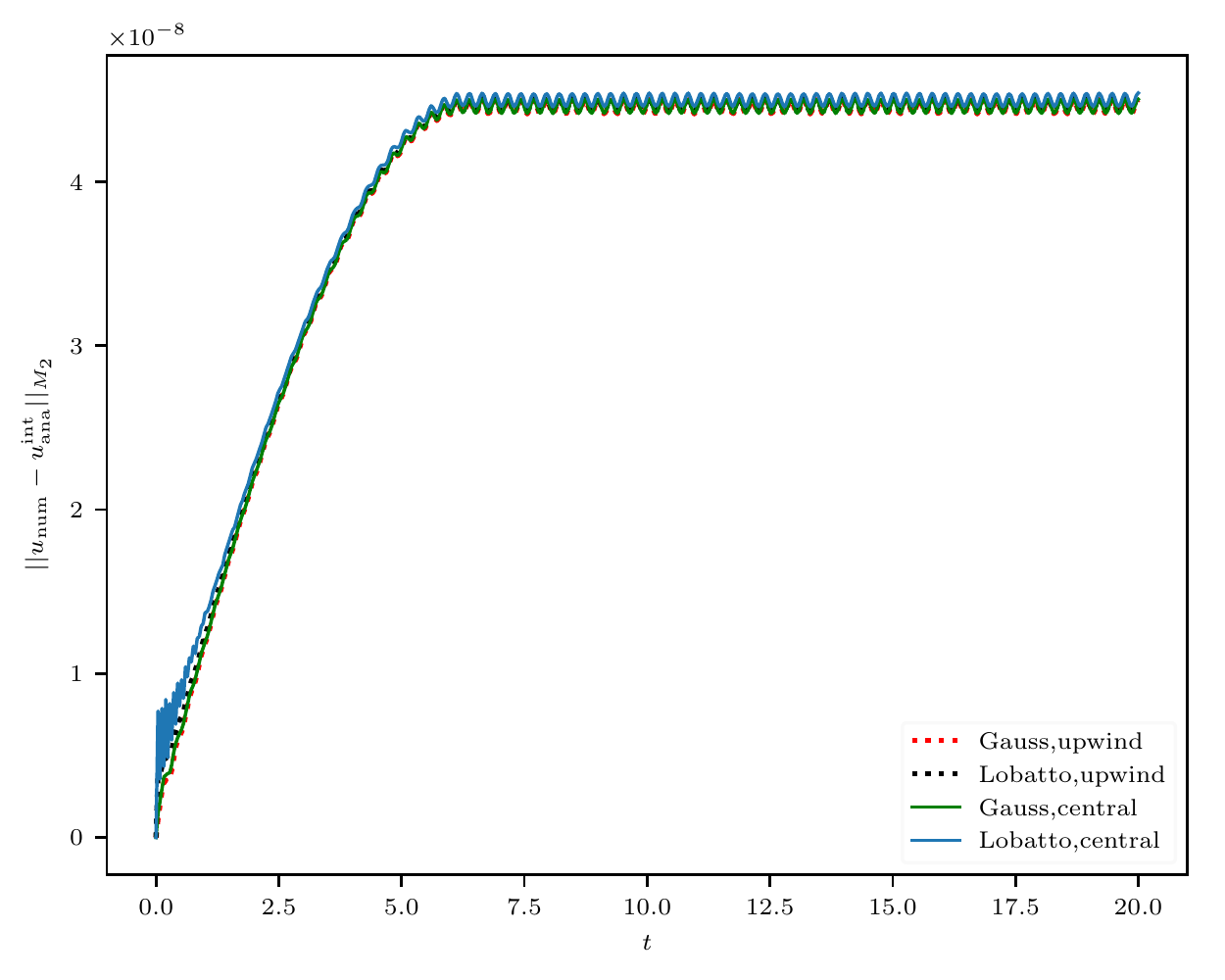}
    \caption{$N=8$ }
  \end{subfigure}%
  \caption{Error as a function in time. The dashed lines are always the calculation with 
  the upwind flux. $K=50$ and $t=20$.  }
  \label{fig:sin_test_polynome}
\end{figure}
\\Last, but not least,  we also study the convergence speed and  observe spectral accuracy in all cases, see figure \ref{fig:error}.
 \begin{figure}[!htp]
\centering 
    \includegraphics[width=0.5\textwidth]{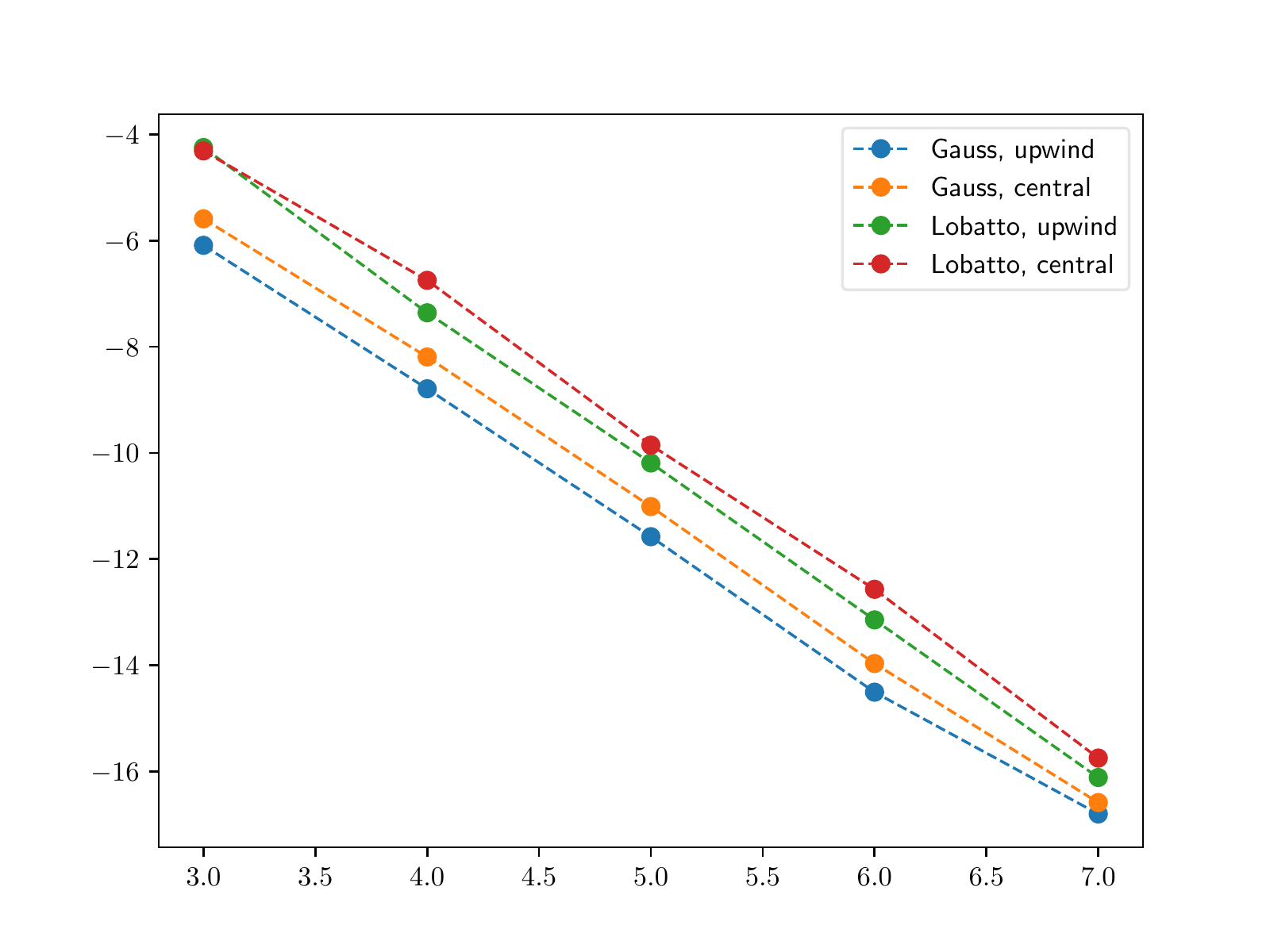}
  \caption{Convergences  in time asymptotic errors (last value) as  functions of $N$ for $K=50$. }
  \label{fig:error}
\end{figure}
This suggests that the approximation errors in $\Ep(t)$ decay faster
than $\frac{1}{\delta_0}$ grows,
since with inequality \eqref{eq:Errorbound3} one predicts that the time asymptotic 
error is bounded by $\Ep(t)/\delta_0$.
This matches also with the investigation in \cite{kopriva2017error}.\\
We are not only considering the FR scheme which is equivalent to the DG framework. We also investigate 
the schemes of Huynh and the SD methods.\\
The correction terms are written in table \ref{ta:correction_terms} and we demonstrate the error 
behaviors in figure \ref{fig:sin_test_SD} using the $||\cdot||_{\K_M}$-norm and  the global norm \eqref{eq:globalnorm}. 
 \begin{figure}[!htp]
\centering 
  \begin{subfigure}[b]{0.45\textwidth}
    \includegraphics[width=\textwidth]{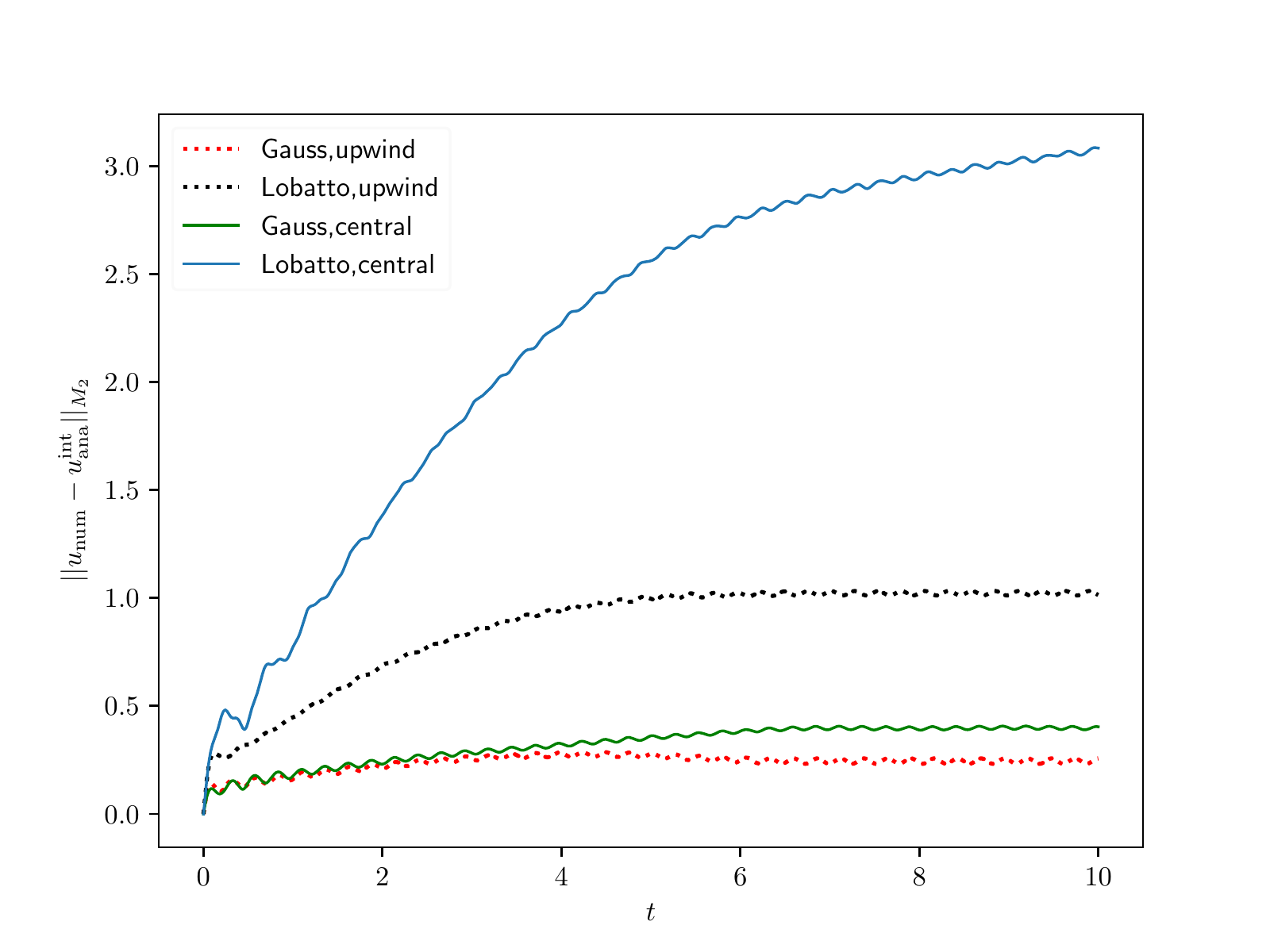}
    \caption{$N=3, K=20, t=10$}
  \end{subfigure}%
  ~
  \begin{subfigure}[b]{0.45\textwidth}
    \includegraphics[width=\textwidth]{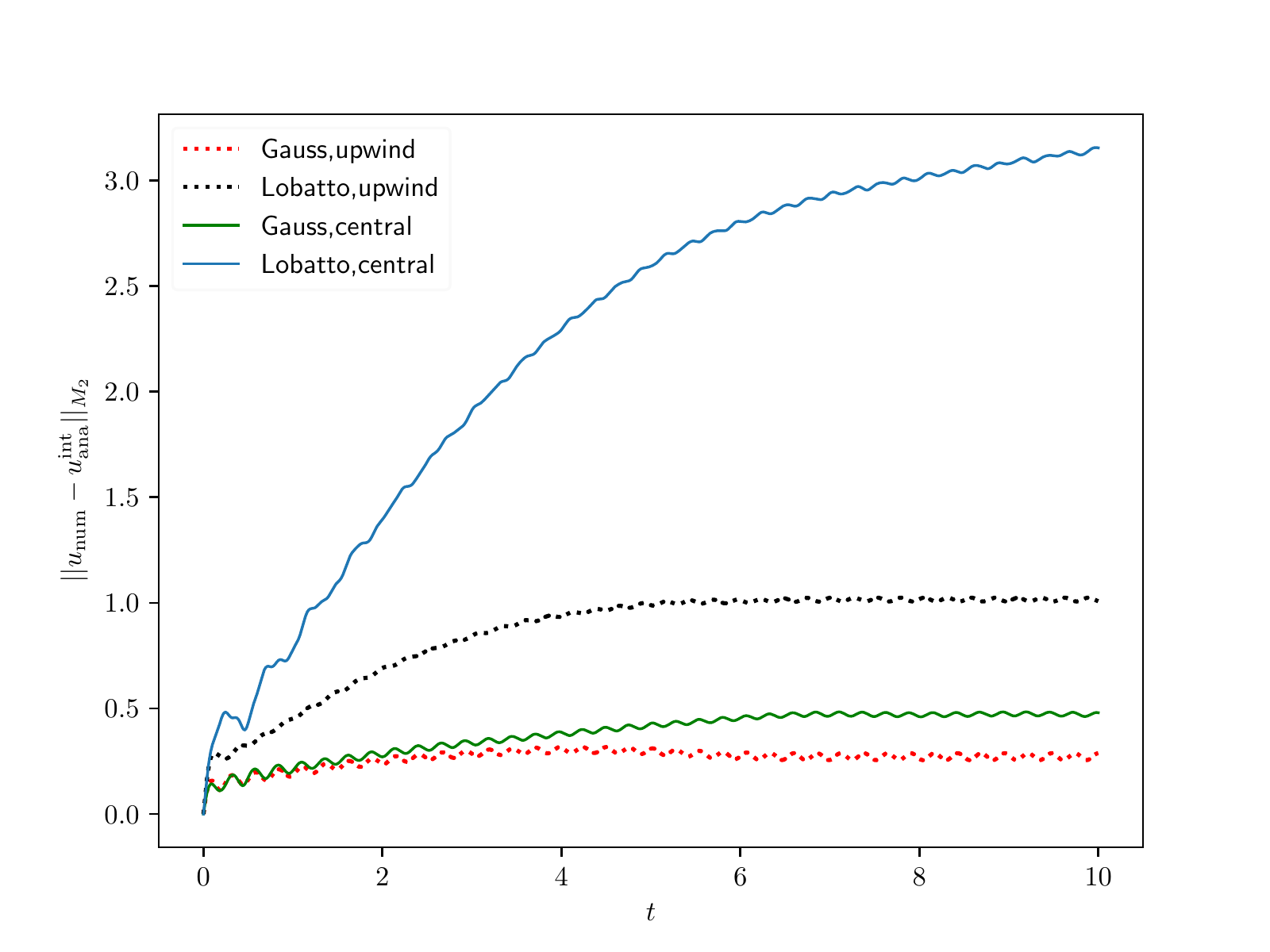}
    \caption{$N=3, K=20, t=10$}
  \end{subfigure}%
   \\~
   \begin{subfigure}[b]{0.45\textwidth}
    \includegraphics[width=\textwidth]{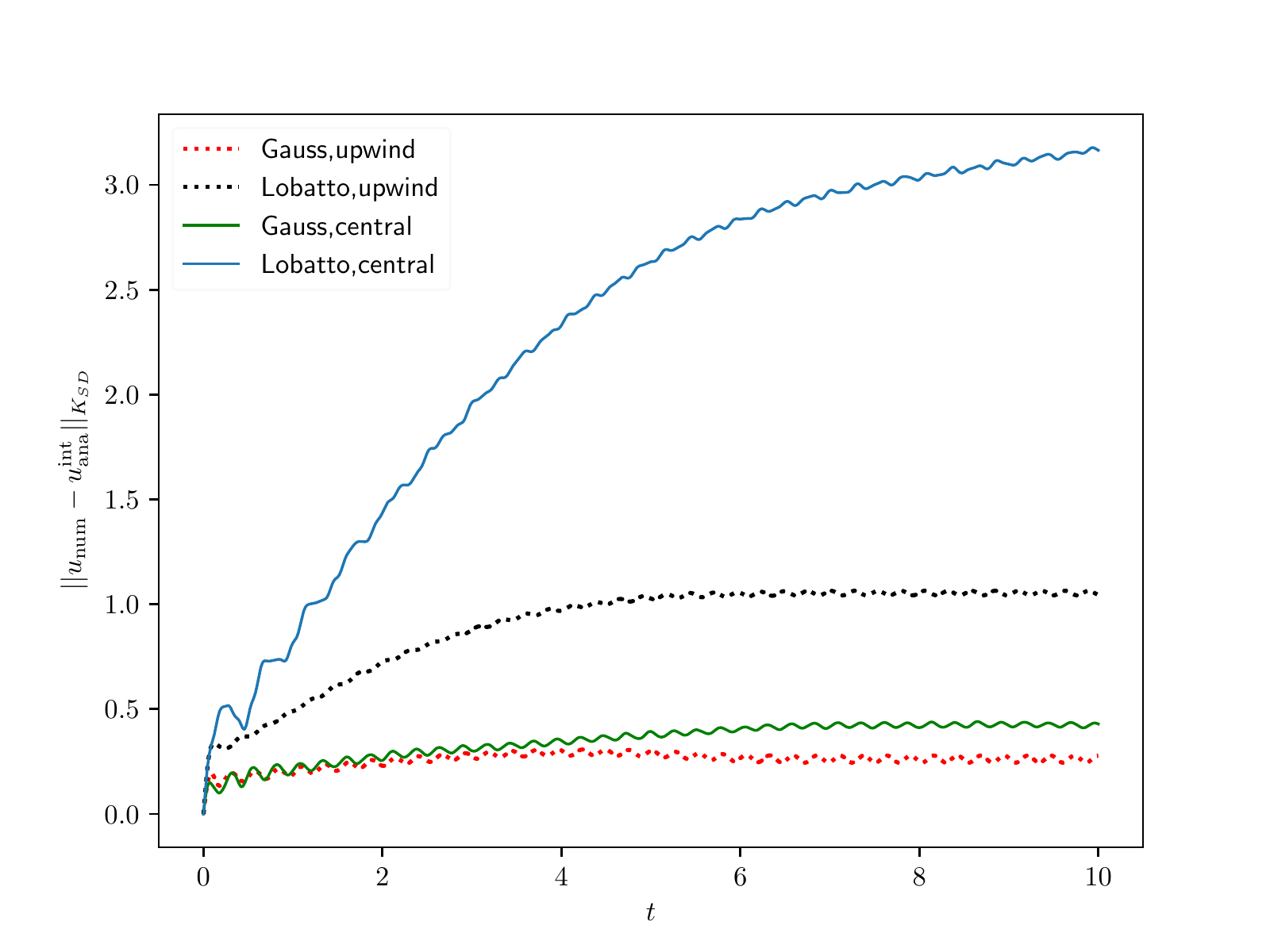}
    \caption{$N=3, K=20, t=10$}
  \end{subfigure}%
  ~
 \begin{subfigure}[b]{0.45\textwidth}
    \includegraphics[width=\textwidth]{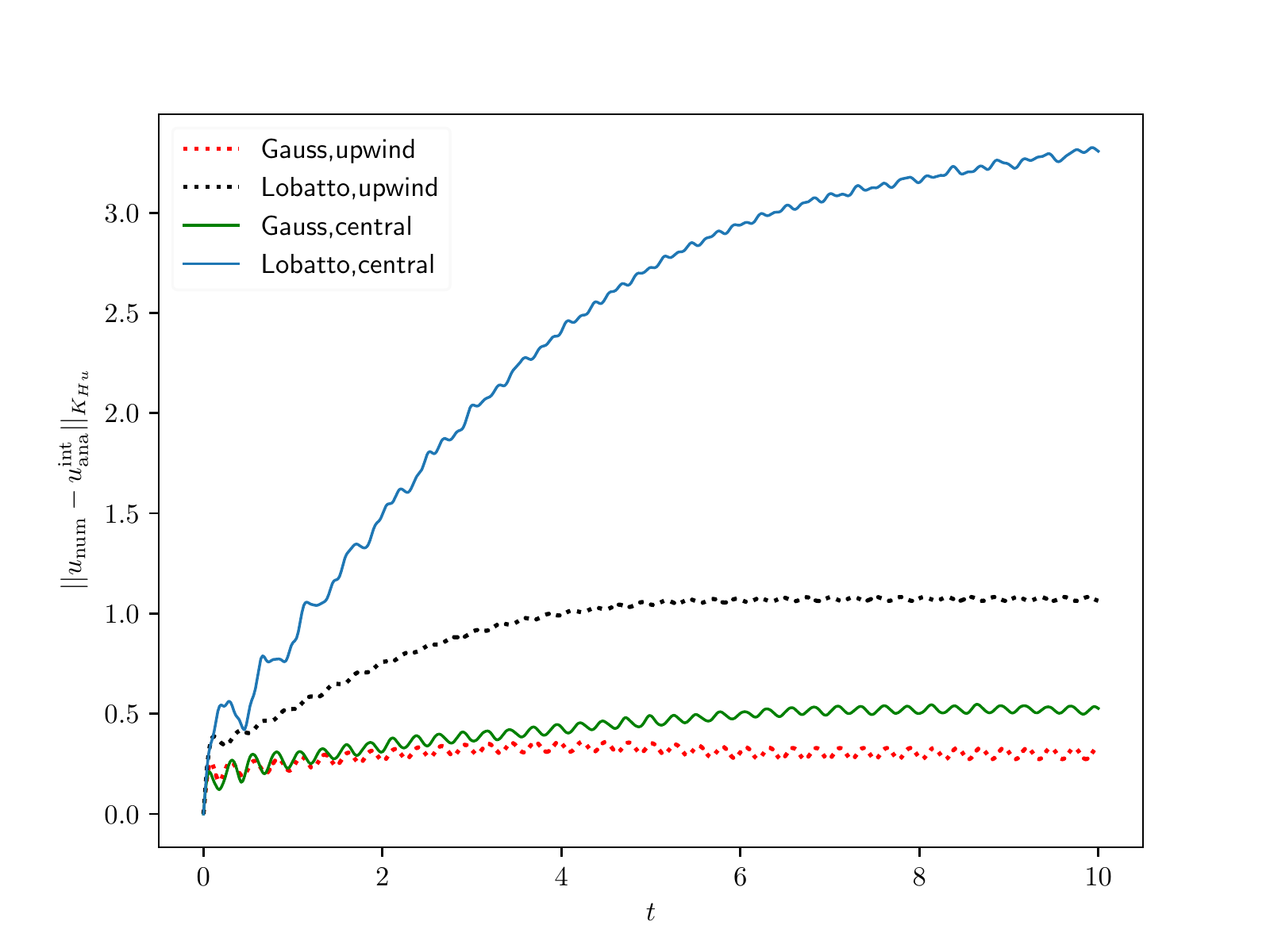}
    \caption{$N=3, K=20,t=10$}
  \end{subfigure}%
  \caption{Error as a function in time. The dashed lines are always the calculation with 
  the upwind flux. Right side SD scheme and on the left side Huynh. (c) and (d) respective $\K_M$-norm.  }
  \label{fig:sin_test_SD}
\end{figure}
We realize that the error is always bounded in all cases (P1) and, indeed,
that the influence of the numerical flux is less important comparing to the chosen basis function 
((P5) and (P6)). 
Here, the usage of Gauß-Legendre nodes demonstrates their good advantage comparing to the usage
of Gauß-Lobatto nodes (P5). 
Nevertheless, we make also one more observation in this case. Different from the DG case, 
our errors show always some oscillations even using the upwind flux (P7). This can be seen 
in figure \ref{fig:Gauss} where only Gauß-Legendre nodes are considered. 
This is due to the fact that using another correction term as for the case $\kappa=0$, 
we get some \emph{over correction} and \emph{under correction}
at the boundaries through our correction functions.
It is not surprising that by using Gauß-Lobatto nodes and a central flux 
we get the worst simulation if the resolution is low order (P1-P7).
\begin{figure}[!htp]

    \includegraphics[width=0.45\textwidth]{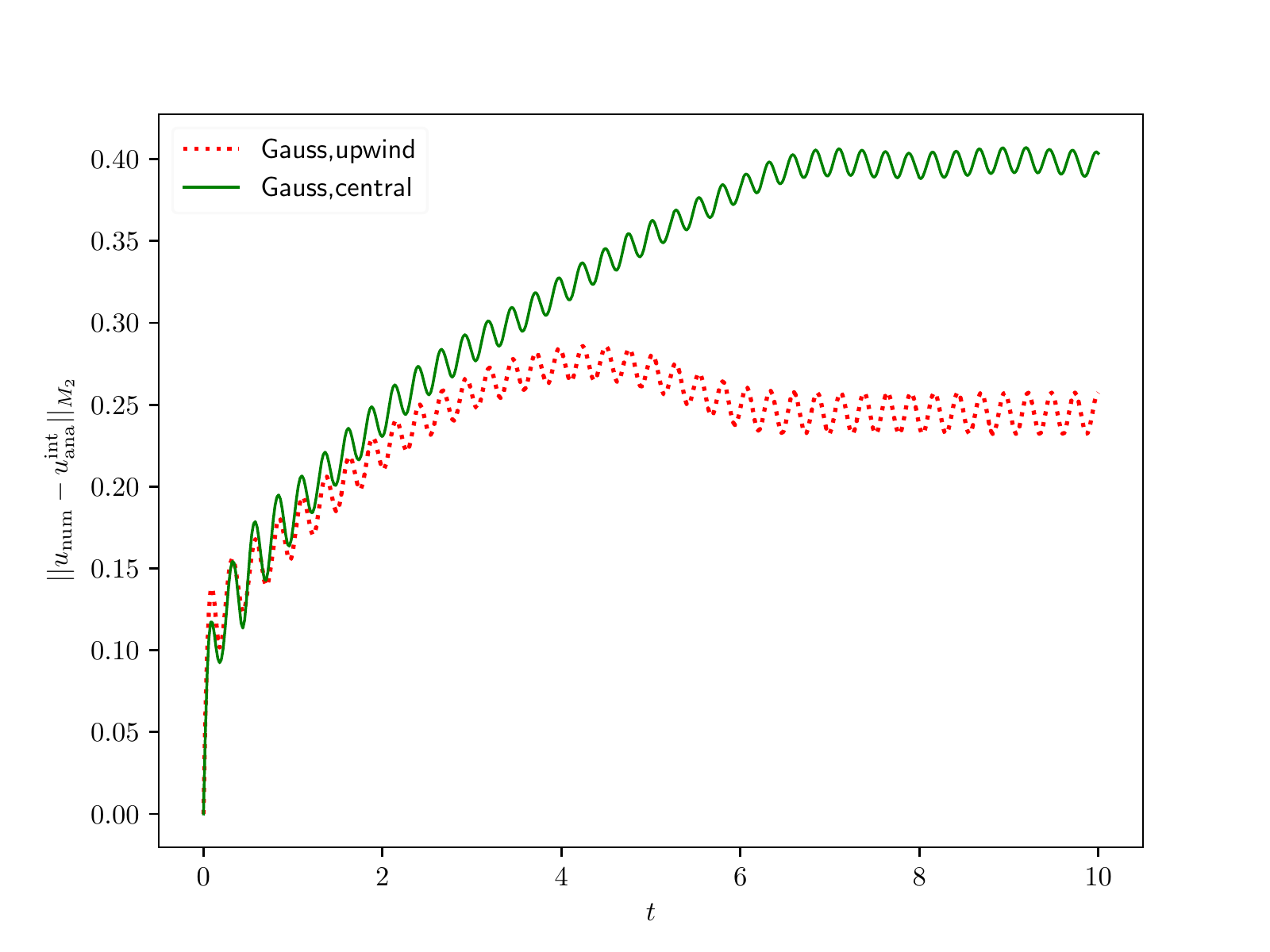}
    \includegraphics[width=0.45\textwidth]{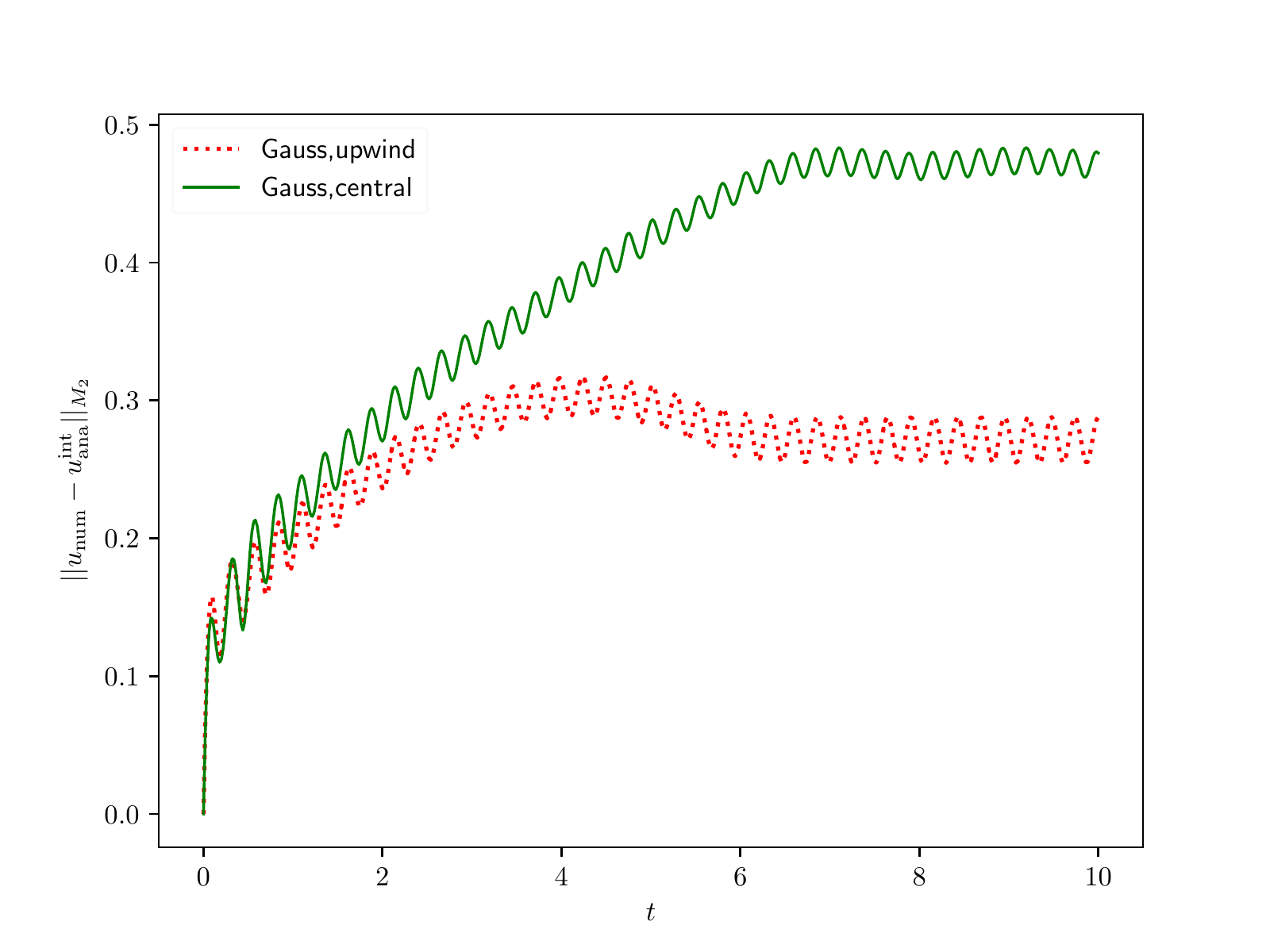}
  \caption{Error as a function in time only with Gauß-Legendre nodes. Left: SD, Right: Huynh}
  \label{fig:Gauss}
\end{figure}

If we decrease the correction terms more rapidly to zero, we obtain some intermediate 
schemes on the way to the DG case ($\kappa=0$). 
In figure \ref{fig:other}, we plot the SD method and 
Huynh scheme by dividing the correction terms with $2^{N-1}$.
We see that these simulations behave like our first test 
(figure \ref{fig:sin_test})
except with noisier behavior. We also realize that using  Gauß-Legendre nodes, 
the scheme demonstrates a higher amplitude in the 
oscillations and, 
therefore, the Gauß-Lobatto nodes seem better (more about this in section \ref{subsec:counter}). 
\begin{figure}[!htp]
    \includegraphics[width=0.45\textwidth]{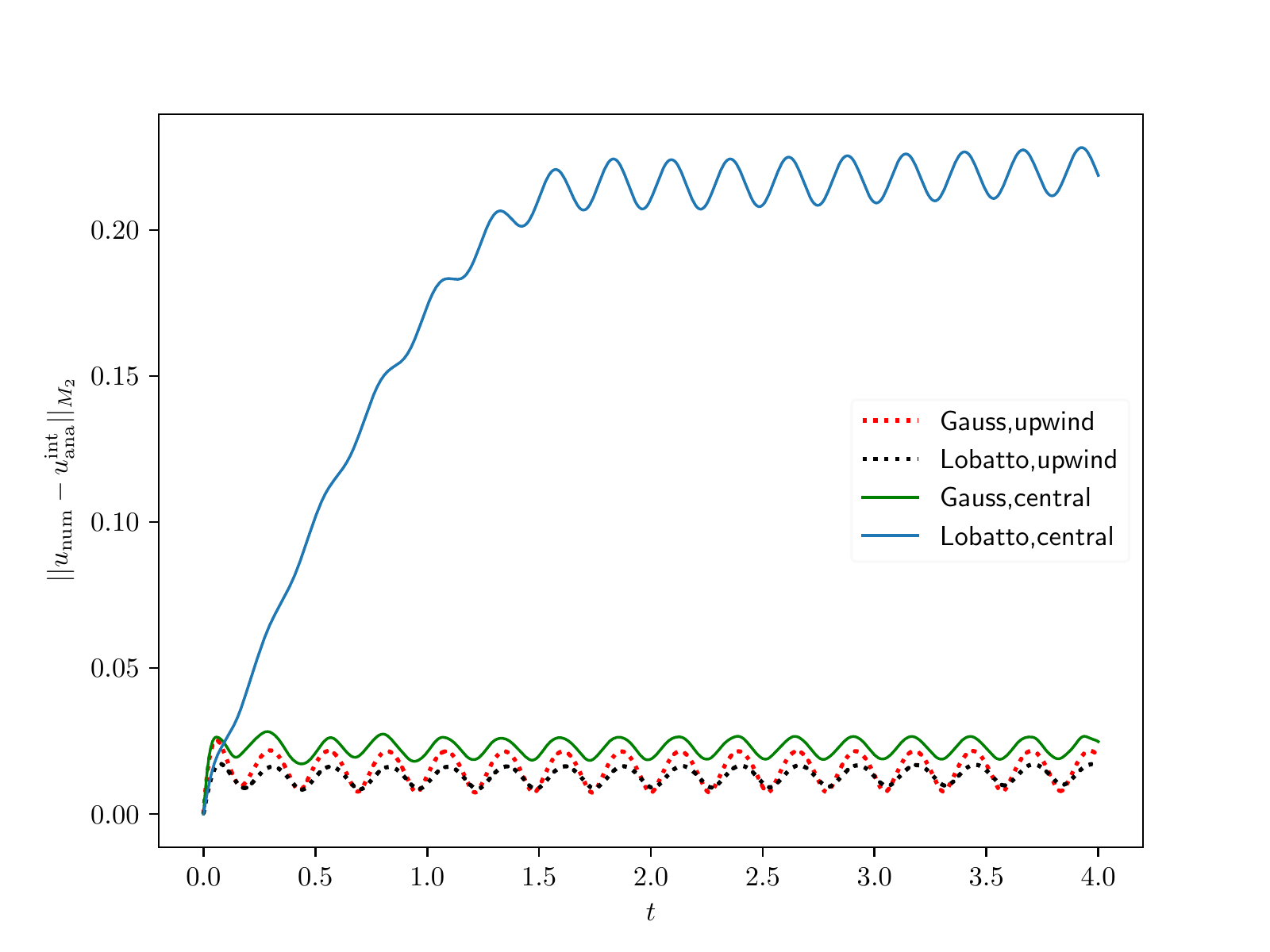}
    \includegraphics[width=0.45\textwidth]{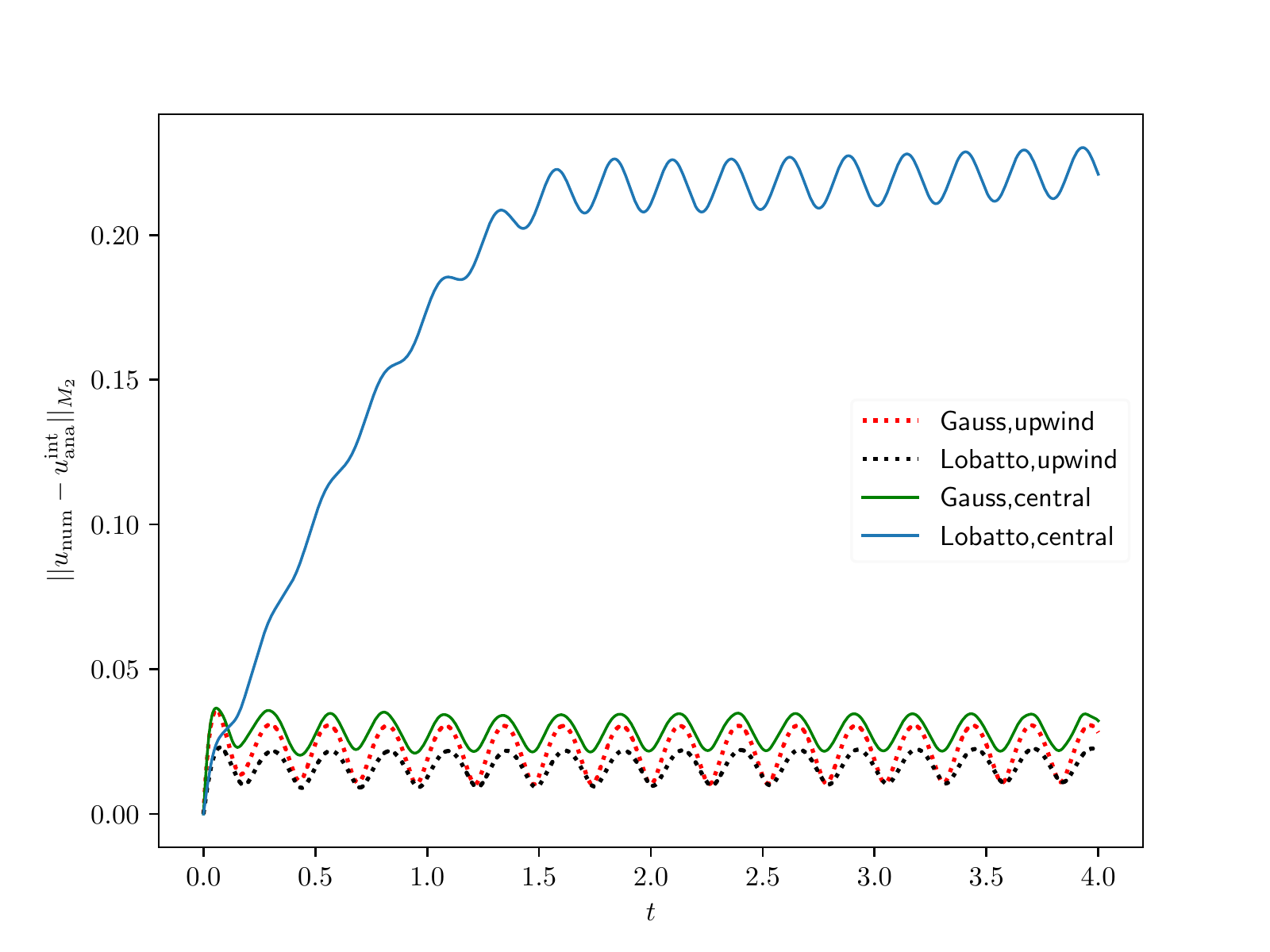}
  \caption{Error  as a function in time only with Gauß-Legendre nodes $K=30$ and $N=4$. \\
  Correction Term: $c_{SD}/2^{N-1}$ (left),  $c_{Hu}/2^{N-1}$ (right)}
  \label{fig:other}
\end{figure}

\subsubsection*{Cosine-Testcase}
As a second testcase, we investigate 
the error behavior for $L=2\pi$ and the initial condition $u_0=\cos (12(x-0.1))$, with the boundary
condition $g(t)$ chosen so that the exact solution is $u(x,t)=\cos(12(x-t-0.1))$.
With this testcase we want to strengthen our conclusions from before.

In figure \ref{fig:cos_test} we illustrate the discrete errors over time for different number of elements
with a fourth and sixth order polynomial approximation.
 \begin{figure}[!htp]
\centering 
  \begin{subfigure}[b]{0.45\textwidth}
    \includegraphics[width=\textwidth]{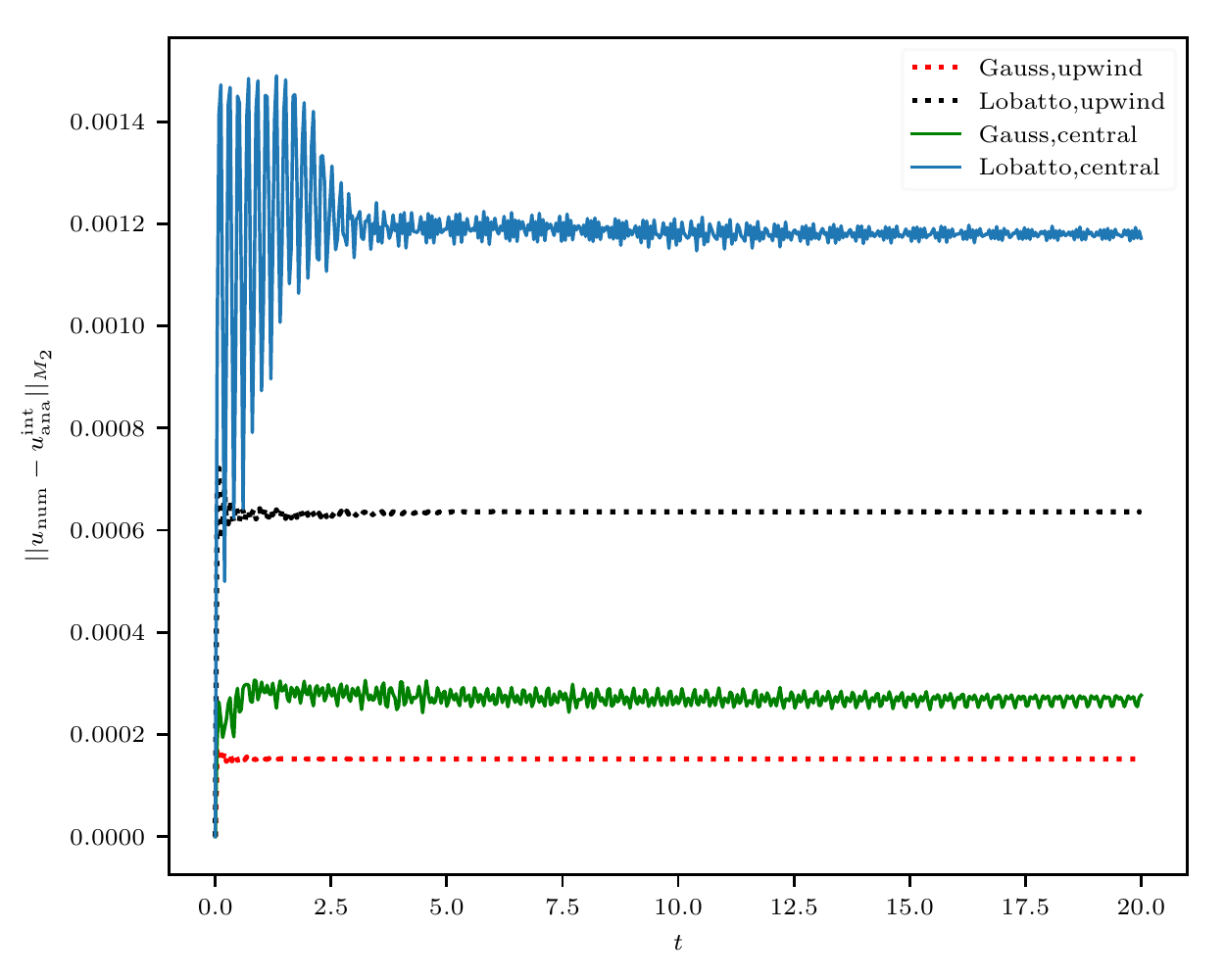}
    \caption{$N=4$}
  \end{subfigure}%
  ~
  \begin{subfigure}[b]{0.45\textwidth}
    \includegraphics[width=\textwidth]{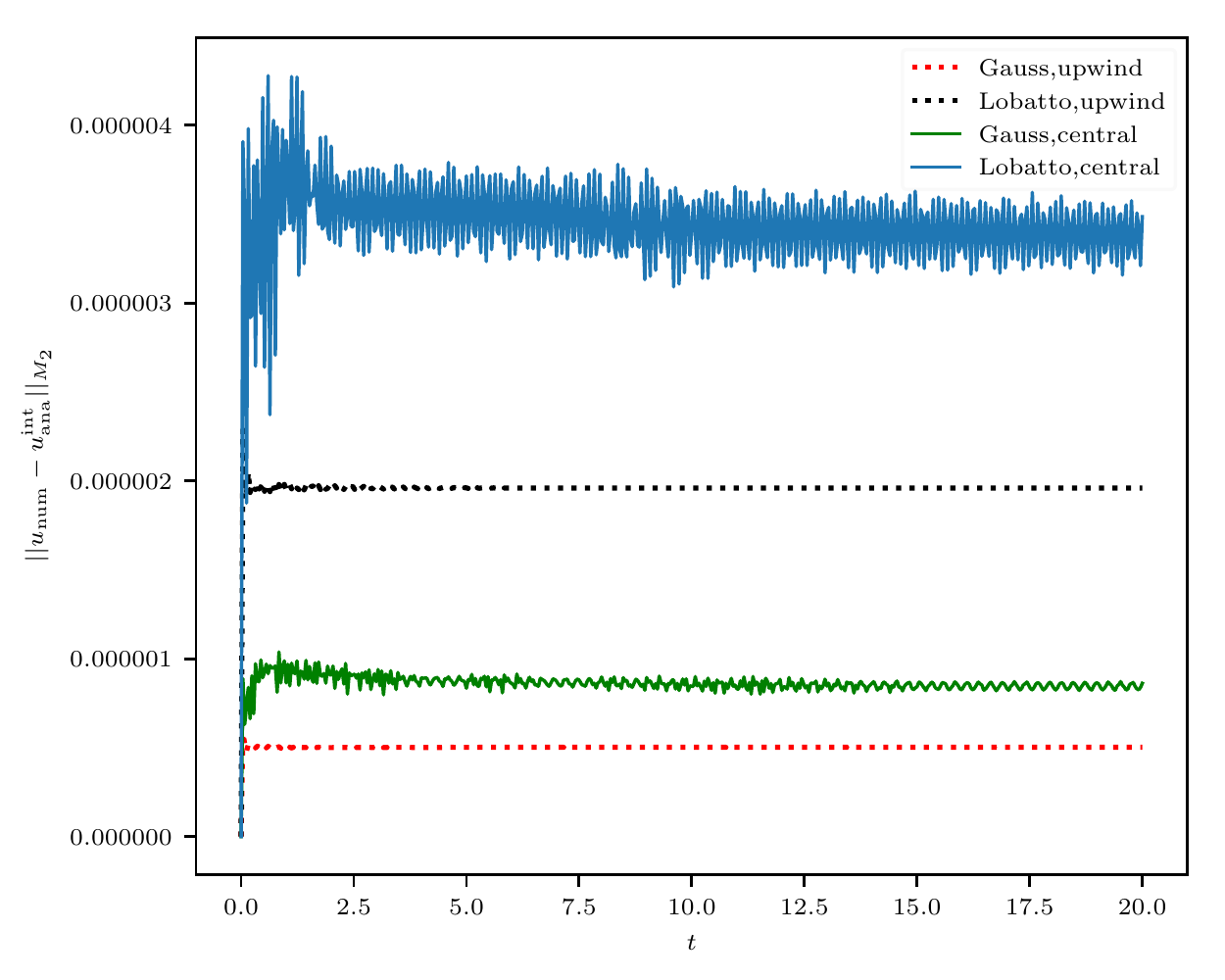}
    \caption{$N=6$}
  \end{subfigure}%
  \caption{Error  as a function in time with $50$-elements.}
  \label{fig:cos_test}
\end{figure}
We make equivalent observations like before and see that using Gauß-Legendre nodes
in our scheme yields  more accurate solutions 
than when using  Gauß-Lobatto nodes (P5). 
Also, the difference  between the upwind flux error and central flux error  is not so large (P6).
The same observations can be made by using the FR schemes from before.

\subsubsection*{$\tilde{\epsilon}_2$-Term}
Finally, we study the $\tilde{\epsilon}_2$ terms for the two testcases. 
In figure \ref{fig:epsilon_1},
the $\tilde{\epsilon}_2$-error is plotted over the time for different polynomial orders as in the $\sin$ testcase.
\begin{figure}[!htp]
    \includegraphics[width=0.45\textwidth]{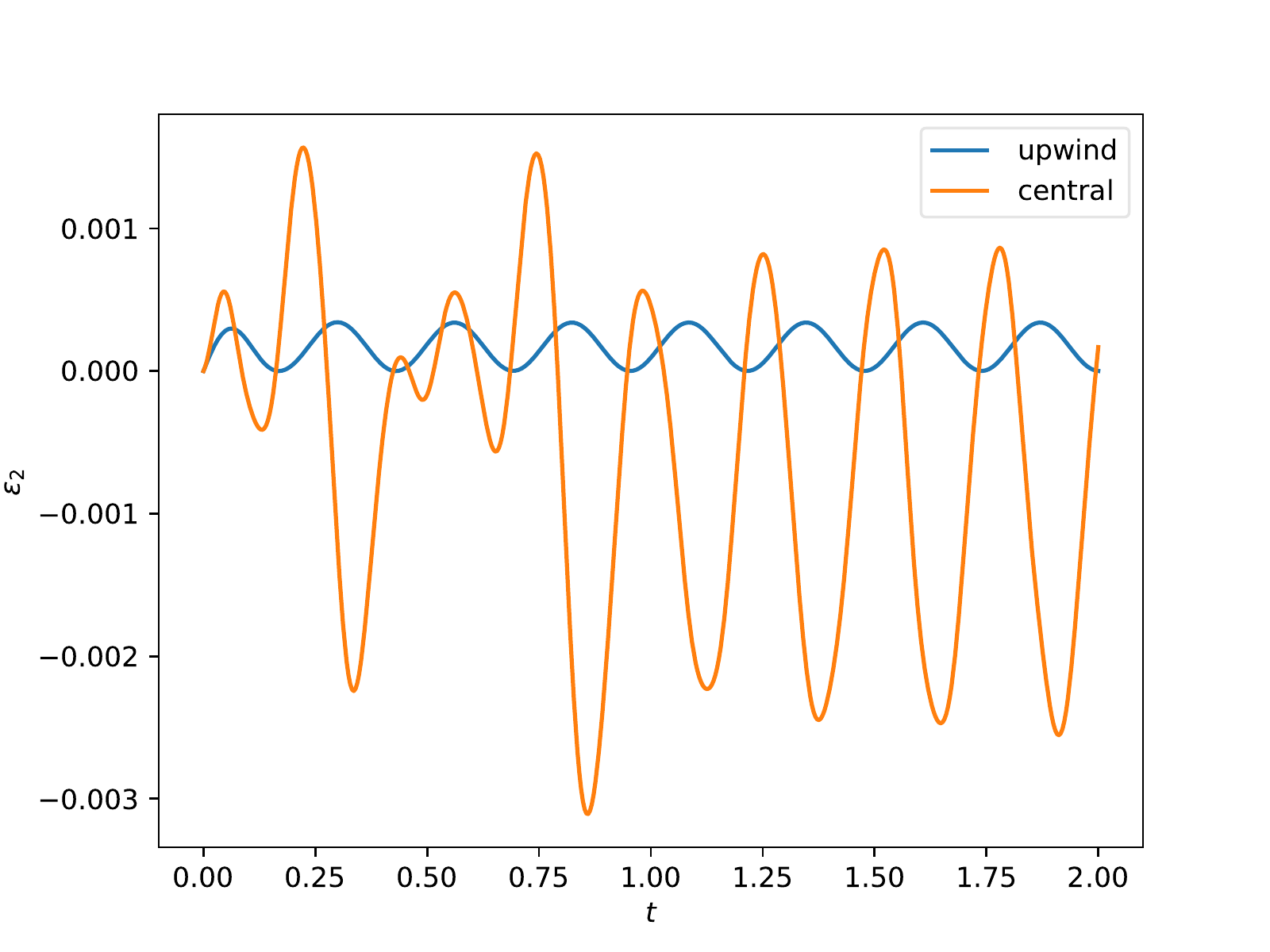}
    \includegraphics[width=0.45\textwidth]{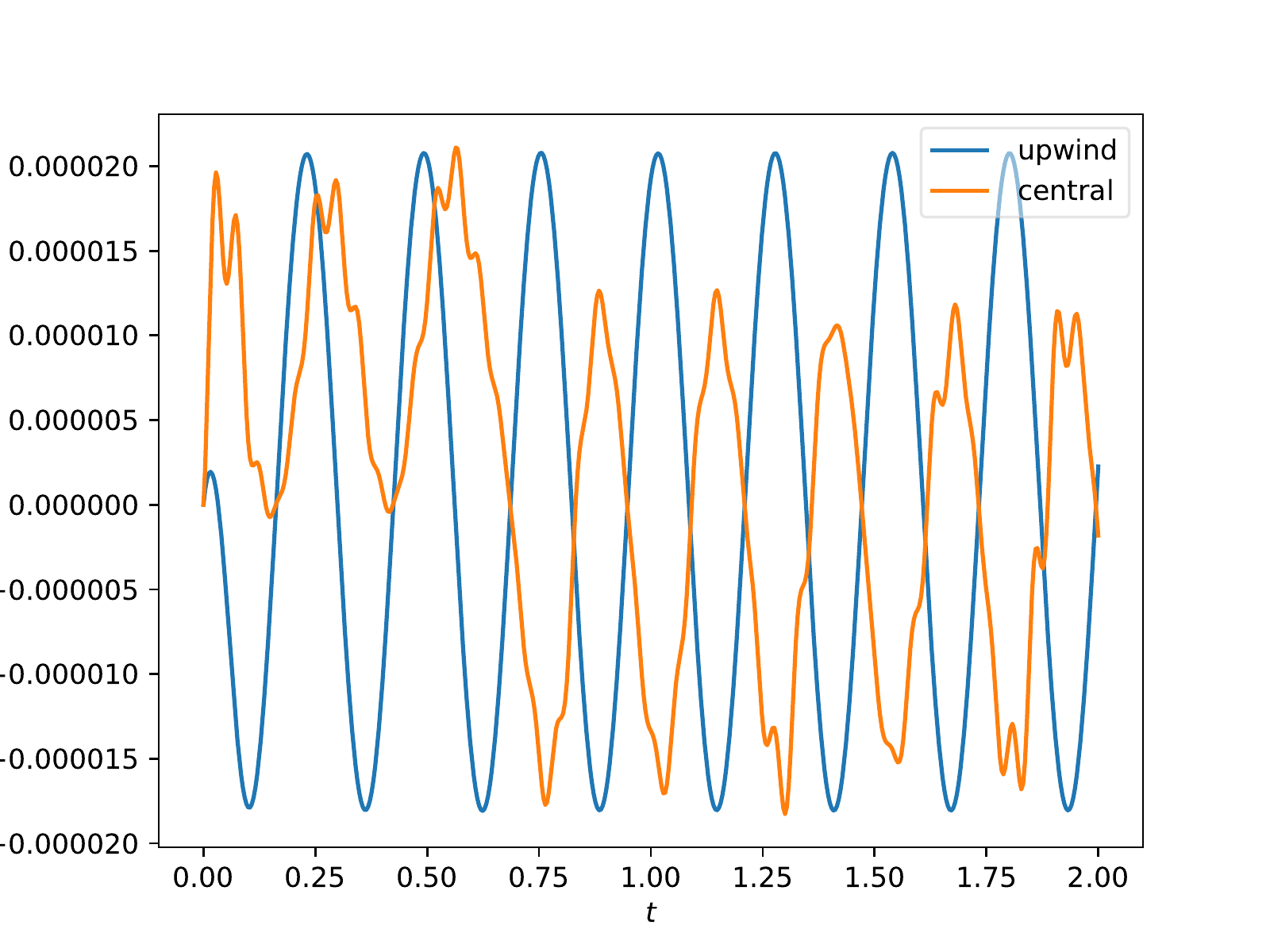}
  \caption{$\tilde{\epsilon}_2$ a function in time $K=30$ and $N=4, 5$. \\
 }
  \label{fig:epsilon_1}
\end{figure}
Both times the error starts positive for the lower order 
approximation ($N=4$), and the upwind flux stays strictly positive whereas the term 
using the central flux shows a higher amplitude in the oscillations
and turns also negative. Even the middle point of the amplitudes is below zero.  
Since the $\epsilon_1$- errors show these oscillations when applying the central flux, 
we also get them in the  $\tilde{\epsilon}_2$-terms. 
This is also not surprising because  $\epsilon_1$  
has a direct influence on $\tilde{\epsilon}_2$. 
Using a fifth order approximation, the amplitude of the error using the upwind flux 
shows a higher amount, but these are symmetrical around zero, whereas applying 
the central  flux the symmetric point  lies in the negative axis. 
Comparing  the amounts of the total error in  figure \ref{fig:sin_test}  with figure \ref{fig:epsilon_1}
demonstrates that the $\tilde{\epsilon}_2$ errors 
are significantly less, but may have an positive influence, especially if Gauß-Legendre nodes are used. 
In  figure \ref{fig:epsilon_2} we have an analogous behavior for
the cosine testcase. 
\begin{figure}[!htp]
    \includegraphics[width=0.45\textwidth]{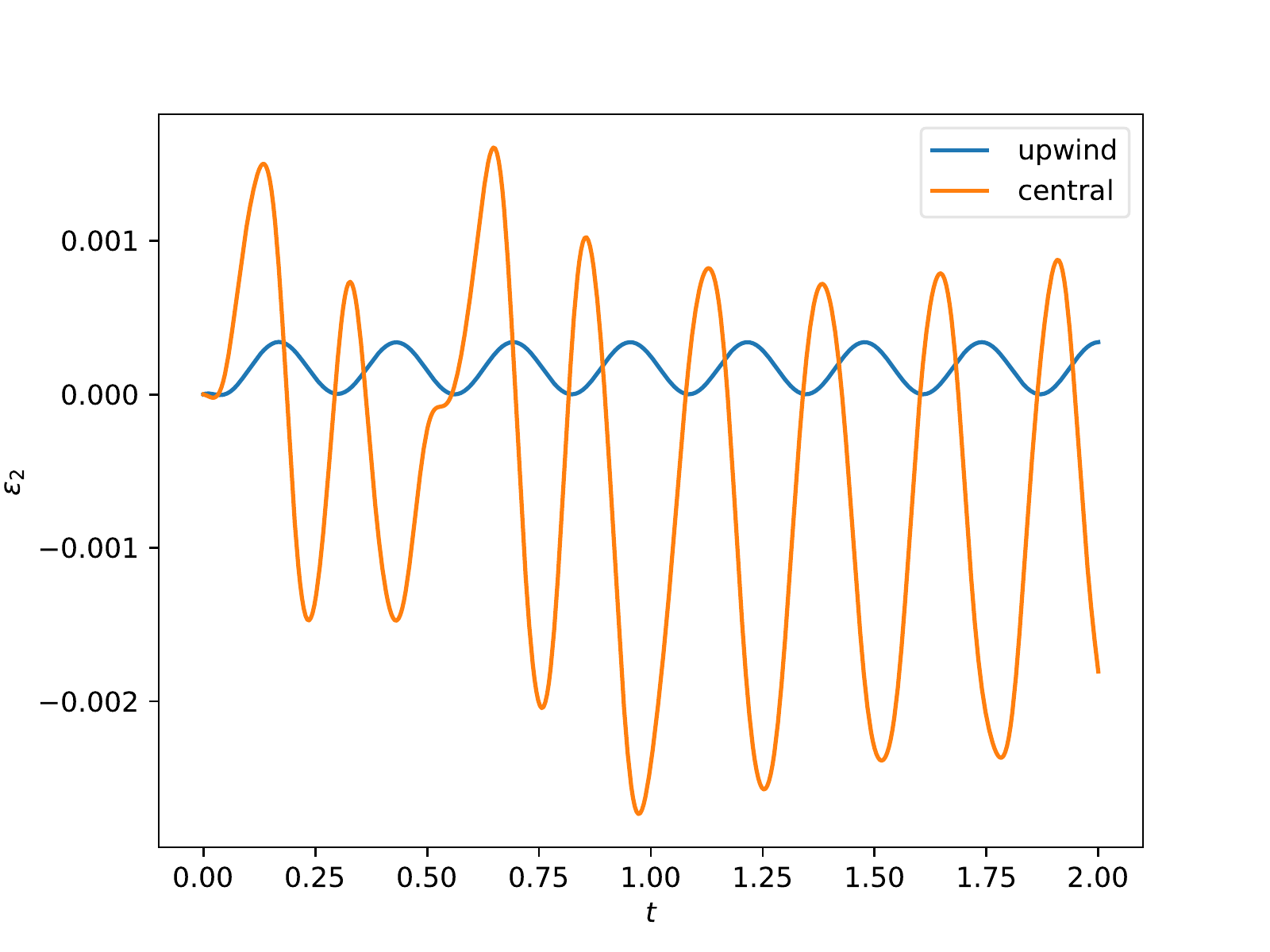}
    \includegraphics[width=0.45\textwidth]{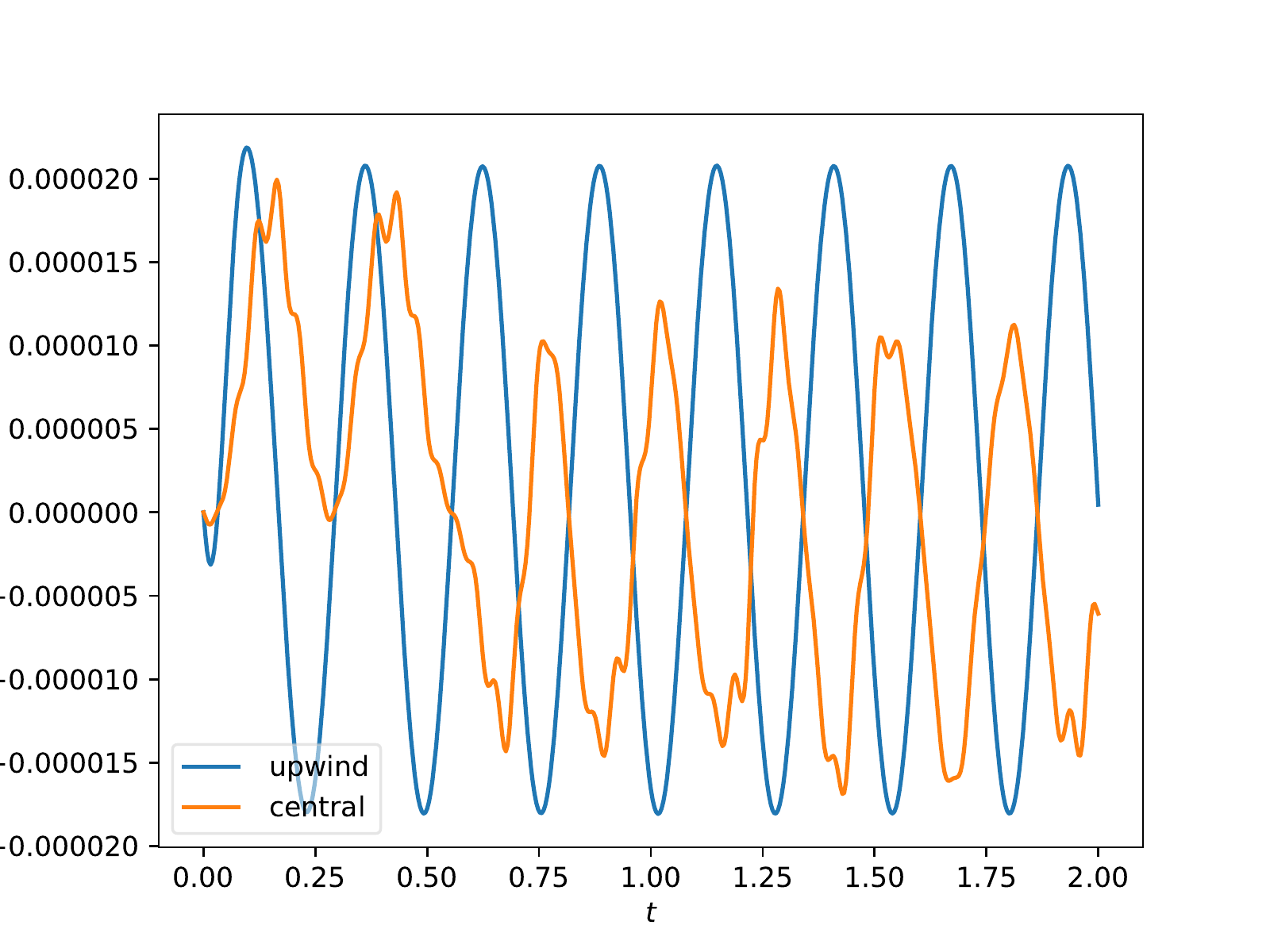}
  \caption{$\tilde{\epsilon}_2$ a function in time $K=30$ and $N=4, 5$. \\
 }
  \label{fig:epsilon_2}
\end{figure}
 All of these results are consistent with our investigation / predictions from before. 
 
 \subsection{Limitations and Counterexamples} \label{subsec:counter}
 We make a series of test calculations and most of the time
 the simulations support our predictions.  Nevertheless, there are several 
examples which 
question some of the  predictions $P2-P7$. We consider and discuss 
in the following several examples when  P2 or P5 are not longer correct.
At the end, we will give further a counterexample if the 
solution is not in the assumed space.

\subsubsection*{P2 is incorrect }
We start in our first example with the sine-testcase and use the 
FR scheme which is equivalent to the DG framework. 
In the error plot \ref{fig:error}, one  realizes that the upwind error
lies under the central error if a polynomial approximation 
of order three is used. In fact, we see this clearly in figure (\ref{fig:sind_test_cont} (a)).
Here, the central error lies above the upwind error, and also the 
asymptotic state is nearly the same. We assume that 
the noisy state is
\emph{periodic} with the central flux.
We may interpret this as using polynomial  order $3$ in our schemes is  too 
inaccurate for the approximation with the Gauß-Lobatto basis.
Then, applying an upwind flux adds too much dissipation into our calculation,
and this unpredictable behavior 
contradicts P2. 
In figure \ref{fig:sind_test_cont} (b), we get a similar error behavior as before 
if we decrease the number of elements\footnote{In 
\cite{gassner2011comparison}  the influence of the dispersion and dissipation 
errors of Gauß-Legendre and Gauß-Lobatto 
is investigate also in respect to the number of elements.}
$K$. 
With the higher jumps at the element interfaces, the upwind flux yields a more inexact 
numerical solution. We may conclude that we need an adequate number of elements 
to get the predicted results in the Gauß-Lobatto case.  
   \begin{figure}[!htp]
\centering 
  \begin{subfigure}[b]{0.4\textwidth}
    \includegraphics[width=\textwidth]{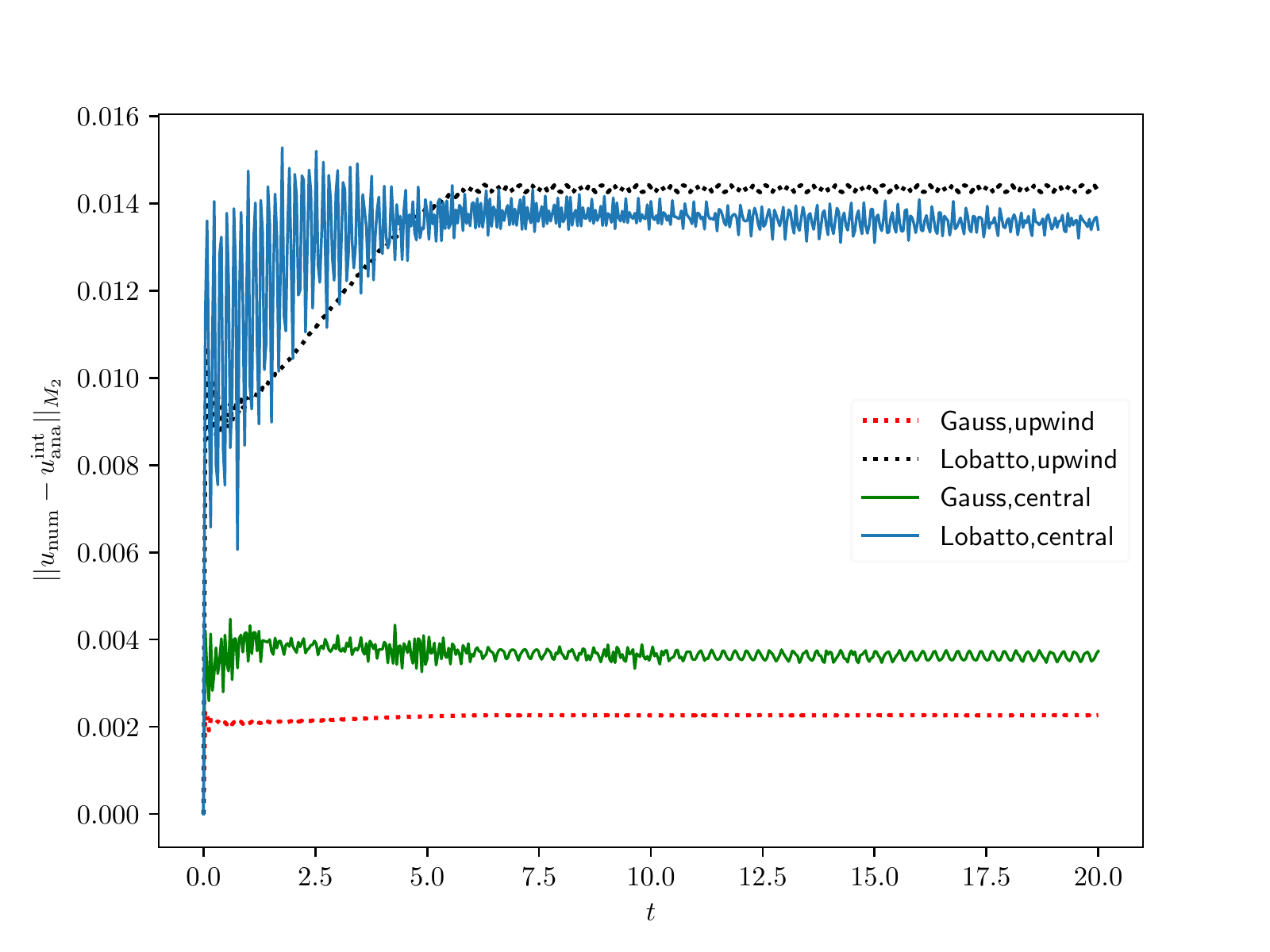}
    \caption{$N=3$, $K=50$, $t=20$}
  \end{subfigure}%
  ~
  \begin{subfigure}[b]{0.4\textwidth}
    \includegraphics[width=\textwidth]{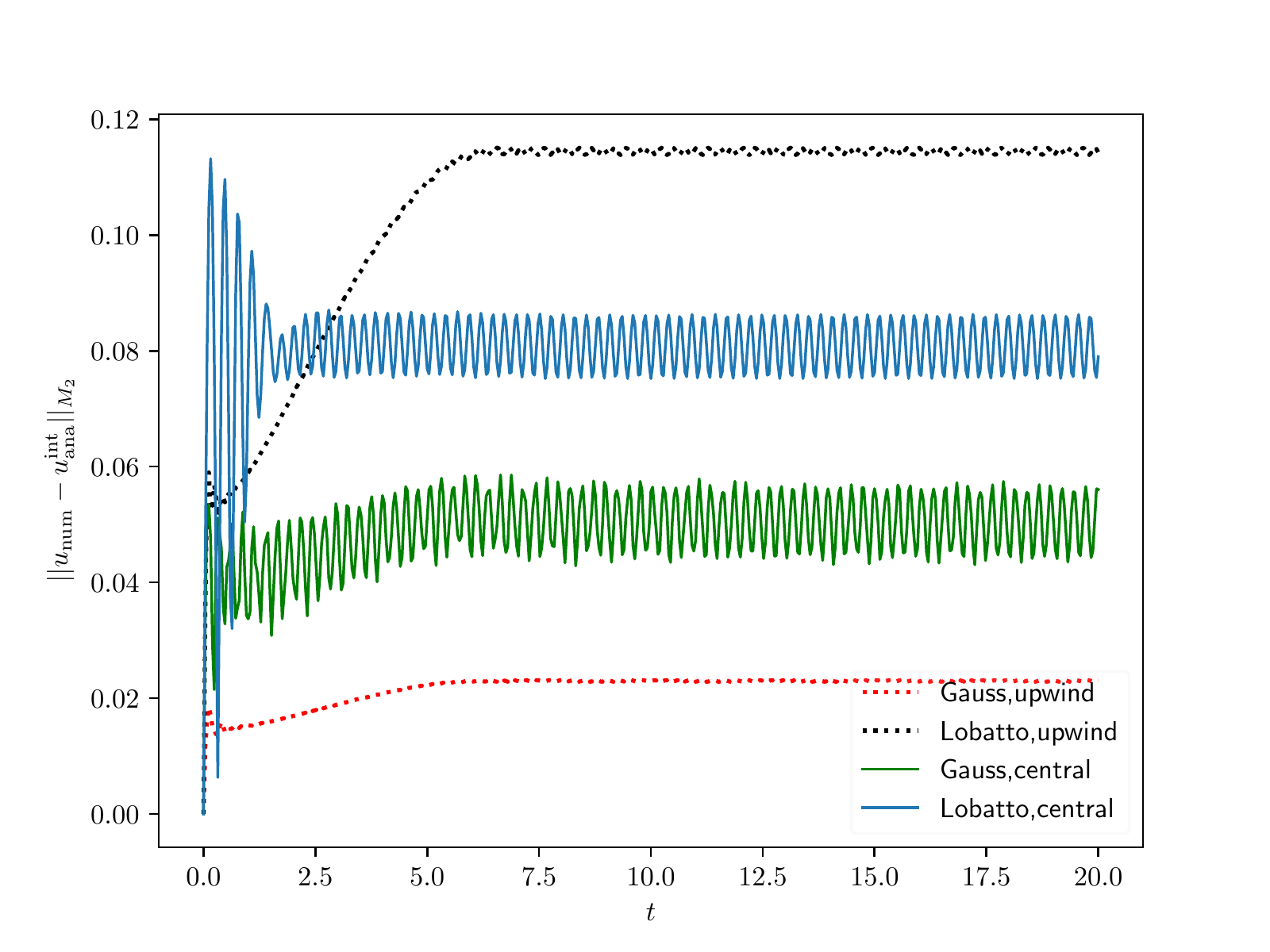}
    \caption{$N=4$, $K=20$, $t=20$}
  \end{subfigure}%
  \caption{Error as a function in time.}
  \label{fig:sind_test_cont}
\end{figure}
However, the numerical errors (upwind and central) with the Gauß-Legendre basis show the suspected 
behaviors from our results and one may interpret that as an advantage by using this basis type, 
but this is not the case. 
In figure \ref{fig:cos_test_cont} (a) we see the numerical errors of the cosine-case when using
polynomial order $3$ and 
$20$ elements. Here, the errors applying Gauß-Lobatto nodes behave in accordance to  (P2) and the 
errors with Gauß-Legendre basis do not.
We suppose that by using Gauß-Legendre basis and an upwind flux the jumps between the element interfaces
is too high and we get this effect.
If we again increase the number of elements, and so, the numbers of degrees of freedom, we realize a change 
in the error behaviors (compare \ref{fig:cos_test_cont} (a)-(c)). Nevertheless, 
the absolute error when using Gauß-Legendre nodes is comparatively low. 
   \begin{figure}[!htp]
\centering 
  \begin{subfigure}[b]{0.33\textwidth}
    \includegraphics[width=\textwidth]{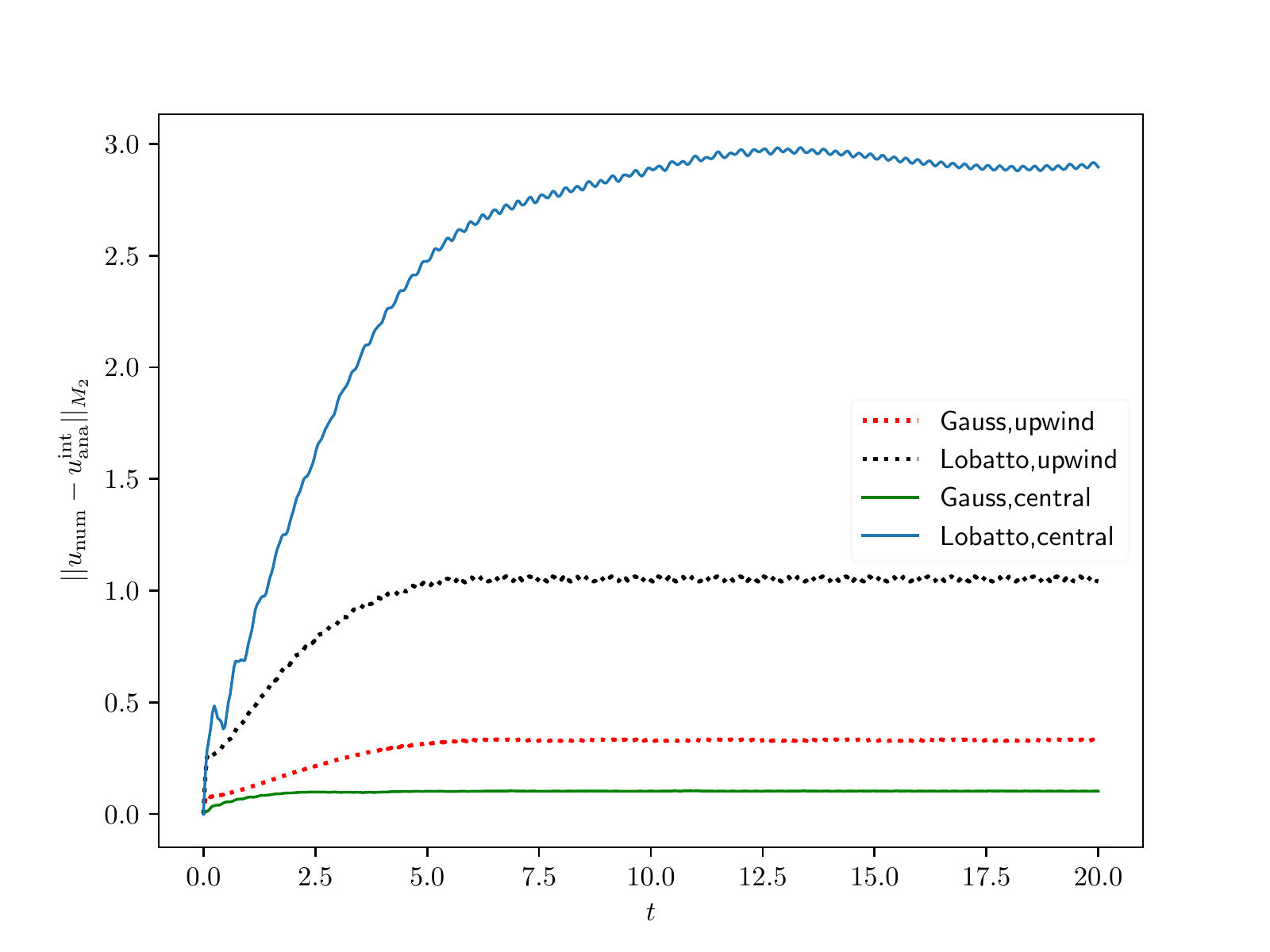}
    \caption{$N=3$, $K=20$, $t=20$}
  \end{subfigure}%
  ~
  \begin{subfigure}[b]{0.33\textwidth}
    \includegraphics[width=\textwidth]{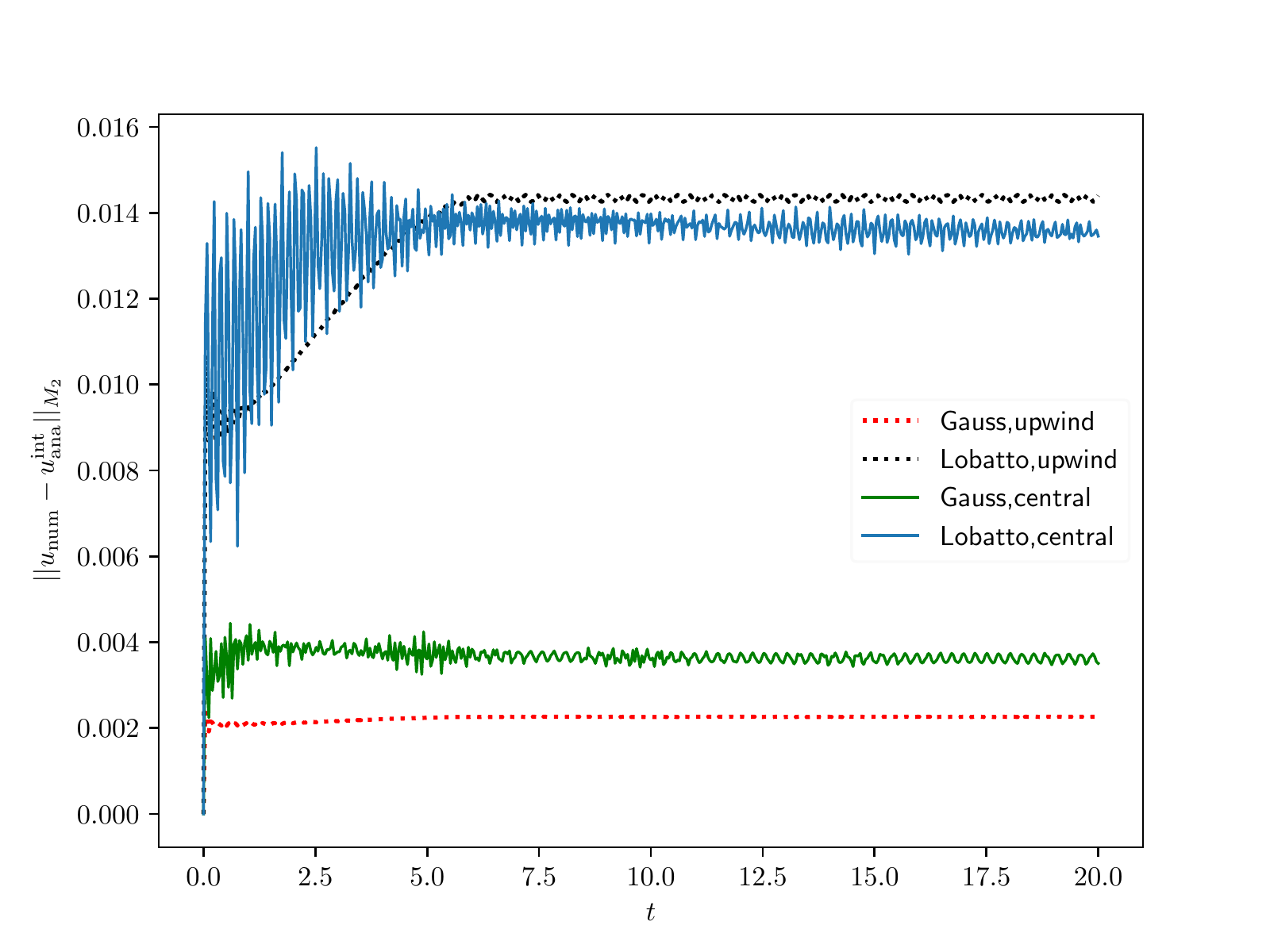}
    \caption{$N=3$, $K=50$, $t=20$}
  \end{subfigure}%
    ~
  \begin{subfigure}[b]{0.33\textwidth}
    \includegraphics[width=\textwidth]{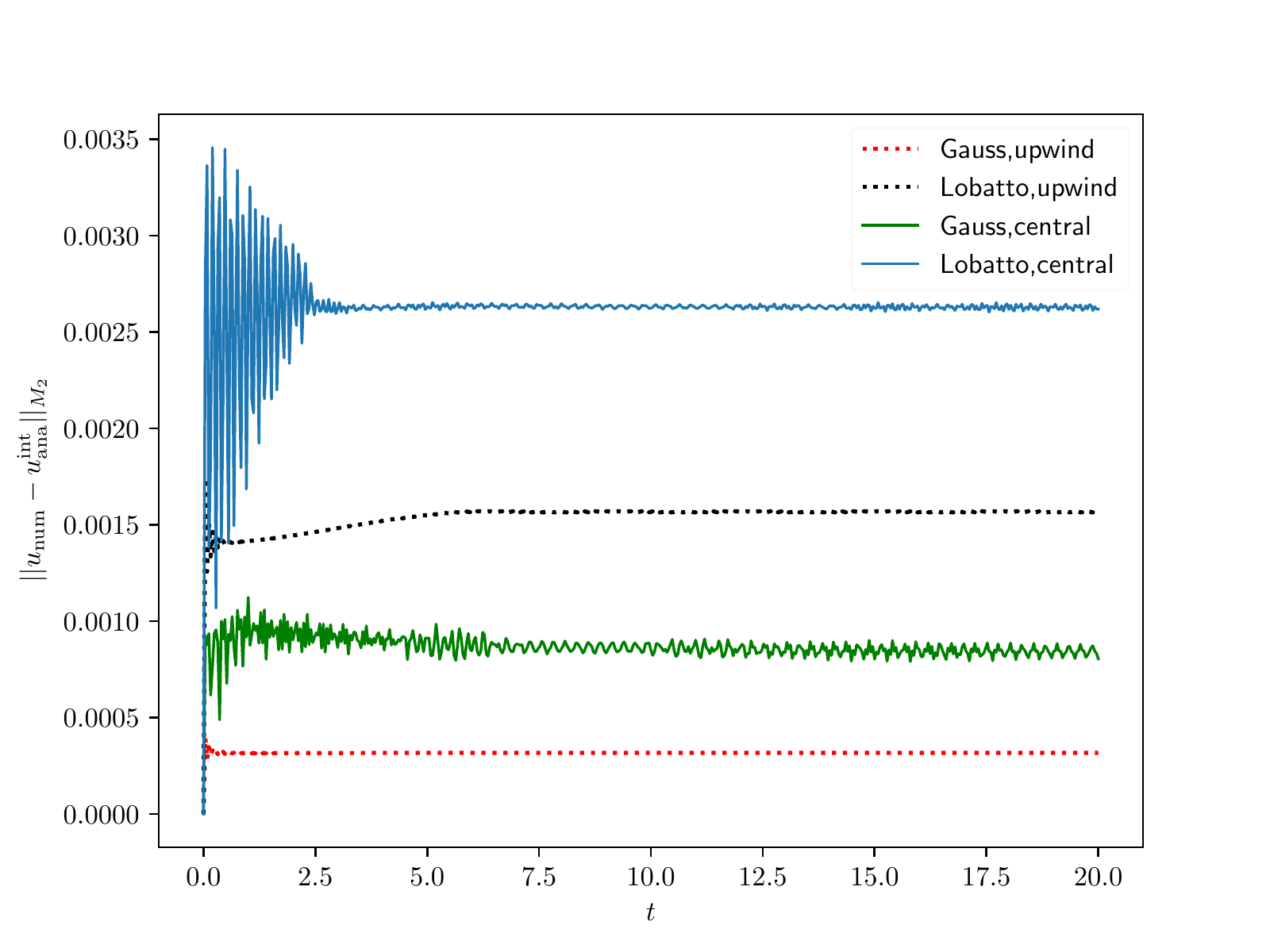}
    \caption{$N=3$, $K=80$, $t=20$}
  \end{subfigure}%
  \caption{Error  as a function in time.}
  \label{fig:cos_test_cont}
\end{figure}
 This limitation is noticed by focusing on the scheme with the  correction matrix
 $\mat{C} = \mat{M}[^{-1}] \mat{R}[^T] \mat{B}$.
 As it is already known \cite{vincent2011newclass, vincent2015extended, castonguay2012newclass},
 the most accurate 
 results are obtained if $\kappa=0$.
 Simultaneously, these results can also be  seen as an example of the even-odd phenomena for central fluxes. 
 It is observed in the literature \cite{chen2017entropy, gassner2013skew, hindenlang2019order,winters2015comparison}
 when the numerical flux function at the interface is  symmetric (with 
 the central flux), then the convergence order for an order $N$ scheme is $N+1$ if $N$ is odd, or only $N$
 if $N$ is even. 
 If the numerical flux adds dissipation, then this even-odd behavior is usually gone.
 Since the choice of the flux has influence of the convergence order, and it  has so on the errors. 

 \subsubsection*{P5 is not correct }
 Nevertheless, for other correction terms we have also some limitations for our prediction P5.
 Already in figure \ref{fig:other}, we realize that error using Gauß-Lobatto nodes and 
 an upwind flux is less 
 compared to the case applying Gauß-Legendre nodes, 
 since the amplitudes of the oscillations are not as high.
 Indeed, for FR schemes with $\kappa\neq0$ we have some 
 over and under corrections at the boundaries of every element. 
 If we increase the order or accuracy these amplitudes will decrease and, since the Gauß-Lobatto nodes 
 include the boundaries, it will yield to an error which is equal or less than applying Gauß-Legendre nodes.

 Simultaneously, we have to realize that using Gauß-Legendre leads, already in low order computations, to 
 comparatively good results (compare figures \ref{fig:other} and
 \ref{fig:Less}). 
     \begin{figure}[!htp]
\centering 
  \begin{subfigure}[b]{0.45\textwidth}
     \includegraphics[width=\textwidth]{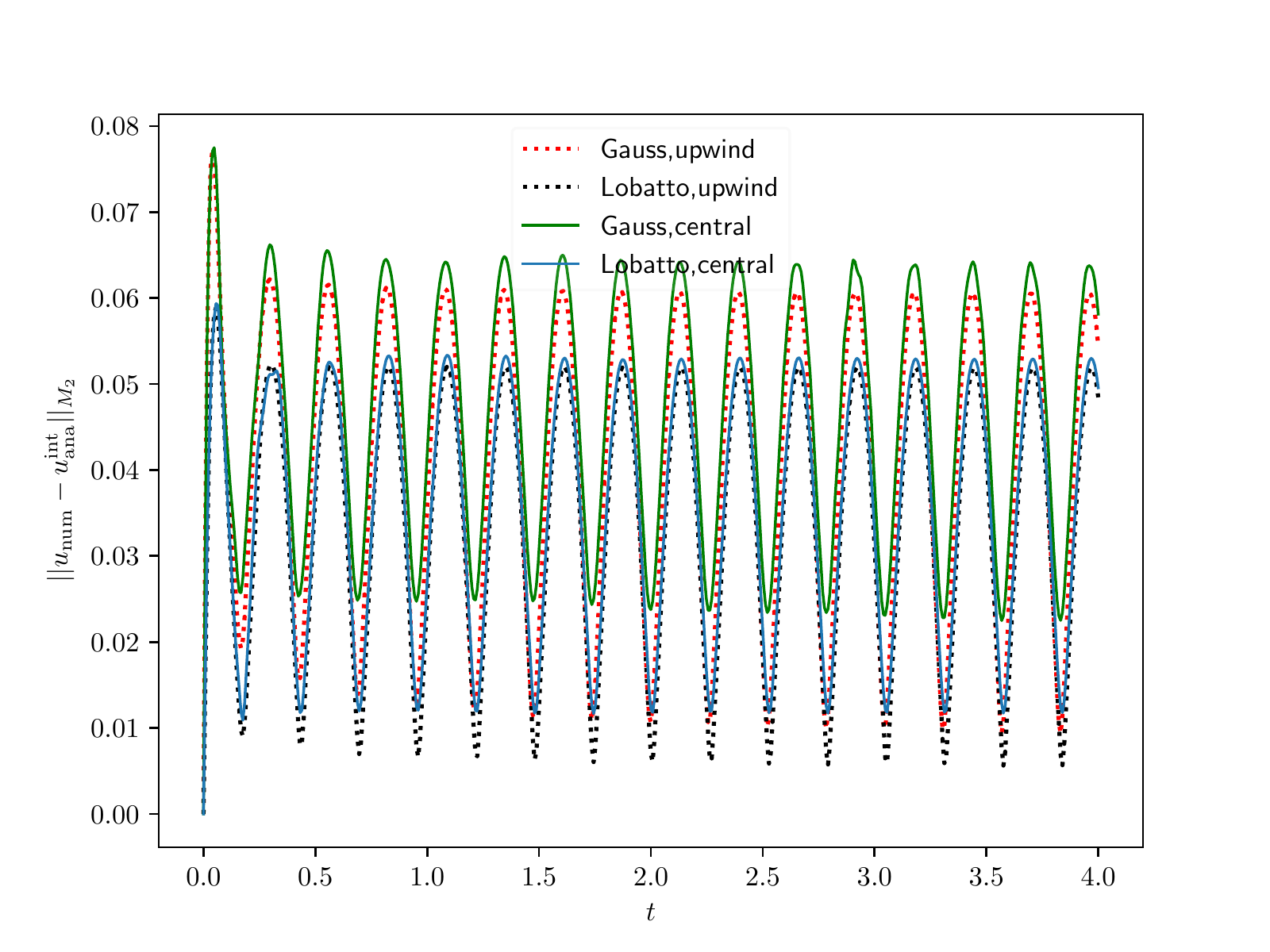}
    \caption{$N=6$, $K=30$, $t=4$, SD}
  \end{subfigure}%
  ~
  \begin{subfigure}[b]{0.45\textwidth}
    \includegraphics[width=\textwidth]{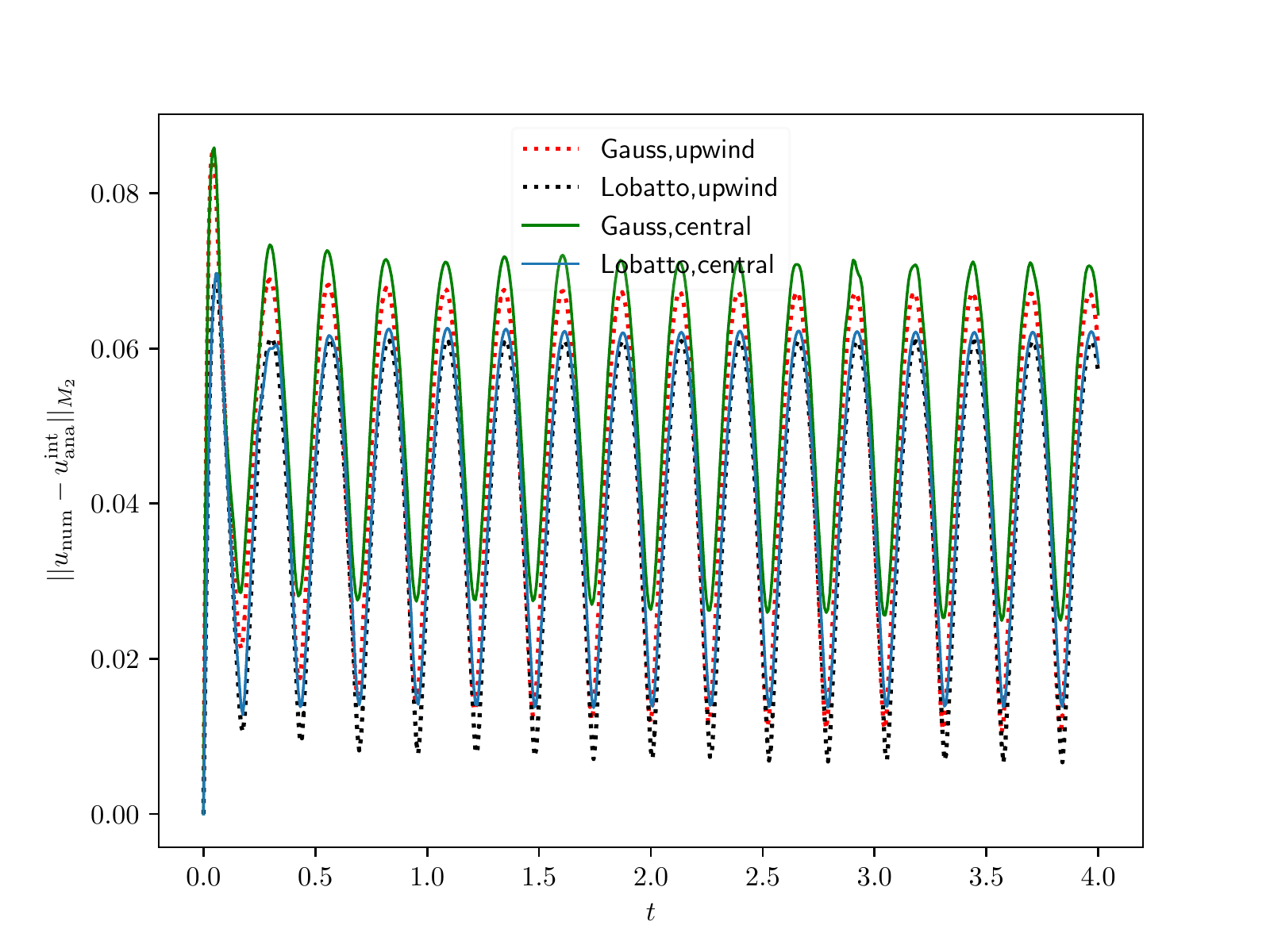}
    \caption{$N=6$, $K=30$, $t=4$, Huynh}
  \end{subfigure}%

  \caption{Error  as a function in time.}
  \label{fig:Less}
\end{figure}

 \subsubsection*{Counterexample }
 In section \ref{sec:model} we mentioned an example where the norm of 
 solution  $||u||_{H_{\kappa,N}^m}$ is not uniformly
 bounded in time. 
 We select our initial and boundary conditions in such way
 that we get as the solution
 $u(t,x)=(x-t)^8$. Using $\kappa=0$ the  figure
 \ref{fig:counter}(a) shows the unbounded increase of the errors. 
 However, we have to mention that the relative errors remain bounded \ref{fig:counter}(b).

    \begin{figure}[!htp]
\centering 
  \begin{subfigure}[b]{0.45\textwidth}
    \includegraphics[width=\textwidth]{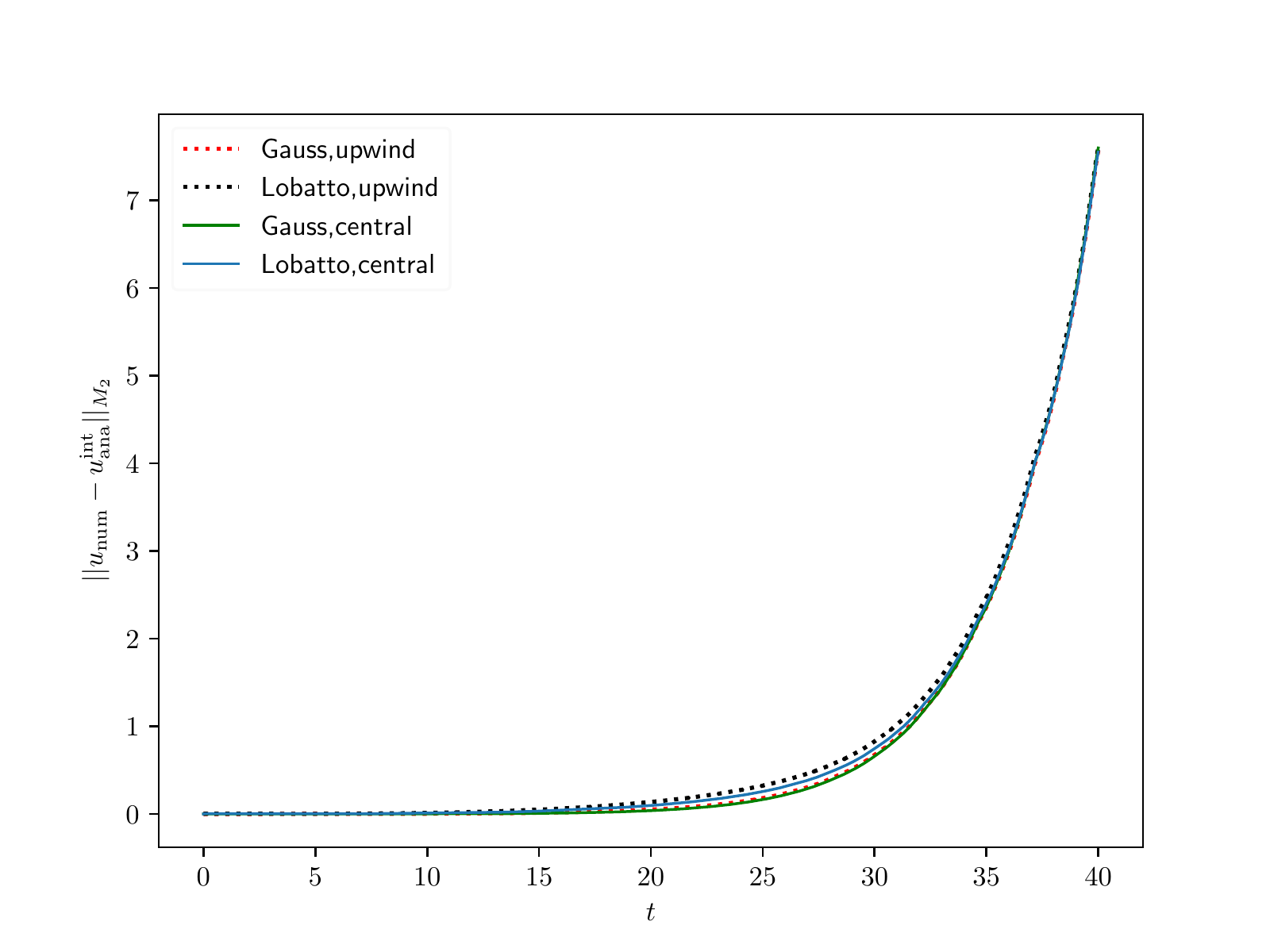}
    \caption{$N=4$, $K=50$, $t=40$}
  \end{subfigure}%
    ~
  \begin{subfigure}[b]{0.4\textwidth}
    \includegraphics[width=\textwidth]{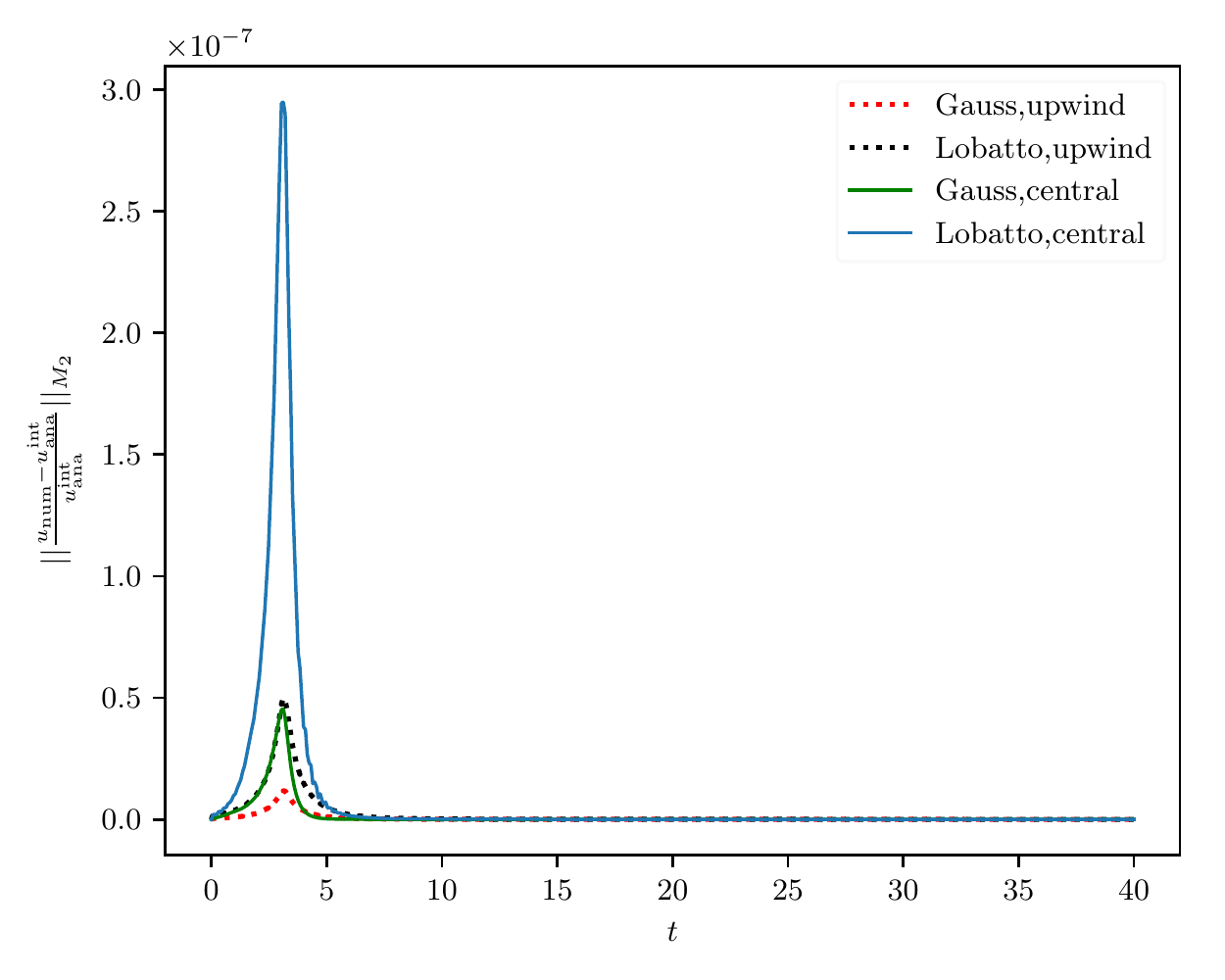}
    \caption{$N=4$, $K=50$, $t=40$}
  \end{subfigure}%
  \caption{Left: Errors as  functions in time. Right: Relative errors}
  \label{fig:counter}
\end{figure}

\section{Summary and Conclusion}
\label{Sec:Summary}

In this paper, we transfer the results about the bounded error growth 
from the discontinuous Galerkin spectral element method \cite{kopriva2017error} 
to the more general framework of SBP-FR methods.
Furthermore, we extend the investigation by including the Gauß-Legendre basis, 
where \cite{kopriva2017error} considers only the Gauß-Lobatto basis. 
Indeed, for both bases (Gauß-Lobatto / Gauß-Legendre), the numerical flux used
at element boundaries affects the error growth behavior. If an adequate number
of elements is used, the upwind flux leads to better results.
The asymptotic values are smaller and are reached in a shorter time period. 
At once, also the selection of basis has a big influence and in our opinion is
even more important. 
Using Gauß-Legendre basis, the error reaches the asymptotic value faster 
and to a lower value than by using Gauß-Lobatto basis.
Also, the impact of the different numerical fluxes (central / upwind) 
when applying Gauß-Legendre basis is less important than in the Gauß-Lobatto case,
especially using a low order polynomial approximation. 
These effects decrease when the order of polynomial approximation is increased and/or using more 
elements (which also increase the resolution). \\
The investigation implies that the usage of Gauß-Legendre basis has some advantages compared to 
Gauß-Lobatto and should be preferred. 
However, there are several points which we have to mention yet.\\
We investigate a trivial model problem \eqref{eq:Model} where the flux function is simple; 
$f(u)=u(x,t)$. Already by using the more complicated flux $f(u)=a(x)u(x,t)$ several 
problems arise  in the discretization by using Gauß-Legendre nodes, see \cite{manzanero2017insights} for details.
The reason is that Gauß-Legendre points do not include the boundary points in one element, and we get
some aliasing effect if we are not careful in the discretization.  
In \cite{ranocha2018generalised}, the author proves a way to solve these issues by
applying further correction terms to approximate the boundary terms correctly.
Recently,
the authors have investigated the long time error behavior of a DG method in this context in  \cite{offner2019error}.
For non-linear flux functions stability problems rise automatically. The aliasing effect 
is quite stronger and to remedy these issues, further correction terms are needed \cite{ranocha2016summation}.
By including the boundary points, these correction terms are simpler and better understood. 
\cite{ranocha2017shallow}  provides the correction terms for the shallow water equation using 
Gauß-Legendre nodes  and 
recently, the concept of \emph{decoupled SBP} operators introduced by Chan \cite{chan2018discretely} 
is used to build those correction terms for the Euler equations.
The numerical study in \cite{chan2018efficient}
demonstrates also some advantages of applying Gauß-Legendre nodes which supports our theoretical analysis here.
In our calculations, the time integration analysis was neglected, but in practice
it is also an important issue. In \cite{gassner2011comparison}, the authors already investigate 
the time-step restriction 
by using Gauß-Lobatto or
Gauß-Legendre nodes in the DGSEM and find out that Gauß-Lobatto nodes have favored properties.\\
The above mentioned issues are not unimportant. However, 
due to our analysis and the numerical results (also in \cite{chan2018efficient})
the usage of Gauß-Legendre basis should be taken into account. In our tests the asymptotic error values 
are  reached  faster and to a smaller amount.  Nevertheless, further studies
are  necessary. First, one has to analyze the impact of $\tilde{\epsilon}_2$ 
not only numerically but also analytically.
Secondly, we must study 
what happens with the approximation error if the initial conditions have jumps or even more complex
flux functions are considered. 


\section{Appendix}\label{sec:appendix}

\subsubsection*{Connection to DG}

As it is described several times, there is a close connection between the DG 
and the FR framework. 
Therefore, we repeat the main aspect and present one example for a better understanding.
FR schemes use in their discretizations  
of  \eqref{eq:scalar_CL}
no weak /variational or integral form. A differential form \eqref{eq:SBP CPR} is applied.
The main idea of the FR schemes is that the numerical fluxes at the boundaries are corrected by 
correction functions in such manner that basic properties (e.g. conservation) hold.
In \cite{vincent2011newclass}, the authors develop a class of energy stable FR schemes
depending on a single scalar parameter. The correction functions 
are given for the left and right boundary in one element by the following formula:
{\small
\begin{equation}\label{eq:correction_functions}
 c_{LB}(\xi)= \frac{(-1)^p}{2} \l[ L_p(\xi) -\l(\frac{\lambda_p L_{p-1}(\xi)+L_{p+1}(\xi)}{1+\lambda_p} \r) \r], \quad  
 c_{RB}(\xi)= \frac{1}{2} \l[ L_p(\xi) +\l(\frac{\lambda_p L_{p-1}(\xi)+L_{p+1}(\xi)}{1+\lambda_p} \r) \r],
\end{equation}}
where $L_p$ is the $p$-Legendre polynomial and $\lambda_p=\frac{\kappa(p+1) 2^{2p} (p!)^4 }{((2p)! p!)^2} $ 
is a term with the free parameter $\kappa$. 
The translation about these correction functions and our notation can be found 
in \cite{ranocha2016summation}. \\
As mentioned before, we get into the DG framework by selecting $\kappa=0$.
Then, the corrections functions are the right and left Radau polynomials
and the application of these polynomials as correction functions is essential. 
We present the following example
from \cite[pages 23-25]{huynh2007flux}. We strongly recommend also 
the review paper \cite{huynh2014high} where this connection is also pointed out.

\begin{ex}\label{Ex:FR}

We are considering a scalar conservation law{\small
 \begin{equation}\label{eq:advection}
  \partial_t u + \partial_x f(u)=0.
 \end{equation}}
 A DG formulation for the problem on the standard interval $I=[-1,1]$   
 is given by {\small
\begin{equation}\label{DGEx}
 \int_I   \frac{\Delta x_k}{2} (\partial_t U) \phi(\xi)  \d \xi
+  \phi(1) \fnum_{up}(1)-\phi(-1) 
\fnum_{up}(-1)
-\int_I  \partial_{\xi} \phi(\xi) \d \xi=0,
\end{equation}}
where $U$, $F$ are polynomials that approximate $u$ and $f$.
$\fnum_{up}$ is the upwind flux and $\phi$ is the test function (polynomials of
degree $N$).
Instead of solving \eqref{DGEx}, we want to eliminate the 
test function $\phi$. Therefore, we apply again integration by parts 
and reformulate \eqref{DGEx} to {\small
\begin{equation}\label{eq:DGEx2}
 \int_I   \frac{\Delta x_k}{2}  (\partial_t U) \phi(\xi) \d \xi
+ \phi(1) [\fnum_{up}(1)-F(1)]-\phi(-1) 
[\fnum_{up}(-1)-F(-1)]
+\int_I (\partial_{\xi} F) \phi(\xi)\d \xi=0.
\end{equation}}
The right and left Radau polynomials  of degree $N+1$ 
have the property that for any polynomial $\phi$ of degree $N$ 
or less the following equations  
%
{\small
\begin{equation}\label{eq:radauproperty}
 -\phi(-1)=\int_{-1}^{1} (\partial_{\xi}c_{LB}(\xi))\phi(\xi) \d \xi \text{ and } 
  \phi(1)=\int_{-1}^{1} (\partial_{\xi} c_{RB}(\xi))\phi(\xi) \d \xi 
\end{equation}}
are fulfilled.
With this property \eqref{eq:radauproperty}
we are able to factor out $\phi$ in \eqref{eq:DGEx2} and obtain
{\small
\begin{equation*}
 \int_I \l(  \frac{\Delta x_k}{2}  (\partial_t U)+ 
 (\partial_{\xi} \hat{F}) \r) \phi(\xi) \d \xi=0
\end{equation*}}
with
{\small
\begin{equation*}
  \hat{F}(\xi)= F(\xi)+ [\fnum_{up}(1)-F(1)]c_{RB}(\xi)+
[\fnum_{up}(-1)-F(-1)]c_{LB}(\xi).
\end{equation*}}
Switching to the global coordinate,
{\small
\begin{equation}\label{eq:DG_Ex_final}
 \int_{x^{k-1}}^{x^k} \l(  (\partial_t U)+ 
 (\partial_{x} \hat{F}) \r) \phi(x) \d x=0.
\end{equation}}
Since the equation \eqref{eq:DG_Ex_final} holds for any polynomial $\phi$ of degree $N$,
it is equivalent to 
{\small
\begin{equation}\label{eq:Flux_diff}
  \partial_t U+ 
 \partial_{x}\hat{F} =0,
\end{equation}}
which is nothing else than the flux reconstruction scheme
with the Radau polynomials as correction functions.
Finally, we showed that the DG scheme is equivalent to this FR method.
\end{ex}

This connection  is also pointed out in the review article \cite{huynh2014high}.
Furthermore, the relation to the DGSEM can also be seen by comparing
the work \cite{gassner2013skew} and our introduction  in section \ref{sec2:CPR}.
Here, also the used notations to describe the methods are quite similar.

\subsubsection*{Stability of FR schemes}
Finally, we like to mention that Jameson utilizes
in his investigation in \cite{jameson2010proof}
a kind of broken Sobolev norm.\\
Let $m\in \N_0$. The norm of the Sobolev space $H^m((-1,1))$
is given by 
{\small
\begin{equation*}
 ||u||_{H^{m}((-1,1))}:= \l(\sum\limits_{j=0}^m ||u^{(j)}||_{\L^2((-1,1))}^2\r)^\frac{1}{2},
\end{equation*}}
where the derivatives $u^{(j)}$ are taken in a weak sense.
In \cite{jameson2010proof},
the norm 
{\small
\begin{equation}\label{eq_broken_Sobloev_norm}
||u||^2_{H_{\kappa,N}((-1,1))}:=\int_{-1}^1  u^2+\kappa (u^{(N)})^2 \d \xi,
\end{equation}}
is used. Here, $N$ represents the order of accuracy.
The solution space (a polynomial space) is equipped with the norm 
\eqref{eq_broken_Sobloev_norm} and linear stability is studied in this context.
This norm is also used in \cite{vincent2011newclass} where the one-parameter family of FR schemes is 
developed, the correction functions \eqref{eq:correction_functions}
are defined, and coefficients $\kappa$ are determined to embed the known schemes into their setting.\\
As it can be found in \cite{vincent2011newclass}, 
the nodal DG framework is obtained by selecting $\kappa\equiv0$.
For spectral difference and Huynh scheme, we have
{\small
\begin{equation*}
 \kappa_{SD}(N)=\frac{2N}{(2N+1)(N+1)(a_NN!)^2} \qquad \text{ and } 
 \kappa_{Hu}(N)=\frac{2(N+1)}{(2N+1)N(a_NN!)^2}.
\end{equation*}
}
These formulas have been used to calculate the values in table \ref{ta:correction_terms}.\\
Instead of working with this continuous broken Sobolev norm 
\eqref{eq_broken_Sobloev_norm} from \cite{jameson2010proof},
we apply the discrete counterpart and use this in our investigation.
More details can be found in \cite{ranocha2016summation}.

  {\footnotesize{
\bibliographystyle{abbrv}
\bibliography{literature}
}}

\end{document}